\documentclass[11pt,letterpaper]{amsart}
\usepackage{amsmath,amssymb,amsthm,
color, euscript, enumerate}
\usepackage{dsfont}                                                                          
\usepackage{longtable}
\usepackage{tikz}
\usepackage{caption}
\usepackage{multirow}
\usepackage{graphicx}

\newcommand{\RR}{{\mathbb R}}

\newcommand{\CC}{{\mathbb C}}
\newcommand{\be}{\begin{equation}}
\newcommand{\ee}{\end{equation}}
\newcommand{\ba}{\begin{array}}
\newcommand{\ea}{\end{array}}
\newcommand{\bea}{\begin{eqnarray}}
\newcommand{\eea}{\end{eqnarray}}

\usepackage{graphicx}
\usepackage{color}
\newcommand{\corr}[1]{{\color{black}{#1}}}
\usepackage{transparent}

\newtheorem{theorem}{Theorem} [section]

\newtheorem{lemma}{Lemma}[section]
\newtheorem{proposition}{Proposition}[section]
\newtheorem{definition}{Definition} [section]

\DeclareMathOperator{\sgn}{sgn}
\DeclareMathOperator{\res}{res}

\numberwithin{equation}{section}
\numberwithin{table}{section}

\allowdisplaybreaks

\begin{document}
\begin{center}
{\large   \bf STABILITY OF KADOMTSEV-PETVIASHVILI MULTI-LINE SOLITONS }
 
\vskip 15pt

{\large  Derchyi Wu}

\vskip 5pt

{ Institute of Mathematics, Academia Sinica, 
Taipei, Taiwan}

e-mail: {\tt mawudc@gate.sinica.edu.tw}

{\today}
\end{center}
\vskip 10pt

\vskip 10pt
{\bf ABSTRACT}

\begin{enumerate}
\item[]{\small We prove the long-standing inverse scattering theory (IST) of  perturbed Kadomtsev Petviashvili multi-line solitons. Our work is the first rigorous IST of a multi-dimensional integrable system when both continuous and discrete scattering data are present, and the support of continuous scattering data does not degenerate into contours in the complex plane. As an application,   an $L^\infty$-stability theorem of the Kadomtsev Petviashvili multi-line solitons is justified.}

\end{enumerate} 

\tableofcontents
\section{Introduction}\label{S:goal}
The Korteweg-de Vries (KdV) equation
\be\label{E:KdV}
-4u_{x_3}+u_{x_1x_1x_1}+6uu_{x_1}=0  
\ee
 is an asymptotic model for the propagation of one-dimensional small amplitude  long wave surface waves. Via asymptotic characterization, variational techniques, and nonlinear analysis, it is shown that the $L^\infty$-stability of the $1$-solitons 
\be\label{E:KdV-soliton}
u_c(x_1,x_3)=  2c^2\textrm{sech}^2 \left(c(x_1+c^2x_3)\right),
\ee
or $n$-solitons (setting $x_2\equiv 0$ in \eqref{E:line-tau}), as   solutions  of the KdV equation, is not true but   orbital stability in various Sobolev spaces holds \cite{Ben72,Bo75,Wein86,MS93,MMT02,MV03,MM05,AMV13,KV22}. Integrability is hardly or less exploited in these PDE approaches. Hence methods work for non integrable systems 
  with low regularity as well but detailed description of the evolution profile  is lost (cf \cite{Wein86}). 
 
To study the soliton stability under weak transversal perturbations of the KdV equation, Kadomtsev and Petviashvili derived the two-dimensional models (KP equations),
\be\label{E:KPII-intro}
(-4u_{x_3}+u_{x_1x_1x_1}+6uu_{x_1})_{x_1}\pm 3u_{{x_2}{x_2}}=0.    
\ee
and conjectured that solitons \eqref{E:KdV-soliton} are unstable for the KPI equation (corresponding to $-$ where surface tension is present), and stable for the KPII equation (corresponding to $+$ where surface tension is absent) \cite{KP70}. The phenomena were analysed formally using the IST through the Gelfand-Levitan-Marchenko equation \cite{Z75,Bur84}. Rigorous theory is obtained through PDE approaches: precisely, for the KPI equation, global well posedness of data being $L^2(\mathbb R^2)$-, $L^2(\mathbb R\times T)$-perturbations of \eqref{E:KdV-soliton}  are proved  \cite{IK05,MST07},   orbital stability and instability  of \eqref{E:KdV-soliton} in $L^2(\mathbb R\times T)$-perturbations are justified  \cite{RT09,RT12}; for the KPII equation, global well posedness of data being $H^s(\mathbb R\times T)$-, $H^s(\mathbb R^2)$-perturbations of \eqref{E:KdV-soliton} have been solved  \cite{MST11}, orbital stability in $L^2(\mathbb R\times T)$-perturbations is proved   \cite{MT12},  and orbital stability, instability theories in $L^2(\mathbb R^2)$-perturbations are derived  \cite{M15,M19}. We refer to \cite{KS21} and the reference therein for more details on this topic.

The goal in this paper is to prove the following $L^\infty$-stability  of KPII multi-line  solitons:  
\begin{theorem}     \label{T:intro-stability}
If $u_0(x_1,x_2)=u_s(x_1,x_2,0)+v_0(x_1,x_2)$, $u_s(x)=u_s(x_1,x_2,x_3)$ is a $  {\mathrm{Gr}(N,M)_{> 0}}$ KP soliton, $  \sum_{|l|\le {d+8}} |{ (1+|x_1|+|x_2|)}\partial_x^lv _0|_{L^1\cap L^\infty} \ll 1$, and $ {d\ge 0}$,  then there exists $u :\mathbb R\times\mathbb R\times \mathbb R^+\to \mathbb R$ such that
 \begin{gather}
(-4u_{x_3}+u_{x_1x_1x_1}+6uu_{x_1})_{x_1}+3u_{{x_2}{x_2}}=0, \ \
u(x_1,x_2,0)=u_0 (x_1,x_2), \label{E:theorem-1.1}\\
\sum_{0\le l_1+2l_2+3l_3\le d+4}|\partial^l_x\left[u(x )-u_s(x )\right]|_{L^\infty }\le C \sum_{|l|\le {d+8}} |{ (1+|x_1|+|x_2|)}\partial_x^lv _0|_{L^1\cap L^\infty}.\label{E:theorem-1.1-1}
\end{gather} 
\end{theorem} 
Here and throughout the paper, $\partial_x^l\equiv \partial_{x_1}^{l_1}\partial_{x_2}^{l_2}\partial_{x_3}^{l_3}$, $  l_j$ are non negative integers, $|l|=l_1+l_2+l_3$, 
 $C$ denotes a uniform constant  which is independent of $x$  and the spectal variable $\lambda$.  To introduce the notion of {\bf $ \mathbf{\textrm{Gr}(N,M)_{> 0}} $ KP solitons}, note that one major breakthrough in the KPII theory was given by Sato. He realized that solutions can be written in terms of points of an infinite-dimensional Grassmannian \cite{S81,SS82,S89a,S89b}. In particular,  a real finite dimensional version of the Sato theory concerns      $   {\mathrm{Gr}(N,M)_{\ge 0}}$ KP  solitons  which are regular in the entire $x_1x_2$-plane with peaks localized and non decaying along certain line segments and rays for each fixed time $x_3$. They can be constructed by
 
\setlength{\unitlength}{1cm}
\begin{picture}(10,3.7)
\put(2,0.6){\includegraphics[trim=0mm 0mm 0mm 0mm,clip,height=3.0cm]{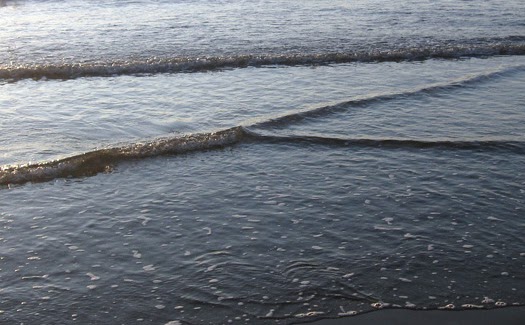}}
\put(8,0.6){\includegraphics[trim=0mm 0mm 0mm 0mm,clip,height=3.0cm]{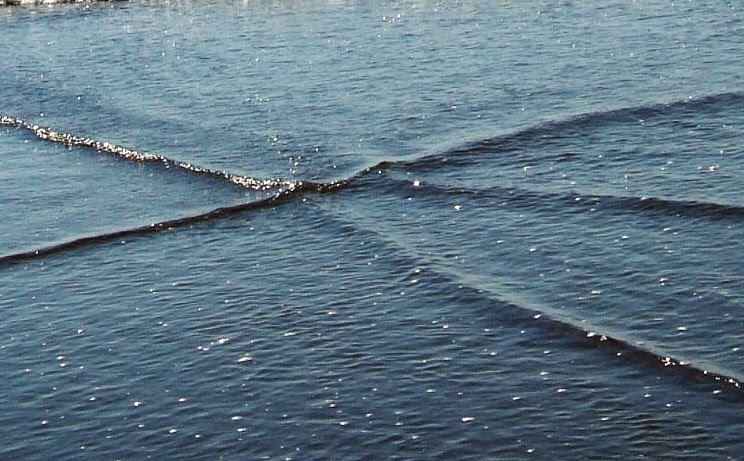}}
\put(2.5,0){\parbox{10cm}{\scriptsize Photos of $\mathrm{Gr}(1,3)_{>0}$ (left) and $\mathrm{Gr}(2,5)_{\ge 0}$(right) taken in Nuevo Vallarta, Mexico by Mark J. Ablowitz}}
\end{picture}

\be\label{E:line-tau}
u_s(x)= 2\partial^2_{x_1}\ln\tau(x) 
\ee   \cite{BC07,BK03,KW13} where the $\tau$-function is  the Wronskian determinant
\begin{align}
\tau(x)=&\left|
\left(
\begin{array}{cccc}
a_{11} &a_{12} & \cdots & a_{1M}\\
\vdots & \vdots &\ddots &\vdots\\
a_{N1} &a_{N2} & \cdots & a_{NM}
\end{array}
\right)
\left(
\begin{array}{ccc}
E_{1} & \cdots & \kappa_1^{N-1}E_1\\
E_{2} & \cdots & \kappa_2^{N-1}E_2\\
\vdots & \ddots &\vdots\\
E_{M} & \cdots & \kappa_M^{N-1}E_M\\
\end{array}
\right)
\right|\label{E:line-grassmannian} \\
=&\sum_{1\le j_1< \cdots< j_N\le M}\Delta_{j_1,\cdots,j_N}(A)E_{j_1,\cdots,j_N}(x),\nonumber
\end{align}
with $\kappa_1<\cdots<\kappa_M$,  
 $A=(a_{ij})\in {\mathrm{Gr}(N, M)_{\ge 0}}$ (full rank $N\times M$ real matrices with non negative minors),   $E_j(x)=\exp\theta_j(x)=\exp( \kappa_j x_1+\kappa_j^2 x_2+\kappa_j^3 x_3)$, $\Delta_{j_1,\cdots,j_N}(A)$  the $N\times N$ minor of the matrix $A$ whose columns are labelled by the   index set $J=\{j_1< \cdots< j_N\}\subset\{1,\cdots,M\}$, and $
 E_{j_1,\cdots,j_N}(x)=\Pi_{l<m}(\kappa_{j_m}- \kappa_{j_l})\exp ( \sum_{n=1}^N\theta_{j_n}(x)  )$.   Important progress in combinatoric properties, wave interaction,   resonant theories,   asymptotics characterization, and classification theories of $   {\mathrm{Gr}(N,M)_{\ge 0}}$ KP solitons have been developed  \cite{CK09,K17,K18}. Finally, $ \textrm{Gr}(N,M)_{> 0}$ KP solitons for which all minors $\Delta_{j_1,\cdots,j_N}$ are positive,  {namely, fulfill the totally positive (TP) condition}, form a dense subset of ${ \mathrm{Gr}(N,M)_{\ge 0}}$ KP solitons.

Notice that the solitary waves \eqref{E:KdV-soliton} are the  simplest $ {\mathrm{Gr}(1,2)_{> 0}}$ KP  solitons 
\be\label{E:line-tau-oblique}
 u_s(x) 
=  \frac{(\kappa_1-\kappa_2)^2}2\textrm{sech}^2\frac{\theta_1(x)-\theta_2(x)-\ln a}2 
\ee by setting $\kappa_1=-\kappa_2=c$ and $a=1$. 
   But the famous counter example, given by 
\[
\begin{gathered}
 u_j(x_1, 0)=2c_j^2\mbox{sech}^2(c_j x_1),\quad  c_1\ne c_2,\ \ |c_1-c_2|\ll 1,\\
|u_1(x_1,\frac{1}{c_2-c_1})-u_2(x_1,\frac{1}{c_2-c_1})|\sim \mathcal O(c_1^2), 
\end{gathered}\]
   which shows $L^\infty$-stability for the KdV $1$-soliton cannot hold and a phase shift in the orbital stability  is necessary, does not fulfill the smallness assumption of Theorem \ref{T:intro-stability} and is no longer a counterexample for $L^\infty$-stability for $ \textrm{Gr}(N,M)_{> 0}$ KP solitons. Precisely, in contrast to the KdV equation,   Theorem \ref{T:intro-stability} implies  velocities, amplitude, and phases of the leading term of the KPII solution agree with those of the $   {\mathrm{Gr}(N,M)_{> 0}}$ KP  soliton    under small perturbations.

We prove Theorem \ref{T:intro-stability}  through the inverse scattering theory (IST). 
Based on the Lax pair  
\begin{equation}\label{E:KPII-lax-1}
\left\{ 
{\begin{array}{l}
 (-\partial_{x_2}+\partial_{x_1}^2+u )\Phi(x,\lambda)=0,\\
 (-\partial_{x_3}+ \partial_{x_1}^3+\frac 32u\partial_{x_1}+\frac 34u_{x_1}+\frac 34\partial_{x_1}^{-1}u_{x_2}-\lambda^3   )\Phi (x,\lambda)=0, 
\end{array}}
\right.
\end{equation} a $\overline\partial$ formulation (Cauchy integral equation approach) of the IST for  decaying initial data was completed   \cite{ABF83,Li86,GN88,Gr97,W87}. 
Pioneering research on the IST for perturbed $ {\mathrm{Gr}(1,2)_{> 0}}$ KP  solitons was done   \cite{VA04,BP301}. Novel fundamental  contributions on the direct problem of the  IST for perturbed $  {\mathrm{Gr}(N,M)_{\ge 0}}$ KP  solitons, in particular, introducing the Sato theory, constructing the Green function, discovering multi-valued properties of the Green function at $\kappa_j$, deriving boundedness of the discrete part of the Green function, as well as verifying   the $\mathcal D^\flat$-symmetry of the Sato eigenfunctions and $\mathcal D^\sharp$-symmetry of the eigenfunctions (see \eqref{E:intro-sym-N-D-new}) \cite{P00,BP302,BP309,BP310,BP211,BP212,BP214}, have been established by Boiti, Pempinelli, Pogrebkov, and Prinari. Building upon Boiti et al's work, rigorous direct scattering theory for perturbed $  {\mathrm{Gr}(N,M)_{> 0}}$ KP  solitons is carried out in  \cite{Wu18,Wu20,Wu21}. We sketch the theory here and defer precise definitions to Section \ref{S:preliminaries}:    for
\be\label{E:intro-ini-data}
  \begin{array}{l}
  u_0(x_1,x_2)=u_s(x_1,x_2,0)+v_0(x_1,x_2),\\
  \textit{$ u_s(x)$  a ${\mathrm{Gr}(N,M)_{> 0}}$ KP soliton,  } 
  \textit{$   \sum_{|l|\le { d+8}} |{ (1+|x_1|+|x_2|)}\partial_x^lv _0|_  {L^1\cap L^\infty} \ll 1$,  $d\ge 0$,}\\
    \textit{$z_1=0$, $\{z_n,\kappa_j\}_{1\le n\le N,1\le j\le M}$  distinct reals, $
 \det  (
\frac 1{\kappa_{k}-z_{h}}
 )_{1\le k, h\le N}\ne 0$,}
 \end{array}
 \ee we  
 \begin{itemize}
     \item [(1)] derive, for $\forall\lambda\in\CC\backslash\{z_n,\kappa_j\}$, the unique existence of the eigenfunction  $m_0$ satisfying
   \begin{gather}
(-\partial_{x_2}+\partial_{x_1}^2 +2\lambda\partial_{x_1}
+u_0(x_1,x_2))m_0(x_1,x_2 ,\lambda)=0, \ \label{E:intro-Lax}\\
\lim_{|x|\to\infty}m_0(x_1,x_2,\lambda)=  \widetilde\chi(x_1,x_2,0,\lambda)=\frac{(\lambda-z_1)^{N-1}}{\Pi_{2\le n\le N}(\lambda-z_n)}\chi(x_1,x_2,0,\lambda); \label{E:intro-bdry} 
\end{gather}

       \item [(2)]  construct the {\bf forward scattering transform}
\be\label{E:intro-SD} 
\mathcal S(u_0,\{ z_n\})=(\{ z_n\},\{\kappa_j\}, \mathcal  D,s _c(\lambda)) 
\ee   which is continuous at $u_s$ and $\mathcal S(u(x),\{ z_n\})$ evolves linearly in $x_3$ if $u$ is a KPII solution;
\item [(3)] justify, for $\forall\lambda\in\CC\backslash\{z_n\}$,   the  Cauchy integral equation and the $\mathcal D$-symmetry  
\begin{gather}   
 {  {   m}_0(x_1,x_2, \lambda) =1+\sum_{n=1}^N\frac{   m_{0,z_n, \res }(x_1,x_2  )}{\lambda -z_n }  +\mathcal C  T_0
     m_0 ,}  \label{E:intro-CIE-0} \\
{ (e^{\kappa_1x_1+\kappa_1^2x_2 }m_0(x_1,x_2,\kappa^+_1),\cdots,e^{\kappa_Mx_1+\kappa_M^2x_2 }m_0(x_1,x_2,\kappa^+_M))\mathcal D=0.}\label{E:intro-sym-0}
     \end{gather} 
 \end{itemize} 
 
In this paper we complete the IST for perturbed $  {\mathrm{Gr}(N, M)_{> 0}}$ KP solitons   by solving the inverse scattering problem, 
\begin{theorem}\label{T:intro-inverse}
There exist an eigenfunction space $W$ and a positive constant $\epsilon_0\ll 1$, such that for any {\bf $d$-admissible scattering data} $\mathcal S =(\{z_n\},\{\kappa_j\},\mathcal D,s_c(\lambda))$, the system  of the Cauchy integral equation ({\bf CIE}) and the {\bf $\mathcal D$-symmetry},
\begin{gather}   
 {  {   m}(x, \lambda) =1+\sum_{n=1}^N\frac{   m_{z_n, \res }(x  )}{\lambda -z_n }  +\mathcal C  T
     m ,\ \lambda\neq z_n,}  \label{E:intro-CIE} \\
{ (e^{\kappa_1x_1+\kappa_1^2x_2+\kappa_1^3x_3}m(x,\kappa^+_1),\cdots,e^{\kappa_Mx_1+\kappa_M^2x_2+\kappa_M^3x_3}m(x,\kappa^+_M))\mathcal D=0 },\label{E:intro-sym}
     \end{gather} are uniquely solved in $W$ satisfying
\be\label{E:intro-asymp}
  { \sum_{0\le l_1+2l_2+3l_3\le d+5}| \partial^l_{x}\left[m(x ,\lambda)-\widetilde \chi (x ,\lambda)\right]|_{W}\le C\epsilon_0}.
 \ee

Moreover,      
\begin{gather}
 \left(-\partial_{x_2}+\partial_{x_1}^2+2 \lambda\partial_{x_1}+ u (x) \right)  m (x ,\lambda)=0 ,\label{E:intro-Lax-u}\\\
 u(x )\equiv - 2  \partial_{x_1}\sum _{n=1}^N    m_{z_n,\res}(x  )+\frac i\pi\partial_{x_1}\iint  T  m  \ d\overline\zeta\wedge d\zeta ,\label{E:intro-u-rep}\\
\sum_{0\le l_1+2l_2+3l_3\le d+4}|\partial^l_x\left[u(x )-u_s(x )\right]|_{L^\infty }\le C \epsilon_0,\label{E:intro-Lax-u-asymp}
\end{gather}
and $u :\mathbb R\times\mathbb R\times \mathbb R^+\to \mathbb R$ solves the KPII equation
\be\label{E:intro-KP-ist}
(-4u_{x_3}+u_{x_1x_1x_1}+6uu_{x_1})_{x_1}+3u_{{x_2}{x_2}}=0.
\ee
 
\end{theorem}

Precise definition of the {\bf eigenfunction space $W=W_x $}, the {\bf continuous scattering operator $T$}, $d$-admissibility of $\mathcal S$  with relations to  $\epsilon_0$,  $\widetilde\chi$, and $u_s$ are deferred to Section \ref{S:preliminaries}. We shall define the {\bf inverse scattering transform} $\mathcal S^{-1} (\{z_n\},\{\kappa_j\},\mathcal D,s_c )$ by the representation formula \eqref{E:intro-u-rep}. Continuity of $\mathcal S^{-1}$ at each $d$-admissible scattering data $(\{z_n\},\{\kappa_j\},\mathcal D,s_c )$  then follows from \eqref{E:intro-Lax-u-asymp}. Finally, for initial data $u_0$ satisfying \eqref{E:intro-ini-data}, the scattering data $\mathcal S (u_0, \{z_n\})$ defined by \eqref{E:intro-SD} is $d$-admissible and $m_0 \in W_{x_1,x_2,0} $. Hence Theorem \ref{T:intro-stability} follows from continuities of the direct scattering transform $\mathcal S$ and the inverse scattering transform $\mathcal S^{-1}$.

The system of \eqref{E:intro-CIE} and \eqref{E:intro-sym} represents both analytic and algebraic aspects of the IST toward an understanding of the KPII equation. Unique solvability of the system is verified by
proving the convergence of the iteration sequence 
\begin{align}
 \phi^{(k)}(x ,\lambda)  
\equiv&  1+\sum_{n=1}^N\frac{\color{black}\phi^{(k)}_{z_n,\res} (x   )}{\lambda -z_n}  +\mathcal CT \phi^{(k-1)}(x, \lambda)  ,\ \ k>0,\label{E:recursion-iteration-introduction}\\
\phi^{(0)}(x ,\lambda)\equiv &\widetilde\chi_{z_n,\kappa_j,A},\label{E:recursion-iteration-bdry-introduction}
 \end{align}
where  residues $\phi^{(k)}_{z_n,\res} $ are constructed by using the $\mathcal D$-symmetry 
\be\label{E:recursion-iteration-D-introduction}
(e^{\kappa_1x_1+\kappa_1^2x_2+\kappa_1^3x_3}\phi^{(k)}(x,\kappa^+_1), \cdots,e^{\kappa_Mx_1+\kappa_M^2x_2+\kappa_M^3x_3}\phi^{(k)}(x,\kappa^+_M))\mathcal D=0,
\ee and evaluation of CIE \eqref{E:recursion-iteration-introduction} of $\phi^{(k)} $ at $\kappa_j^+=\kappa_j+0^+$. Applying 
$d$-admissibility, the $\mathcal D$-symmetry, 
and Sato theory (see Proposition \ref{P:N-alg-sym}), convergence of the iteration sequence reduces to proving the  {\bf Cauchy integral operator (CIO)} estimate  (see Theorem \ref{T:CIO}),
\be\label{E:intro-main-est-0}
|(\mathcal CT)^n \corr{\widetilde \chi}|_{W}\le (C\epsilon_0)^n|\corr{\widetilde \chi}|_{W} ,\ \ \forall n\ge 0,
\ee for \corr{$\widetilde \chi\equiv\widetilde\chi_{z_n,\kappa_j,A}$}, 
in particular, small $L^\infty$-estimates  of the CIO near $\kappa_j$ (see Theorem \ref{T:basic}): 
 \be\label{E:intro-main-est}
 \begin{gathered}
|E_{\kappa_j}(\mathcal CTE_{\kappa_j})^n \corr{\widetilde \chi}|_{W}\le (C\epsilon_0)^n|E_{\kappa_j}\corr{\widetilde \chi}|_{W},\ \ \forall n\ge 0,\ 1\le j\le M,\\
|E_{\kappa_j}\phi|_W 
\equiv   | \phi^\flat     |_{ L^\infty(D_{\kappa_j})}+| \phi^\sharp     |_{\color{black} C^\mu_{\tilde\sigma}(D_{\kappa_j,\frac{1}{\tilde\sigma})}\cap L^\infty(D_{\kappa_j})}.
\end{gathered}
\ee
 Here $E_{\kappa_j}$ are the characteristic functions of   $D_{\kappa_j}$,   $E_{\kappa_j}\phi=E_{\kappa_j}(\phi^\flat+\phi^\sharp )$, 
$\phi^\flat  $ is a multi-valued bounded function, and $\phi^\sharp $, vanishing at $\kappa_j$,  is \corr{a bounded function on $D_{\kappa_j}$ and  $\tilde\sigma$-rescaled H$\ddot{\mbox{o}}$lder continuous  on $D_{\kappa_j,\frac{1}{\tilde\sigma}}$} (see Definition \ref{D:phase} for details).    
  They   are major estimates of this paper and are non trivial  because the CIO near $\kappa_j$ is oscillatory, non homogeneous, non symmetric,  and singular.

 We take advantage of the two dimensional property of   $\mathcal C$ and  special features of $T$ to prove \eqref{E:intro-main-est}.  
  More precisely,    scaling invariant properties enable us to decompose the leading terms of   CIO's near $\kappa_j$ to be sums of uniformly compact \corr{$\tilde s$-}domain and non uniformly compact \corr{$\tilde s$-}domain integrals (see \eqref{E:scale-mu}). For the uniformly compact \corr{$\tilde s$-}domain integrals, multi-valued properties and estimates can be derived by applying Stokes' theorem and H$\ddot{\mbox o}$lder interior estimates.  For the non uniformly compact \corr{$\tilde s$-}domain integrals, we write them as iterated integrals in polar coordinates. For instance, one type of  non uniformly compact \corr{$\tilde s$-}domain integrals  are written as, $\lambda=\kappa_j+ re^{i\alpha}=\kappa_j+ \frac{\tilde r}{{\tilde\sigma}}e^{i\alpha}$, $  \zeta=\kappa_j+  se^{i\beta}=\kappa_j+ \frac{\tilde s}{{\tilde\sigma}}e^{i\beta}$ (see \eqref{E:scale-homo}),  
\[
I_5= -\frac { \theta(\tilde r-1)}{2\pi i}  \int_{-\pi}^\pi  d\beta[\partial_ {\beta}   \ln (1-\gamma_j |\beta| )] 
       \int_{0 }^{\widetilde\sigma \delta} \frac{  e^{-i\wp(x,\kappa_j+\frac {\tilde s}{{\tilde\sigma}}e^{i\beta})}\widetilde\chi(x,\kappa_j+\frac {\tilde s}{{\tilde\sigma}}e^{-i\beta})}{\tilde s -  \tilde r e^{ i (\alpha-\beta) }} d\tilde s     
\] 
where $\wp( x,\zeta)=i((\overline\zeta- \zeta){x_1}+(\overline\zeta^2-\zeta^2){x_2}+(\overline\zeta^3-\zeta^3){x_3})$ is the {\bf phase function} of  $T$.   Therefore taking the advantage of $\widetilde s$-meromorphic properties,    estimates can be derived by  applying  the deformation method  and  stationary point analysis of $\wp$.

Notice that, due to singular structures of  $\mathcal CT$, the IST fails to carry the decay information of  the initial data to the KPII solution $u(x)$, but the above argument will be bootstrapped to carry the regularities to any order.
  
The paper is organized as follows:  notation and definitions are provided in Section \ref{S:preliminaries}. In Section \ref{S:CIO-0}, after explaining and highlighting features of the approach, we derive CIO estimates near $\kappa_j$ in Subsection \ref{SS:CIO-j}. Complete CIO estimates  on $W$ are provided in Subsection \ref{SS:CIO}. 

In Section \ref{S:eigenfunction}, we prove  \eqref{E:intro-CIE}-\eqref{E:intro-asymp} of Theorem \ref{T:intro-inverse} by justifying convergence of the iteration sequence \eqref{E:recursion-iteration-introduction}, \eqref{E:recursion-iteration-bdry-introduction}. 

In Section \ref{S:potential}, applying the heat operator to the iteration sequence \eqref{E:recursion-iteration-introduction}, \eqref{E:recursion-iteration-bdry-introduction}, and taking the limit by  estimates of the CIO and unique solvability of the CIE, we prove \eqref{E:intro-Lax-u}-\eqref{E:intro-Lax-u-asymp}  of Theorem \ref{T:intro-inverse}, construct the inverse scattering transform $\mathcal S^{-1}$, and prove the continuity.

In Section \ref{S:cauchy-kp}, we prove \eqref{E:intro-KP-ist}  by  the Lax pair approach.  Combining with the direct scattering theory \cite{Wu20, Wu21}, the Cauchy problem of perturbed $  {\mathrm{Gr}(N, M)_{> 0}}$ KP solitons is solved and  a uniqueness theorem (independent of $\{z_n\}$) is proved. These together amount  to a proof of Theorem \ref{T:intro-stability}. 

We provide examples in Section \ref{S:eigenfunction}-\ref{S:cauchy-kp} to illustrate how the $L^\infty$-stability of $ {\mathrm{Gr}(1,2)_{> 0}}$ KP  solitons \eqref{E:line-tau-oblique} is induced via the IST.

{\bf Acknowledgments}.  I would like to give special thanks to A. K. Pogrebkov and Y. Kodama for their kind and genius input. 
This research   was   supported by NSC 109-2115-M-001 -003 -. 
\section{Notation, definitions, and preliminaries}\label{S:preliminaries}

For the direct scattering problem,  the boundary data \eqref{E:intro-bdry} for the Lax equation is defined by  $\widetilde\chi  =\widetilde\chi_{\{z_n\},\{\kappa_j\},A}  =\frac{(\lambda-z_1)^{N-1}}{\Pi_{2\le n\le N}(\lambda-z_n)}\chi $ where   
\be\label{E:chi-intro}
\begin{split}
&\chi =\chi_{\{\kappa_j\},A}(x ,\lambda)=\frac{ 1}{\tau(x )}{\sum_{1\le j_1<\cdots<j_N\le M}\Delta_{j_1,\cdots,j_N}(A)(1-\frac{\kappa_{j_1}}{\lambda})\cdots(1-\frac{\kappa_{j_N}}{\lambda})E_{j_1,\cdots,j_N}(x )}
\end{split}
\ee  is the normalized Sato eigenfunction which is a rational function of $\lambda$ whose only singularity is a pole of order $N$ at $\lambda=0$ \cite[Theorem 6.3.8., (6.3.13) ]{D91}, \cite[Proposition 2.2, (2.21)]{K17}. Here $z_n$ are introduced to make poles of $\widetilde\chi$ simple. The simple pole condition is used to reduce  CIO estimates  near $z_n$ to standard Hilbert transform estimates in the inverse problem (see Proposition \ref{P:hilbert}). We shall prove that $z_n$ are auxiliary parameters to solve the Cauchy problem of the KPII equation by showing  different sets of $\{z_n\}$ yield the same solution in Theorem \ref{T:u-WP}. 

Given $u_0$ satisfying  \eqref{E:intro-ini-data},   the forward scattering transform $\mathcal S(u_0,  \{z_n\})=(\{z_n\},\{\kappa_j\}, \mathcal D , s _c(\lambda))$  \cite[Theorem 4]{Wu21} is defined by $\{z_n\}$ and $\{\kappa_j\}$ which are blow-up and multi-valued points of  $m_0$ respectively; $\mathcal  D$ are norming constants between values of $m_0$ at $\lambda=\kappa_j^+=\kappa_j+0^+$ and can be computed by
\be\label{E:intro-sym-N-D-new}
\begin{split}
&  \mathcal  D = \widetilde{\mathcal  D} \times \left({\tiny\begin{array}{ccc}\widetilde{\mathcal  D}_{11}&\cdots &\widetilde{\mathcal  D}_{1N}\\\vdots&\cdots &\vdots\\
\widetilde{\mathcal  D}_{N1}&\cdots &\widetilde{\mathcal  D}_{NN}\end{array}}\right)^{-1}\textit{diag}\,(\kappa_1^N,\cdots,\kappa_N^N) ,\\
&\widetilde{\mathcal  D}=\textit{diag}\,(\frac{\Pi_{2\le n\le N}(\kappa_1-z_n)}{(\kappa_1-z_1)^{N-1}}, \cdots, \frac{\Pi_{2\le n\le N}(\kappa_M-z_n)}{(\kappa_M-z_1)^{N-1}})\mathcal  D^\sharp  
\\
& { \mathcal  D^\sharp}= 
  \left({\mathcal  D}_{ji}^\sharp\right)= \left(\mathcal  D^\flat_{ji}+\sum_{l=j}^M\frac{c_{jl}\mathcal  D^\flat_{li}}{ 1-c_{jj}} \right), \\
& {\mathcal  D}^\flat= \textit{diag}\,(   
\kappa^N_1 ,\cdots,\kappa^N_M )\, A^T,  
\end{split}\ee where $c_{jl}=-\int\Psi _j(x_1,x_2,0 ) v_0(x_1,x_2)\varphi_l(x_1,x_2,0  )dx_1dx_2$, $\Psi_j(x)$, $\varphi_l(x)$ are residues of the adjoint eigenfunction at $\kappa_j$ \cite[(3.17))]{Wu21} and values of the Sato eigenfunction at $\kappa_l$ \cite[(2.9)]{Wu21};  $s_c(\lambda)$ is the continuous scattering data, arising from the $\overline\partial$-characterization
\be \label{E:conti-sc-debar}
\partial_{\overline\lambda}m_0(x_1,x_2, \lambda)
=  s_c(\lambda) e^{(\overline\lambda-\lambda)x_1+(\overline\lambda^2-\lambda^2)x_2  }m_0{(x_1,x_2, \overline\lambda)},\ \lambda\notin\RR,\ee and is a nonlinear Fourier transform of the initial perturbation
{\be \label{E:conti-sc}
\begin{split}
  s_c(\lambda) =& \frac { \Pi_{2\le n\le N}(\overline\lambda-z_n)}{(  \overline\lambda-z_1)^{N-1}   }\frac {\sgn(\lambda_I)}{2\pi i} \iint e^{-[(\overline\lambda-\lambda)x_1+(\overline\lambda^2-\lambda^2)x_2 ] }  \\
  \times& \xi(x_1,x_2,0, \overline\lambda) v_0(x_1,x_2)m_0(x_1,x_2, \lambda)dx_1dx_2,  
\end{split}
\ee}and $\lambda=\lambda_R+i\lambda_I,\,\overline\lambda=\lambda_R-i\lambda_I$, $ \xi(x,\lambda)
$ is the normalized Sato adjoint eigenfunction \cite[(2.6)]{Wu21}. 
 Algebraic and analytic constraints for scattering data  
 $\mathcal S(u_0,\{z_n\})$ are
\begin{align}
&\mathcal  D=\left({\tiny\ba{ccc}\kappa_1^N&\cdots&0\\
\vdots&\ddots&\vdots\\
0&\cdots&\kappa_N^N\\
\mathcal D_{N+1 1}&\cdots&\mathcal D_{N+1 N}\\
\vdots&\ddots&\vdots\\
\mathcal D_{M 1}&\cdots&\mathcal D_{M N}
\ea}\right),  \label{E:intro-sym-N-D}\\ 
& s_c(\lambda)=
\left\{
{\ba{ll}
 {\frac{ \frac {i}{ 2} \sgn(\lambda_I)}{\overline\lambda-\kappa_j}\frac{\gamma_j}{1-\gamma _j|\alpha|}}+\sgn(\lambda_I)  h_j(\lambda),&\lambda\in   D^ \times_{ \kappa_j },\\
\sgn(\lambda_I) {  \hbar_n}(\lambda),&\lambda\in    D^\times _{ z_n},
\ea}
\right.\label{E:intro-s-c-N}
\end{align}
 and   
 \begin{align}
 &  {\begin{array}{l} 
 |(1-\sum_{j=1}^ME_{{\kappa_j}}  )  \sum_{|l|\le {d+8}}|\left(|\overline\lambda-\lambda|^{l_1}   +| \overline\lambda^2-\lambda^2|^{l_2}\right) s_c (\lambda)|  _{ 
  L^\infty} \\
   +   \sum_{j=1}^M(|\gamma_j|+|h_j|_{ L^\infty(D_{\kappa_j})})+\sum_{n=1}^N|\hbar_n|_{ C^1(D_{z_n})}  \\
    + | {\corr{\textit{diag}\,(q_1, \cdots, q_M)^{-1}}\times\mathcal D \times \corr{\textit{diag}\,(q_1, \cdots, q_N)}  - { \mathcal D}^\flat  }|_{L^\infty}\\
  \le  {C\sum_{|l|\le  { d+8}} |{  (1+|x_1|+|x_2|) } \partial_{x }^{l }       v_0|_{L^1\cap L^\infty}} , \end{array}}\label{E:s-c-ana}\\
  &s_c(\lambda)=  \overline{s_c( \overline\lambda)},
  h_j(\lambda)=-\overline{h_j( \overline\lambda)}, 
 \hbar_n(\lambda)=-\overline{\hbar_n( \overline\lambda)},\ \label{E:s-c-reality}\\
 &q_j=\frac{\Pi_{2\le n\le N}(\kappa_j-z_n)}{(\kappa_j-z_1)^{N-1}} \textit{ for }1\le j\le M.\nonumber
 \end{align} \cite[Theorem 3]{Wu21}  (cf \cite{Wu24} for more detailed explanation about  $|h_j|_{ L^\infty(D_{\kappa_j})}$, $|\hbar_n|_{ C^1(D_{z_n})}$). 
 Here $D_{z,a\delta}=\{   \lambda=z+re^{i\alpha}:0\le r\le {a \delta},|\alpha|\le\pi\}$, $D_{z,a\delta}^\times=D_{z,a\delta}\backslash\{z\} $, $\delta=\frac 12\inf \{|z-z'|:z, z'  \in  \{ z_n,\, \kappa_j \},\ z\ne z'\}$,   $  E_{z,a\delta}(\lambda)\equiv 1$ on $D_{z,a\delta}$, $  E_{z,a\delta} (\lambda)\equiv 0$ elsewhere.  We suppress the $a\delta$-dependence for simplicity if $a=1$.

 Let us note that the forward scattering transform $\mathcal S$ is continuous at any $  {\mathrm{Gr}(N, M)_{> 0}}$ KP soliton $u_0$ by equipping topologies of $\mathcal S$ and  $u$ by the norms on both sides of the analytic constraint \eqref{E:s-c-ana}. 

A linearization theorem for the scattering data $\mathcal S(u(x),z_n)$ is proved:   if $\Phi= e^{ \lambda  x_1+ \lambda ^2x_2}    m(x, \lambda)$ satisfies the Lax pair \eqref{E:KPII-lax-1} and  
\begin{gather}
\partial_{\overline\lambda}  m(x, \lambda)=  {  {  s}_c(\lambda,x_3)}e^{(\overline\lambda-\lambda)x_1+(\overline\lambda^2-\lambda^2)x_2}   m(x,\overline\lambda)  ,\label{E:debar-evol}\\
  (e^{\kappa_1x_1+\kappa_1^2x_2}m(x,\kappa^+_1),\cdots,e^{\kappa_Mx_1+\kappa_M^2x_2}m(x,\kappa^+_M)){ \mathcal D(x_3)}=0,\label{E:D-evol}
\end{gather} then  
\be\label{E:linearization-D-evol}
\begin{gathered}
 z_n(x_3)\equiv z_n,\quad  \kappa_j(x_3)\equiv \kappa_j,\\ {  s}_c(\lambda, x_3)=   {e^{ (\overline\lambda^3-{ \lambda}^3)x_3}}{  s}_c(\lambda ),\quad
   {\mathcal {\mathcal D}}_{mn}(x_3)=  { e^{(\kappa_m^3-\kappa_n^3)x_3}}  {\mathcal D}_{mn}  
   \end{gathered}\ee  \cite[Lemma 4.1]{Wu20}, \cite[Theorem 5]{Wu21}. 	Notice that, $\kappa_j $, speeds of the line solitons, are fixed under the perturbation. This is distinct from the KdV equation. For KdV solitons, $\kappa_j$ change generically no matter how small the perturbation is. What's more, the number of $\kappa_j$ of the KdV equation can be different under perturbation.

For the inverse problem, introduce  
\begin{definition}\label{E:generic-sd}
Given $0<\epsilon_0\ll 1$,  $d\ge 0$, and a $  {\mathrm{Gr}(N, M)_{> 0}}$ KP soliton $u_s$ defined by $\{\kappa_j\}, A$, a scattering data  $ {\mathcal S} =(\{z_n\},\{\kappa_j\}, \mathcal D ,s _c(\lambda))$ is called    $d$-admissible   if  
\begin{align}
&\mathcal  D=\left( \textit{\tiny${\ba{ccc}\kappa_1^N&\cdots&0\\
\vdots&\ddots&\vdots\\
0&\cdots&\kappa_N^N\\
\mathcal D_{N+1 1}&\cdots&\mathcal D_{N+1 N}\\
\vdots&\ddots&\vdots\\
\mathcal D_{M 1}&\cdots&\mathcal D_{M N}
\ea}$}\right),  \label{E:sym-N-D}\\ 
& s_c(\lambda)=
\left\{
{\ba{ll}
 {\frac{ \frac {i}{ 2} \sgn(\lambda_I)}{\overline\lambda-\kappa_j}\frac{\gamma_j}{1-\gamma _j|\alpha|}}+\sgn(\lambda_I)  h_j(\lambda),&\lambda\in   D^ \times_{ \kappa_j },\\
\sgn(\lambda_I) {  \hbar_n}(\lambda),&\lambda\in    D^\times _{ z_n},
\ea}
\right.\label{E:s-c-N}
\end{align} and 
\begin{align}
&\det  (
\frac 1{\kappa_{k}-z_{h}}
 )_{1\le k, h\le N}\ne 0,\  z_1=0,\{z_n,\kappa_j\}  \textit{ distinct real},\label{E:intro-sym-N-D-flat}\\
&
\epsilon_0\ge     (1-\sum_{j=1}^ME_{{\kappa_j}}  )  \sum_{|l|\le {d+8}}|\left(|\overline\lambda-\lambda|^{l_1} +| \overline\lambda^2-\lambda^2|^{l_2}\right) s_c (\lambda)|   _{    
 L^\infty}\label{E:epsilon-0-BPP-new}
 \\
 &\qquad +  \sum_{j=1}^M(|\gamma_j|+|h_j|_{ L^\infty(D_{\kappa_j})})+\sum_{n=1}^N|\hbar_n|_{ C^1(D_{z_n})}\nonumber\\
 &\qquad+   | {\corr{\textit{diag}\,(q_1, \cdots, q_M)^{-1}}\times\mathcal D \times \corr{\textit{diag}\,(q_1, \cdots, q_N)}  - { \mathcal D}^\flat  }|_{L^\infty},\nonumber\\
 &s_c(\lambda)=  \overline{s_c( \overline\lambda)},
  h_j(\lambda)=-\overline{h_j( \overline\lambda)}, 
 \hbar_n(\lambda)=-\overline{\hbar_n( \overline\lambda)},\label{E:N-reality-inv}\\
& { \mathcal D}^\flat = \textit{diag}\,(   
\kappa^N_1 ,\cdots,\kappa^N_M )\, A^T,\ \ q_j=\frac{\Pi_{2\le n\le N}(\kappa_j-z_n)}{(\kappa_j-z_1)^{N-1}} \textit{ for }1\le j\le M. \label{E:intro-sym-N-D-flat-new}
\end{align}     
\end{definition}
From the direct scattering theory discussed above, if $u_0$ satisfying  \eqref{E:intro-ini-data}  then  the forward scattering transform 
\be\label{E:ini-admissible}
\begin{gathered}
\textit{$\mathcal S(u_0,  \{z_n\})$ is   $d$-admissible  with $
\epsilon_0
  \le   C\sum_{|l|\le {d+8}} |{  (1+|x_1|+|x_2|) }\partial_{x }^{l }  v_0|_{L^1\cap L^\infty}$}.
  \end{gathered}
\ee
 
 Finally, 
$\mathcal C$ is the Cauchy integral operator,  $T$  is the {continuous scattering operator}   
\begin{align}
\mathcal C \phi  (x,\lambda)
\equiv=& -\frac{1}{2\pi i}\iint_{\CC}\frac{\phi(x,\zeta)}{\zeta-\lambda}d\overline\zeta\wedge d\zeta,\label{E:ct-operator} \\
 T \phi  (x,\lambda)
\equiv& {  s}_c(\lambda  )e^{(\overline\lambda-\lambda)x_1+(\overline\lambda^2-\lambda^2)x_2+(\overline\lambda^3-\lambda^3)x_3  }\phi(x, \overline\lambda),\label{E:cauchy-operator}
\end{align}  
and $T_0\phi(x_1,x_2,\lambda)=T\phi|_{(x_1,x_2,x_3=0)}$.

\begin{definition}\label{D:phase} Given  $\{z_n,\kappa_j \}$, the {eigenfunction  space} ${ W }=W_x$  consists of $\phi$ satisfying 
\begin{itemize}
\item [$(a)$] $\phi (x, \lambda)=\overline{ \phi (x, \overline\lambda)};$ 
\item [$(b)$] $(1-  \sum_{n=1}^NE_{z_n} )\phi(x, \lambda)\in L^\infty;$ 
\item [$(c)$] for $\lambda \in D_{z_n}^\times$, 
$
\phi(x, \lambda)=\frac{ {\phi_{z_n,\res}(x)}}{\lambda-z_n }  +\phi_{z_n,r}(x, \lambda)$, $ \phi_{z_n,\res}$,   $\phi_{z_n,r}  \in L^\infty( D_{z_n})$; 
\item[$(d)$] for $\lambda =\kappa_j+re^{i\alpha}\in D_{\kappa_j}^\times$, $\phi=\phi^\flat+\phi^\sharp $, 
$\phi^\flat =\sum_{l=0}^\infty \phi_l(X)(-\ln(1-\gamma_j|\alpha|))^l \in L^\infty(D_{\kappa_j}) $, $\phi^\sharp \in C_{\tilde\sigma}^{\mu } (D_{\kappa_j,\color{black}\frac{1}{\tilde\sigma}}) \cap L^\infty (D_{\kappa_j})$, $\phi^\sharp(x,\kappa_j)=0$. 
\end{itemize}Here $C_{\tilde\sigma}^{\mu } (D_{\kappa_j,\color{black}\frac{1}{\tilde\sigma}})=C (D_{\kappa_j,\color{black}\frac{1}{\tilde\sigma}})\cap H_{\tilde\sigma}^{\mu } (D_{\kappa_j,\color{black}\frac{1}{\tilde\sigma}})$, the rescaled H$\ddot{\mbox{o}}$lder space $H_{\tilde\sigma}^{\mu } (D_{\kappa_j,\color{black}\frac{1}{\tilde\sigma}})$, for $z\in\RR$,  consists of  functions $\phi(x,\lambda)\equiv \phi(r,\alpha,X)$,  $X=(X_1,X_2,X_3)$,  satisfying
\be\label{E:rescaled-holder}
|\phi|_{ H_{\tilde\sigma}^\mu (D_{z,\color{black}\frac{1}{\tilde\sigma}})} \equiv  {\color{black}\sup_{\scriptsize{\ba{c}  \tilde r_1,  \tilde r_2\le {\color{black}1,}
 |\alpha_1|,|\alpha_2|\le \pi\ea}} \frac{|\phi(\frac{ \tilde r_1}{\tilde\sigma}, \alpha_1,X)-\phi(\frac{ \tilde r_2}{\tilde\sigma},\alpha_2,X)|}{ |\tilde r_1e^{i\alpha_1}-\tilde r_2e^{i\alpha_2}|^\mu} }\quad<\infty
\ee for $
\lambda_j=z+r_j e^{i\alpha_j}=z+\frac{\tilde r_j}{\tilde\sigma} e^{i\alpha_j}\in D_z$ with the rescaling parameter 
\be\label{E:sigma-parameter}
\tilde\sigma=   \max\{ 1, |X_1|, \sqrt{|X_2|}, \sqrt[3]{|X_3|}\},
\ee and $X_k$ are the coefficients of the phase function of  $T$, 
\be\label{E:phase}
\begin{split}
\wp( x,\lambda)=& i[(\overline\lambda- \lambda){x_1}+(\overline\lambda^2-\lambda^2){x_2}+(\overline\lambda^3-\lambda^3){x_3}] \\
=&  X_1r\sin\alpha  +X_2r^2 \sin2\alpha  +X_3r^3 \sin3\alpha  
\equiv \wp(r,\alpha,X),\\
 X_1( x,z)= &2(x_1+2x_2z +3x_3z^2),\   
 X_2(x,z)=   2(x_2 +3x_3z),\  
  X_3(x,z)= 2x_3  
 \end{split}
 \ee  
for $
\lambda=z+re^{i\alpha} \in D_z$.    
  For   $\phi\in W$,   define 
\begin{multline}\label{E:N-W}
|\phi|_W 
\equiv  |(1-\sum_{n=1}^N E_{z_n})\phi|_{L^\infty}  \\
 +\sum_{n=1}^N (| \phi_{z_n,\res}|_{L^\infty}+ | \phi_{z_n,r}|_{L^\infty(D_{z_n})}) 
  + \sum_{j=1}^M (| \phi^\flat     |_{ L^\infty(D_{\kappa_j})}+| \phi^\sharp     |_{ C^\mu_{\tilde\sigma}(D_{\kappa_j,\frac{1}{\tilde\sigma}})\cap L^\infty (D_{\kappa_j})})  . 
  \end{multline}

\end{definition} The direct scattering theory \cite[Theorem 4]{Wu21}  (or see \cite{Wu24} for more detailed explanation) implies, if $u_0$ satisfying  \eqref{E:intro-ini-data}  then    
\be\label{E:ini-eigen} 
\partial_x^lm_0 \in W_{(x_1,x_2,0)} \quad \textit{ for   }\ 0\le l_1+2l_2+3l_3\le {d+5}.
\ee

\section{Estimates for the Cauchy integral operator $\mathcal CT$}\label{S:CIO-0}

\subsection{The Cauchy integral operator $\mathcal CT$ near $\kappa_j$}\label{SS:CIO-j}

\subsubsection {The strategy}\label{SS:strategy} \hfill \\
 
To satisfy the $\mathcal D$-symmetry  \eqref{E:intro-sym-0}, an $L^\infty$-estimates of the CIO near $\kappa_j$ is inevitable. Hence the goal of this subsection is to establish 
\begin{theorem}\label{T:basic} Suppose $\mathcal S=(\{z_n\},\{\kappa_j\}, \mathcal D,s _c)$ is $d$-admissible, then, for $\forall n\ge 0$, $1\le j\le M$,
\[\begin{array}{r}
    \sum_{0\le l_1+2l_2+3l_3\le d+5} |E_{\kappa_j}\partial_x^l(\mathcal C  T E_{\kappa_j})^n \corr{ \widetilde \chi}| _W
 +  \sum_{0\le l_1+2l_2+3l_3\le d+5} |(1-E_{\kappa_j})\partial_x^l(\mathcal C  T E_{\kappa_j})^n \corr{ \widetilde \chi}| _{   L^\infty }   \\ 
 \le  ( C\epsilon_0)^n \sum_{0\le l_1+2l_2+3l_3\le d+5}|E_{\kappa_j}\partial_x^l\corr{ \widetilde \chi}| _W.   \end{array}\]
\end{theorem}

Because of singularities at $\lambda$ , $\kappa_j$,  no apparent symmetries for a cancellation of the leading term of $s_c$, and the highly oscillatory non homogeneous phase function  $\wp$ of the continuous scattering operator $T$, 
proof of Theorem \ref{T:basic} is not trivial and is the main technical estimate of this paper. Moreover,  argument in the proof will be used in later sections.

We make several observations and remarks  to explain the strategy and the topology $W$ equipped at $\kappa_j$.  
 To start, note the Cauchy integral operator is not singular.  
\begin{lemma}\label{L:vekua} If $\phi\in L^p(D_z)$, $p>2$, then for $\nu=\frac {p-2}p$,
\begin{eqnarray*}
|\mathcal CE_z\phi|_{L^\infty} \le   C|\phi|_{L^p(D_z)},\qquad
|\mathcal CE_z\phi|_{H^\nu(D_z)} \le   C|\phi|_{L^p(D_z)} . 
\end{eqnarray*}  
\end{lemma}
\begin{proof} Please see \cite[Theorem 1.19]{V62} for the details. 
\end{proof}
Hence    estimating $|\mathcal C TE_{\kappa_j}\corr{ \widetilde \chi}|_{C^\mu_{\tilde\sigma}(D_{\kappa_j,\corr{\frac1{\tilde\sigma}}})}$  can be reduced to estimating
 $|\mathcal C \widetilde\gamma_je^{-i\wp }E_{\kappa_j} \corr{ \widetilde \chi}|_{C^\mu_{\tilde\sigma}(D_{\kappa_j,\corr{\frac1{\tilde\sigma}}})}$ from \eqref{E:intro-s-c-N}, \eqref{E:s-c-ana}, \eqref{E:cauchy-operator}. Here, in terms of the polar coordinates 
 $ \zeta= \kappa_j+se^{i\beta}$, $0\le s \le \delta$,   $-\pi< \beta \le\pi  $,
\be\label{E:gamma-N}
\widetilde\gamma_j(\zeta)=\widetilde\gamma_j(s,\beta)= \frac{ \frac {i}{ 2} {\sgn}(\zeta_I)}{\overline\zeta-\kappa_j}\frac{\gamma_j}{1-\gamma_j|\beta|}  =-\frac {\frac i2   \partial_\beta \ln (1-\gamma_j|\beta|) }{\overline\zeta-\kappa_j}={-}\partial_{\overline\zeta}\ln (1-\gamma_j {|\beta|}).     
\ee 

The following lemma implies that  leading terms of $\mathcal C \widetilde\gamma_j e^{-i\wp }E_{\kappa_j}\corr{ \widetilde \chi} $  as well as the iteration  can be integrated and the outcomes are multi-valued functions at $\kappa_j$. 
\begin{lemma}\label{L:stokes} Fixed $\lambda=\kappa_j+re^{i\alpha}, \,\zeta=\kappa_j+se^{i\beta}\in D_{\kappa_j}^\times$,  for any non negative   integer $ l$, 
\begin{align*}
&\mathcal C\widetilde\gamma_jE_{\kappa_j}[-{ \ln (1-\gamma_j|\beta|)}]^l =\frac{[- { \ln (1-\gamma_j|\alpha|)}]^{l+1}}{l+1} -\frac {1}{2\pi i}\oint_{|\zeta-\kappa_j|= \delta}\frac{\frac 1{l+1}[-\ln (1-\gamma_j|\beta|)]^{l+1}} { \zeta-\lambda }d\zeta.
\end{align*}
\end{lemma} 
\begin{proof} Taking $l=0$, in view of \eqref{E:gamma-N} and applying Stokes' theorem,  
\begin{align}
 & \mathcal C\widetilde\gamma_jE_{\kappa_j}1   
=    -\frac {1}{2\pi i}\lim_{\epsilon\to 0}\left\{ \iint_{D_{\kappa_j}/D_{\kappa_j,\epsilon}\cup D_{\lambda,\epsilon}}\frac{ {-}\partial_{\overline\zeta}\ln (1-\gamma_j|\beta|)}{ \zeta-\lambda }d\overline\zeta\wedge d\zeta \right.\nonumber \\
&\left.+\iint_{ D_{\kappa_j,\epsilon} }\frac{\frac i2{\sgn}(\zeta_I) \frac{\gamma_j}{1-\gamma_j|\beta|}}{ (\overline\zeta-\kappa_j)(\zeta-\lambda) }d\overline\zeta\wedge d\zeta + \iint_{  D_{\lambda,\epsilon}}\frac{\frac i2{\sgn}(\zeta_I)\frac{\gamma_j}{1-\gamma_j|\beta|}}{ (\overline\zeta-\kappa_j)(\zeta-\lambda) }d\overline\zeta\wedge d\zeta\ \right\}\nonumber\\
&=    {-} { \ln (1-\gamma_j|\alpha|)}  {+}\frac {1}{2\pi i}\oint_{|\zeta-\kappa_j|= \delta}\frac{\ln (1-\gamma_j|\beta|)} { \zeta-\lambda }d\zeta .\nonumber
\end{align} Other identities for $l>0$ can be proved similarly.
\end{proof}

Moreover, dilating the polar coordinates near $z=\kappa_j$, 
\be\label{E:scale-homo}
\begin{split}
 \lambda=z+ re^{i\alpha}=z+ \frac{\tilde r}{{\tilde\sigma}}e^{i\alpha},&\quad  \zeta=z+  se^{i\beta}=z+ \frac{\tilde s}{{\tilde\sigma}}e^{i\beta},\hskip.23in r,\,s\le\delta,\,|\alpha|,|\beta|\le \pi,\\
  \tilde\lambda=z+ \tilde r e^{i\alpha},&\quad\tilde\zeta=z+\tilde se^{i\beta},\hskip1in
 \tilde r,\,{\tilde s} \le {\tilde\sigma}\delta  
 \end{split}
 \ee we have the scaling invariant property   
 \begin{align}
& \mathcal C_\lambda\widetilde\gamma_j e^{-i\wp(x,\zeta)} E_{\kappa_j,\delta}  f(\kappa_j+se^{i\beta}) \label{E:scaling-inv} \\
=&-\frac 1{2\pi i}\iint_{D _{ \kappa_j,\delta } }\frac {  \frac{ \frac {i}{ 2} {sgn}(\zeta_I)}{ \overline\zeta-\kappa_j}  \frac {+\gamma_j}{1-\gamma_j |\beta|} e^{ (\overline\zeta-\zeta) x_1+(\overline\zeta^2-\zeta^2)x_2+(\overline\zeta^3-\zeta^3)x_3  } f(  \overline\zeta)}{ \zeta-\lambda}d\overline\zeta\wedge d\zeta\nonumber\\
=&-\frac 1{2\pi i}\iint_{D _{ \kappa_j,\delta } }\frac {  \frac{ \frac {i}{ 2} {sgn}(\zeta_I)}{ \overline\zeta-\kappa_j}  \frac {+\gamma_j}{1-\gamma_j |\beta|} e^{ -i(X_1s\sin\beta  +X_2s^2 \sin2\beta +X_3s^3 \sin3\beta) } f(  \overline\zeta)}{ \zeta-\lambda}d\overline\zeta\wedge d\zeta\nonumber\\
=&-\frac 1{2\pi i}\iint_{D _{ \kappa_j,\tilde\sigma\delta } }\frac {   \frac{ \frac {i}{ 2} {sgn}({\tilde\zeta}_I)}{ \overline{\tilde\zeta}-\kappa_j}  \frac {+\gamma_j}{1-\gamma_j |\beta|} e^{ -i(\frac{X_1}{\tilde\sigma}\tilde s\sin\beta  +\frac{X_2}{\tilde\sigma^2}{\tilde s}^2 \sin2\beta  +\frac{X_3}{\tilde \sigma^3}\tilde s^3 \sin3\beta) }f(  \kappa_j+\frac{\overline{{\tilde\zeta}}-\kappa_j}{  \tilde\sigma})}{ {\tilde\zeta}-\tilde\lambda}d\overline{\tilde\zeta}\wedge d\tilde\zeta\nonumber\\
 =& \mathcal C_{\tilde\lambda}\widetilde\gamma_je^{-i\wp(\frac{\tilde s}{\tilde\sigma},\beta, X)}E_{\kappa_j,\tilde\sigma\delta}    f(\kappa_j+\frac{\tilde s}{\tilde\sigma} e^{i\beta}),\nonumber
\end{align}
 with $\tilde\sigma,\,X$ defined by \eqref{E:sigma-parameter}, \eqref{E:phase}. Therefore,
\begin{align}
&|\mathcal C_\lambda\widetilde\gamma_j e^{-i\wp(x,\zeta)} E_{\kappa_j}  f(\kappa_j+se^{i\beta}) |_{C_{\tilde\sigma}^{\mu } (D_{\kappa_j,\frac 1{\tilde\sigma}})\cap L^\infty(D_{\kappa_j }) }\label{E:scale-invariant}\\
&\hskip.9in= |\mathcal C_{\tilde\lambda}\widetilde\gamma_je^{-i\wp(\frac{\tilde s}{\tilde\sigma},\beta, X)}E_{\kappa_j,\tilde\sigma\delta}   f(\kappa_j+\frac{\tilde s}{\tilde\sigma} e^{i\beta})|_{C ^{\mu } (D_{\kappa_j,1})\cap L^\infty(D_{\kappa_j,\tilde\sigma\delta})}.\nonumber
\end{align} which explains by 
     introducing the rescaled H$\ddot{\mbox{o}}$lder continuous norm for $W$ at $\kappa_j$ (Definition \ref{D:phase}), one can tame the highly oscillatory properties of the CIO at $\kappa_j$.

More precisely,  using the change of variables \eqref{E:scale-homo} and  the scaling invariant property \eqref{E:scaling-inv},  
\begin{align}
&  \mathcal C_\lambda\widetilde\gamma_je^{-i\wp(x,\zeta)}E_{\kappa_j} f\equiv I_1+I_2+I_3+I_4+I_5,\label{E:scale-mu}
\end{align}
with
\begin{align}
I_1=&-\frac {\theta(1-\tilde r)}{2\pi i}\iint_{ \tilde s<2} \frac{  \widetilde \gamma_j(\tilde s, \beta)  f^{ \flat}(\frac {\tilde s}{ {\tilde\sigma} },-\beta,X) }{\tilde\zeta-\tilde \lambda}d\overline{\tilde\zeta} \wedge d\tilde\zeta  ,\label{E:scal-1}\\
I_2=&-\frac {\theta(1-\tilde r)}{2\pi i}\iint_{ \tilde s<2} \frac{  \widetilde \gamma_j(\tilde s, \beta)[ e^{-i\wp(\frac {\tilde s}{ {\tilde\sigma} },\beta,X)}-1]f^{\flat}(\frac {\tilde s}{ {\tilde\sigma} },-\beta,X) }{\tilde\zeta-\tilde \lambda}d\overline{\tilde\zeta} \wedge d\tilde\zeta ,\label{E:scal-2}\\
I_3=&-\frac {\theta(1-\tilde r)}{2\pi i}\iint_{ \tilde s<2} \frac{  \widetilde \gamma_j(\tilde s, \beta) e^{-i\wp(\frac {\tilde s}{ {\tilde\sigma} },\beta,X)}f^{ \sharp}(\frac {\tilde s}{ {\tilde\sigma} },-\beta,X) }{{\tilde\zeta}-\tilde \lambda}d\overline{\tilde\zeta} \wedge d{\tilde\zeta} ,\label{E:scal-3}\\
I_4=&-\frac {\theta(1-\tilde r)}{2\pi i}\iint_{2< \tilde s< {\tilde\sigma} \delta} \frac{ \widetilde \gamma_j(\tilde s, \beta) e^{-i\wp(\frac {\tilde s}{{\tilde\sigma}},\beta,X)}f(\frac {\tilde s}{{\tilde\sigma}},-\beta,X)  }{{\tilde\zeta}-\tilde \lambda}d\overline{\tilde\zeta} \wedge d{\tilde\zeta},\label{E:scal-4}\\
I_5=&-\frac {\theta(\tilde r-1)}{2\pi i}\iint_{  \tilde s< {\tilde\sigma} \delta} \frac{  \widetilde\gamma_j(\tilde s, \beta) e^{-i\wp(\frac {\tilde s}{{\tilde\sigma}},\beta,X)}f(\frac {\tilde s}{ {\tilde\sigma} },-\beta,X) }{{\tilde\zeta}-\tilde\lambda}d\overline{\tilde\zeta} \wedge d{\tilde\zeta},\label{E:scal-5}
\end{align}where $I_k=I_k(x, \lambda)$, and $f(x,\zeta)$ is identified as $f(s,\beta,X)$.

We shall apply Stokes' theorem  and H$\ddot{\mbox{o}}$lder interior estimates to derive estimates for $I_1$, $I_2$, and $I_3$ (see Proposition \ref{P:l-c-mu} in $\S$ \ref{SS:Holder}). 
For $I_4$, $I_5$, \corr{and $f=\widetilde\chi$}, integrals with {slow decaying kernels} on non uniformly compact \corr{$\tilde s$-}domains, we write them as iterated integrals in polar coordinates,
\begin{align*}
I_4=&-\frac {\theta(1-\tilde r)}{2\pi i}\int_{-\pi}^\pi  d\beta[\partial_ {\beta}   \ln (1-\gamma_j |\beta| )] 
       \int_{2< \tilde s< {\tilde\sigma} \delta} \frac{  e^{-i\wp(\frac {\tilde s}{ {\tilde\sigma} },\beta,X)}\corr{\widetilde \chi}(\frac {\tilde s}{{\tilde\sigma}},-\beta,X)}{\tilde s -  \tilde r e^{ i (\alpha-\beta) }} d\tilde s    , \\
I_5=&-\frac { \theta(\tilde r-1)}{2\pi i}  \int_{-\pi}^\pi  d\beta[\partial_ {\beta}   \ln (1-\gamma_j |\beta| )] 
       \int_{0 }^{\widetilde\sigma \delta} \frac{  e^{-i\wp(\frac {\tilde s}{ {\tilde\sigma} },\beta,X)}\corr{\widetilde \chi}(\frac {\tilde s}{{\tilde\sigma}},-\beta,X)}{\tilde s -  \tilde r e^{ i (\alpha-\beta) }} d\tilde s.      
\end{align*}Estimates for them will be derived by applying meromorphic  properties in $\tilde s$, the deformation method, and stationary point analysis of $\wp$. Main ideas are, near stationary points, we deform 
\[\tilde s\in\RR\quad\longrightarrow\quad \tilde se^{i\tau}\in\CC\] such that the segment of $ \tilde s $   turns into a union  of line segments $\Gamma$'s and arcs $S$'s, and
\begin{align}
 (a)\ &\ \textit{on $\Gamma$'s,    $\mathfrak {Re} ({-i\wp(\frac {\tilde se^{i\tau}}{ {\tilde\sigma} },\beta,X)})\le -\frac 1C |\sin (k\beta)|\tilde s^k $ }\label{E:deform-intro-1}  ;\\
 (b)\ &\  \textit{on $\Gamma$'s,  ${|{\tilde\zeta}-\tilde\lambda}|\ge \frac 1C\max\{\tilde s,\tilde r\}$;}\nonumber\\
(c)\ &\ \textit{on $S$'s,  $\mathfrak {Re} ({-i\wp(\frac {\tilde se^{i\tau}}{  \tilde\sigma},\beta,X)})\le 0$. }\nonumber
\end{align} 
  Therefore, if $\tilde\sigma=\sqrt[k]{|X_k|}$,   $k=3,2 $, estimates can be done if $|{{\tilde\zeta}-\tilde\lambda}|\ge \frac 1C$.  When $|{{\tilde\zeta}-\tilde\lambda}|\le \frac 1C$, we shall prove  $I_5$ near the stationary is no longer a singular integral and estimates can be obtained  by means of Lemma \ref{L:vekua} (see Proposition \ref{P:non-homogeneous-F3} in $\S$ \ref{SSS:FN-homo}).

Estimates for  $\tilde\sigma=|X_1|$ are more involved because $e^{\mathfrak {Re} ({-i\wp(\frac {\tilde se^{i\tau}}{ {\tilde\sigma} },\beta,X)})}   $   either decays too slowly or creates a $\frac 1{|\sin\beta|}$ singularity. We shall either take advantage of the scaling invariant properties of the Hilbert transform or   a finer decomposition to overcome difficulties (see Proposition \ref{P:non-homogeneous-F1} and remarks before Definition \ref{D:renormal} for subtle technical consideration in $\S$ \ref{SSS:FN-homo}).

After analysing $I_1$-$I_5$, we prove Theorem \ref{T:basic} in $\S$ \ref{SS:CIO-k-2}.

\subsubsection {Estimates for $I_1$, $I_2$, and $I_3$}\label{SS:Holder}  \hfill\\

Proofs in this subsection reminiscent of estimates of the Beltrami's equation (cf. \cite[\S 8, Chapter I]{V62}) thanks to the leading singular terms can be integrated by Stokes' theorem.

\begin{proposition}\label{P:l-c-mu} Suppose $\mathcal S=(\{z_n\},\{\kappa_j\}, \mathcal D,s _c)$ is $d$-admissible and  $f\in W$. 
\begin{align}
 |  I_1^\flat|_{L^\infty (D_{\kappa_1})}+ |   I_1^\sharp|_{C^{\mu }_{\tilde\sigma}(D_{\kappa_1,\corr{\frac{1}{\tilde\sigma}} })}\le &  C\epsilon_0 |f^\flat |_{L^\infty(D_{\kappa_1})} ,\label{E:I-1-est}\\
 |I_2|_{C^{\mu}_{\tilde\sigma}(D_{\kappa_1,\corr{\frac{1}{\tilde\sigma}}  })}\le & C\epsilon_0 |f^\flat |_{L^\infty(D_{\kappa_1})}  ,\label{E:I-2-est}\\
 | I_3|_{C^\mu_{\tilde\sigma}(D_{\kappa_j,\frac{1}{\tilde\sigma}})}\le & C\epsilon_0    | f^\sharp|_{C^\mu_{\tilde\sigma}(D_{\kappa_j,\frac{1}{\tilde\sigma} })} .\label{E:I-3-est}
\end{align}

\end{proposition}

\begin{proof}  
 
 Firstly, applying Lemma \ref{L:stokes},  
\[
\begin{gathered}
  I_1= \theta(1-\tilde r)F^\flat(x,\lambda)-\frac {\theta(1-\tilde r)}{2\pi i}\oint_{|{\tilde\zeta}-\kappa_1|= 2}\frac{F^\flat(x,{ \zeta} )} { {\tilde\zeta}-\tilde\lambda }d{\tilde\zeta} ,\\
f^\flat(x,\zeta)= \sum_{l=0}^\infty f_l(x)[-\ln(1-\gamma|\beta|)]^l,\\
F^\flat(x,\zeta )= \sum_{l=0}^\infty f_l(x)\frac{[-\ln(1-\gamma|\beta|)]^{l+1}}{l+1} .
\end{gathered}
\]Therefore,
\be\label{E:1-homo-1}
\begin{gathered}
  I_1^\flat= \theta(1-\tilde r){ F^\flat(x, \lambda )}-\frac {\theta(1-\tilde r)}{2\pi i}\oint_{|{\tilde\zeta}-\kappa_1|= 2}\frac{ F^\flat(x,{ \zeta} )} { {\tilde\zeta}-\kappa_1 }d{\tilde\zeta}  ,\ \
  I_1^\sharp=   I_1-  I_1^\flat,\\
|  I_1^\flat|_{L^\infty (D_{\kappa_1})}+ |   I_1^\sharp|_{C^{\mu }_{\tilde\sigma}(D_{\kappa_1,\corr{\frac{1}{\tilde\sigma}} })}\le  C\epsilon_0 |f^\flat |_{L^\infty(D_{\kappa_1})} .
\end{gathered}\ee

Besides,  from  Lemma \ref{L:vekua} and $
|\widetilde\gamma_1(\tilde s, \beta)[ e^{-i\wp(\frac {\tilde s}{ \tilde\sigma },\beta,X)}-1]|<C\epsilon_0$ for $\tilde s<2$, 
\be\label{E:1-homo-5}
|I_2|_{C^{\mu}_{\tilde\sigma}(D_{\kappa_1,\corr{\frac{1}{\tilde\sigma}}  })}\le  C\epsilon_0 |f^\flat |_{L^\infty(D_{\kappa_1})}  .
\ee

For $I_3$,  one has 
$
| \widetilde\gamma_1(\tilde s,\beta)f ^\sharp(\frac {\tilde s}{ \tilde\sigma},\beta,X)|_{L^\infty(D_{\kappa_1 })}\le C \epsilon_0  |f^\sharp|_{H^\mu_{\tilde\sigma}(D_{\kappa_1,\corr{\frac{1}{\tilde\sigma}}})}\tilde s^{\mu-1}$ from $  f^\sharp\in C^{\mu}_{\tilde\sigma}(D_{\kappa_1,\corr{\frac{1}{\tilde\sigma}}})$ and  $  f^\sharp(x,\kappa_1)=0$. Therefore, an improper integral yields
\be \label{E:mean-vekua}
 | I_3|_{L^\infty(D_{\kappa_1})}\le C\epsilon_0  |f^\sharp|_{H^\mu_{\tilde\sigma}(D_{\kappa_1,\corr{\frac{1}{\tilde\sigma}}})} .
 \ee
 
 To derive the $H^\mu_{\tilde\sigma}$-estimate of $I_3$, from Lemma \ref{L:vekua}, let $\tilde \lambda_j=\kappa_1+\tilde r_je^{i\alpha_j}$, $ \tilde r_j \le 1$, $j=1,2$, 
 \be\label{E:phi-sharp}
 \varphi_{f^\sharp}(x,\zeta)=e^{-i\wp(\frac {\tilde s}{ \tilde\sigma },\beta,X) }{f^\sharp}(x, \overline\zeta),
 \ee and decompose  
\begin{align} 
&  I_3(x,\lambda_1)-I_3(x, \lambda_2)\label{E:linear-norm-dec}\\
=& -\frac {\tilde\lambda_1-\tilde\lambda_2}{4\pi i} \iint_{\tilde s\le 2}\widetilde\gamma_1({\tilde\zeta})\frac{ \varphi_{f^\sharp}(\frac {\tilde s}{  \tilde\sigma},\beta,X)
-\varphi_{f^\sharp}(\frac {\tilde r_1}{ \tilde\sigma },\alpha_1,X)
 }{({\tilde\zeta}-\tilde\lambda_1)({\tilde\zeta}-\tilde\lambda_2)} d\bar{\tilde\zeta}\wedge d{\tilde\zeta}\nonumber\\
- &  \frac {\tilde\lambda_1-\tilde\lambda_2}{4\pi i} \iint_{\tilde s\le 2}\widetilde\gamma_1({\tilde\zeta})\frac{ \varphi_{f^\sharp}(\frac {\tilde s}{  \tilde\sigma},\beta,X)
-\varphi_{f^\sharp}(\frac {\tilde r_2}{ \tilde\sigma},\alpha_2,X)
 }{({\tilde\zeta}-\tilde\lambda_1)({\tilde\zeta}-\tilde\lambda_2)} d\bar{\tilde\zeta}\wedge d{\tilde\zeta}\nonumber\\
+&\frac {\varphi_{f^\sharp}(\frac {\tilde r_1}{ \tilde\sigma },\alpha_1,X)}{4\pi i} \iint_{\tilde s\le 2}\widetilde\gamma_1({\tilde\zeta})[\frac{ 1
}{{\tilde\zeta}-\tilde\lambda_2}-\frac{ 1
}{{\tilde\zeta}-\tilde\lambda_1}] d\bar{\tilde\zeta}\wedge d{\tilde\zeta}\nonumber\\
+&\frac {\varphi_{f^\sharp}(\frac {\tilde r_2}{ \tilde\sigma },\alpha_2,X)}{4\pi i} \iint_{\tilde s\le 2}\widetilde\gamma_1({\tilde\zeta})[\frac{ 1
}{{\tilde\zeta}-\tilde\lambda_2}-\frac{ 1
}{{\tilde\zeta}-\tilde\lambda_1}] d\bar{\tilde\zeta}\wedge d{\tilde\zeta}.\nonumber
\end{align}

In view of  $  f^\sharp\in C^{\mu }_{\tilde\sigma}(D_{\kappa_1,\corr{\frac{1}{\tilde\sigma}}})$ and  $  f^\sharp(x,\kappa_1)=0$, we have
\be\label{E:mu-radial}
|\varphi_{f^\sharp}(\frac {\tilde r}{ \tilde\sigma },\alpha,X)|_{L^\infty(D_{\kappa_1 })}\le C |f^\sharp|_{H^\mu_{\tilde\sigma}(D_{\kappa_1,\corr{\frac{1}{\tilde\sigma}} })}\tilde r^\mu .
\ee Along with Lemma \ref{L:stokes}, yields
\begin{align*}
|\frac {\varphi_{f^\sharp}(\frac {\tilde r_1}{ \tilde\sigma},\alpha_1,X)}{4\pi i}  \iint_{\tilde s\le 2}\widetilde\gamma_1({\tilde\zeta})[\frac{ 1
}{{\tilde\zeta}-\tilde\lambda_2}-\frac{ 1
}{{\tilde\zeta}-\tilde\lambda_1}] d\bar{\tilde\zeta}\wedge d{\tilde\zeta}| 
\le & C\epsilon_0 |f^\sharp|_{H^\mu_{\tilde\sigma}(D_{\kappa_1,\corr{\frac{1}{\tilde\sigma}} })}  |\tilde\lambda_1-\tilde\lambda_2|^\mu,   \tilde r_1=\tilde r_2,\\
|\frac {\varphi_{f^\sharp}(\frac {\tilde r_1}{ \tilde\sigma},\alpha_1,X)}{4\pi i} \iint_{\tilde s\le 2}\widetilde\gamma_1({\tilde\zeta})[\frac{ 1
}{{\tilde\zeta}-\tilde\lambda_2}-\frac{ 1
}{{\tilde\zeta}-\tilde\lambda_1}] d\bar{\tilde\zeta}\wedge d{\tilde\zeta}| 
\le & 0,\qquad  \alpha_1=\alpha_2
\end{align*}respectively. Therefore, 
\be \label{E:I_1-mu-3}
|\frac {\varphi _{f^\sharp}(\frac {\tilde r_1}{ \tilde\sigma},\alpha_1,X)}{4\pi i} \iint_{\tilde s\le 2}\widetilde\gamma_1({\tilde\zeta})[\frac{ 1
}{{\tilde\zeta}-\tilde\lambda_2}-\frac{ 1
}{{\tilde\zeta}-\tilde\lambda_1}] d\bar{\tilde\zeta}\wedge d{\tilde\zeta}|
\le  C\epsilon_0 |f^\sharp|_{C^{\mu }_{\tilde\sigma}(D_{\kappa_1,\corr{\frac{1}{\tilde\sigma}}})}  |\tilde\lambda_1-\tilde\lambda_2|^\mu.
\ee

In an entirely similar way, 
\be \label{E:I_1-mu-4}
|\frac {\varphi_{f^\sharp}(\frac {\tilde r_2}{ \tilde\sigma},\alpha_2,X)}{4\pi i} \iint_{\tilde s\le 2}\widetilde\gamma_1({\tilde\zeta})[\frac{ 1
}{{\tilde\zeta}-\tilde\lambda_2}-\frac{ 1
}{{\tilde\zeta}-\tilde\lambda_1}] d\bar{\tilde\zeta}\wedge d{\tilde\zeta}| 
\le  C\epsilon_0 |f^\sharp |_{C^{\mu }_{\tilde\sigma}(D_{\kappa_1,\corr{\frac{1}{\tilde\sigma}}})}  |\tilde\lambda_1-\tilde\lambda_2|^\mu.
\ee

Let us now investigate the first term on the right hand side of \eqref{E:linear-norm-dec}.  Applying Lemma \ref{L:vekua}, it suffices to derive the estimate  for all $\lambda_1$, $\lambda_2$ with $\tilde D\subset \{\tilde s\le 2\}  $ being a disk centred at $\tilde\lambda_1$ with radius $l$ and $l=2|\tilde\lambda_2-\tilde\lambda_1|$ (cf. \cite[5.1]{Ga66}). Write
\begin{align}
&-\frac {\tilde\lambda_1-\tilde\lambda_2}{4\pi i}\iint_{\tilde s\le 2}\widetilde\gamma_1({\tilde\zeta})\frac{ \varphi _{f^\sharp}(\frac {\tilde s}{  \tilde\sigma},\beta,X)
-\varphi _{f^\sharp}(\frac {\tilde r_1}{ \tilde\sigma},\alpha_1,X)
 }{({\tilde\zeta}-\tilde\lambda_1)({\tilde\zeta}-\tilde\lambda_2)} d\bar{\tilde\zeta}\wedge d{\tilde\zeta}\label{E:holder-ext}\\
 =&-\frac {\tilde\lambda_1-\tilde\lambda_2}{4\pi i}  \iint_{\tilde D}\widetilde\gamma_1({\tilde\zeta})\frac{ \varphi _{f^\sharp}(\frac {\tilde s}{  \tilde\sigma },\beta,X)
-\varphi _{f^\sharp}(\frac {\tilde r_1}{ \tilde\sigma},\alpha_1,X)
 }{({\tilde\zeta}-\tilde\lambda_1)({\tilde\zeta}-\tilde\lambda_2)} d\bar{\tilde\zeta}\wedge d{\tilde\zeta}\nonumber\\
 -& \frac {\tilde\lambda_1-\tilde\lambda_2}{4\pi i}\iint_{\{\tilde s\le 2\}/\tilde D}\widetilde\gamma_1({\tilde\zeta})\frac{ \varphi _{f^\sharp}(\frac {\tilde s}{ \tilde\sigma},\beta,X)
-\varphi _{f^\sharp}(\frac {\tilde r_1}{ \tilde\sigma },\alpha_1,X)
 }{({\tilde\zeta}-\tilde\lambda_1)({\tilde\zeta}-\tilde\lambda_2)} d\bar{\tilde\zeta}\wedge d{\tilde\zeta}.\nonumber
\end{align}

Let $\tilde D_0=\{\zeta: |\zeta-\tilde \lambda_1|<\frac{3l}2\} $.  
 \begin{itemize}
 \item  If ${\tilde\zeta}\in \{\tilde s\le 2\}/\tilde D$ and   $\kappa_1\in \tilde D_0$, then
 \[
  \frac 1C\le|\frac{{\tilde\zeta}-\tilde\lambda_1}{{\tilde\zeta}-\tilde\lambda_2}|,|\frac{{\tilde\zeta}-\kappa_1}{{\tilde\zeta}-\tilde\lambda_1}|,|\frac{{\tilde\zeta}-\kappa_1}{{\tilde\zeta}-\tilde\lambda_2}|\le C.
 \] In this case, using $  f^\sharp\in C^{\mu }_{\tilde\sigma}(D_{\kappa_1,\corr{\frac{1}{\tilde\sigma}}})$  and  {\cite[Chapter 1,\S6.1]{V62}}, 
 \begin{align}
&|-\frac {\tilde\lambda_1-\tilde\lambda_2}{4\pi i}\iint_{\{\tilde s\le 2\}/\tilde D}\widetilde\gamma_1({\tilde\zeta})\frac{ \varphi _{f^\sharp}(\frac {\tilde s}{  \tilde\sigma },\beta,X)
-\varphi _{f^\sharp}(\frac {\tilde r_1}{  \tilde\sigma },\alpha_1,X)
 }{({\tilde\zeta}-\tilde\lambda_1)({\tilde\zeta}-\tilde\lambda_2)} d\bar{\tilde\zeta}\wedge d{\tilde\zeta}|\label{E:holder-ext-1}\\
 \le &C\epsilon_0|f^\sharp|_{C^{\mu }_{\tilde\sigma}(D_{\kappa_1,\corr{\frac{1}{\tilde\sigma}}})}|\tilde\lambda_1-\tilde\lambda_2|\iint_{\{\tilde s\le 2\}/\tilde D} \frac{1}{|\tilde\zeta-\lambda_2| |{\tilde\zeta}-\tilde\lambda_1|^{2-\mu} } d\bar{\tilde\zeta}\wedge d{\tilde\zeta} \nonumber\\
  \le & C\epsilon_0|f^\sharp|_{C^{\mu }_{\tilde\sigma}(D_{\kappa_1,\corr{\frac{1}{\tilde\sigma}}})}|\tilde\lambda_1-\tilde\lambda_2|^\mu.\nonumber
\end{align}
 \item If ${\tilde\zeta}\in \{\tilde s\le 2\}/\tilde D$ and $\kappa_1\notin \tilde D_0$  then
 \[
 \begin{gathered}
 \frac 1C\le|\frac{{\tilde\zeta}-\tilde\lambda_1}{{\tilde\zeta}-\tilde\lambda_2}|\le C,\quad
 |\tilde\lambda_1 -\tilde\lambda_2 |\le \frac 1C\min\{| {\tilde\lambda_1}-\kappa_1 |,|{\tilde\lambda_2}-\kappa_1 |\}. 
 \end{gathered}
 \]
In this case, using $  f^\sharp\in C^{\mu }_{\tilde\sigma}(D_{\kappa_1,\corr{\frac{1}{\tilde\sigma}}})$  and  {\cite[Chapter 1,\S6.1]{V62}},  
\begin{align}
&|-\frac {\tilde\lambda_1-\tilde\lambda_2}{4\pi i}\iint_{\{\tilde s\le 2\}/\tilde D}\widetilde\gamma_1({\tilde\zeta})\frac{ \varphi _{f^\sharp}(\frac {\tilde s}{  \tilde\sigma },\beta,X)
-\varphi _{f^\sharp}(\frac {\tilde r_1}{  \tilde\sigma },\alpha_1,X)
 }{({\tilde\zeta}-\tilde\lambda_1)({\tilde\zeta}-\tilde\lambda_2)} d\bar{\tilde\zeta}\wedge d{\tilde\zeta}|\label{E:holder-ext-2}\\
 \le & C\epsilon_0|f^\sharp|_{C^{\mu }_{\tilde\sigma}(D_{\kappa_1,\corr{\frac{1}{\tilde\sigma}}})}|\tilde\lambda_1-\tilde\lambda_2|\iint_{\{\tilde s\le 2\}/\tilde D} \frac{ 1
 }{|\tilde\zeta-\kappa_1| |{\tilde\zeta}-\tilde\lambda_1|^{2-\mu} } d\bar{\tilde\zeta}\wedge d{\tilde\zeta} \nonumber\\
  \le & C\epsilon_0|f^\sharp|_{C^{\mu }_{\tilde\sigma}(D_{\kappa_1,\corr{\frac{1}{\tilde\sigma}}})}|\tilde\lambda_1-\tilde\lambda_2|^\mu.\nonumber
\end{align}
 \end{itemize}
Therefore the second term on the RHS of \eqref{E:holder-ext} is done.

Let $\tilde L(\zeta)=0$ be the line perpendicular to $\overline{\lambda_1\lambda_2}$ and passing through $\frac 12(\lambda_1+\lambda_2)$.  Set
\[
\begin{split}
\tilde D_{\tilde\lambda_1,\pm}=\tilde D\cap\{ \zeta: L(\zeta)L(\lambda_1)\gtrless 0\}.  
\end{split}
\]Therefore, thanks to $ f^\sharp \in C^\mu_{\tilde\sigma}(D_{\kappa_1,\frac{1}{\tilde\sigma}})$,   setting $\eta=\frac{\tilde\zeta-\tilde\lambda_1}{|\tilde\lambda_1-\lambda_2|}$,  $\frac{\tilde\zeta-\kappa_1}{|\tilde\lambda_1-\lambda_2|}=\eta-r_0e^{i\alpha_0}$,  and using   {\cite[Chapter 1,\S6.1]{V62}},   
\begin{align}
&|\frac {\tilde\lambda_1-\tilde\lambda_2}{4\pi i} \iint_{\tilde D_{\tilde\lambda_1,+}}\widetilde\gamma_1({\tilde\zeta})\frac{ \varphi _{f^\sharp}(\frac {\tilde s}{  \tilde\sigma},\beta,X)
-\varphi _{f^\sharp}(\frac {\tilde r_1}{ \tilde\sigma},\alpha_1,X)
 }{({\tilde\zeta}-\tilde\lambda_1)({\tilde\zeta}-\tilde\lambda_2)} d\bar{\tilde\zeta}\wedge d{\tilde\zeta}|\label{E:lambda-1-p}\\
 \le &C\epsilon_0|\tilde\lambda_1-\tilde\lambda_2||f^\sharp|_{C^\mu_{\tilde\sigma}(D_{\kappa_1,\corr{\frac{1}{\tilde\sigma}}})}  |\iint_{\tilde D_{\tilde\lambda_1,+}}  \frac 1{|\tilde\zeta-\kappa_1||{\tilde\zeta}-\tilde\lambda_1|^{1-\mu}|{\tilde\zeta}-\tilde\lambda_2|}d\bar{\tilde\zeta}\wedge d{\tilde\zeta}|\nonumber\\
 \le &C\epsilon_0|\tilde\lambda_1-\tilde\lambda_2|^\mu|f^\sharp|_{C^\mu_{\tilde\sigma}(D_{\kappa_1,\corr{\frac{1}{\tilde\sigma}}})}  |\iint_{\{|\eta|\le 2\}\cap \tilde D_{\tilde\lambda_1,+}}  \frac 1{|\eta-r_0e^{i\alpha_0}||\eta|^{1-\mu} |\eta-e^{i\alpha'}| } d \eta_Rd\eta_I|\nonumber\\
 \le &C\epsilon_0|\tilde\lambda_1-\tilde\lambda_2|^\mu|f^\sharp|_{C^\mu_{\tilde\sigma}(D_{\kappa_1,\corr{\frac{1}{\tilde\sigma}}})}  |\iint_{\{|\eta|\le 2\}\cap \tilde D_{\tilde\lambda_1,+}}  \frac 1{|\eta-r_0e^{i\alpha_0}||\eta|^{1-\mu}  } d \eta_Rd\eta_I|\nonumber\\
\le &C\epsilon_0|f^\sharp|_{C^\mu_{\tilde\sigma}(D_{\kappa_1,\corr{\frac{1}{\tilde\sigma}}})}|\tilde\lambda_1-\tilde\lambda_2|^\mu. \nonumber
\end{align} 

By analogy,
\begin{align*}
&|\frac {\tilde\lambda_1-\tilde\lambda_2}{4\pi i} \iint_{\tilde D_{\tilde\lambda_1,-}}\widetilde\gamma_1({\tilde\zeta})\frac{ \varphi _{f^\sharp}(\frac {\tilde s}{  \tilde\sigma},\beta,X)
-\varphi _{f^\sharp}(\frac {\tilde r_1}{ \tilde\sigma},\alpha_1,X)
 }{({\tilde\zeta}-\tilde\lambda_1)({\tilde\zeta}-\tilde\lambda_2)} d\bar{\tilde\zeta}\wedge d{\tilde\zeta}|\\
 \le &|\frac {\tilde\lambda_1-\tilde\lambda_2}{4\pi i} \iint_{\tilde D_{\tilde\lambda_1,-}}\widetilde\gamma_1({\tilde\zeta})\frac{ \varphi _{f^\sharp}(\frac {\tilde s}{  \tilde\sigma},\beta,X)
-\varphi _{f^\sharp}(\frac {\tilde r_2}{ \tilde\sigma},\alpha_2,X)
 }{({\tilde\zeta}-\tilde\lambda_1)({\tilde\zeta}-\tilde\lambda_2)} d\bar{\tilde\zeta}\wedge d{\tilde\zeta}|\\
 +&|\frac {\tilde\lambda_1-\tilde\lambda_2}{4\pi i} \iint_{\tilde D_{\tilde\lambda_1,-}}\widetilde\gamma_1({\tilde\zeta})\frac{ \varphi _{f^\sharp}(\frac {\tilde r_2}{  \tilde\sigma},\alpha_2,X)
-\varphi _{f^\sharp}(\frac {\tilde r_1}{ \tilde\sigma},\alpha_1,X)
 }{({\tilde\zeta}-\tilde\lambda_1)({\tilde\zeta}-\tilde\lambda_2)} d\bar{\tilde\zeta}\wedge d{\tilde\zeta}|\\
  \le &C\epsilon_0|f^\sharp|_{C^\mu_{\tilde\sigma}(D_{\kappa_1,\corr{\frac{1}{\tilde\sigma}}})}|\tilde\lambda_1-\tilde\lambda_2|^\mu\\
 +&|\frac {\tilde\lambda_1-\tilde\lambda_2}{4\pi i} \iint_{\tilde D_{\tilde\lambda_1,-}}\widetilde\gamma_1({\tilde\zeta})\frac{ \varphi _{f^\sharp}(\frac {\tilde r_2}{  \tilde\sigma},\alpha_2,X)
-\varphi _{f^\sharp}(\frac {\tilde r_1}{ \tilde\sigma},\alpha_1,X)
 }{({\tilde\zeta}-\tilde\lambda_1)({\tilde\zeta}-\tilde\lambda_2)} d\bar{\tilde\zeta}\wedge d{\tilde\zeta}|.
\end{align*}
Applying $ f^\sharp \in C^\mu_{\tilde\sigma}(D_{\kappa_1,\corr{\frac{1}{\tilde\sigma}}})$,  Stokes' theorem,   and  
 $| {\tilde\zeta-\tilde\lambda_1}| $, $ | {\tilde\zeta-\tilde\lambda_2}|\sim| {\tilde\lambda_1}-\tilde\lambda_2|$ on the boundary of $\tilde D_{\tilde\lambda_1,-}$ (assured by $|\tilde\lambda|\le 1$, $\tilde D\subset\{\tilde s<2\}$), 
\begin{align*}
&|\frac {\tilde\lambda_1-\tilde\lambda_2}{4\pi i} \iint_{\tilde D_{\tilde\lambda_1,-}}\widetilde\gamma_1({\tilde\zeta})\frac{ \varphi _{f^\sharp}(\frac {\tilde r_2}{  \tilde\sigma},\alpha_2,X)
-\varphi _{f^\sharp}(\frac {\tilde r_1}{ \tilde\sigma},\alpha_1,X)
 }{({\tilde\zeta}-\tilde\lambda_1)({\tilde\zeta}-\tilde\lambda_2)} d\bar{\tilde\zeta}\wedge d{\tilde\zeta}|\\
 \le &C|f^\sharp|_{C^\mu_{\tilde\sigma}(D_{\kappa_1,\corr{\frac{1}{\tilde\sigma}}})}|\tilde\lambda_1-\tilde\lambda_2|^{1+\mu} \iint_{\tilde D_{\tilde\lambda_1,-}}\widetilde\gamma_1({\tilde\zeta})\frac{ 1 }{({\tilde\zeta}-\tilde\lambda_1)({\tilde\zeta}-\tilde\lambda_2)} d\bar{\tilde\zeta}\wedge d{\tilde\zeta}|\\
 =&C|f^\sharp|_{C^\mu_{\tilde\sigma}(D_{\kappa_1,\corr{\frac{1}{\tilde\sigma}}})}|\tilde\lambda_1-\tilde\lambda_2|^{1+\mu} \iint_{\tilde D_{\tilde\lambda_1,-}} \frac{ \partial_{\overline\zeta}\left[\ln(1-\gamma|\beta|) \frac 1{ {\tilde\zeta}-\tilde\lambda_1}\right]}{ {\tilde\zeta}-\tilde\lambda_2  } d\bar{\tilde\zeta}\wedge d{\tilde\zeta}|\\
\le &C\epsilon_0|f^\sharp|_{C^\mu_{\tilde\sigma}(D_{\kappa_1,\corr{\frac{1}{\tilde\sigma}}})}|\tilde\lambda_1-\tilde\lambda_2|^\mu.
\end{align*}  
Therefore the first term on the RHS of \eqref{E:holder-ext} is done. Thus  
\begin{multline}\label{E:I_1-mu-1}
|\frac {\tilde\lambda_1-\tilde\lambda_2}{4\pi i}\iint_{\tilde s\le 2}\widetilde\gamma_1({\tilde\zeta})\frac{ \varphi _{f^\sharp}(\frac {\tilde s}{  \tilde\sigma },\beta,X)
-\varphi _{f^\sharp}(\frac {\tilde r_1}{ \tilde\sigma },\alpha_1,X)
 }{({\tilde\zeta}-\tilde\lambda_1)({\tilde\zeta}-\tilde\lambda_2)} d\bar{\tilde\zeta}\wedge d{\tilde\zeta}|\\
 \le C\epsilon_0|f^\sharp|_{C^\mu_{\tilde\sigma}(D_{\kappa_1,\corr{\frac{1}{\tilde\sigma}}})}|\tilde\lambda_1-\tilde\lambda_2|^\mu,\quad\textit{for $|\tilde\lambda_j-\kappa_1|\le 1,\ j=1,2$. }
\end{multline}

In an entirely similar way, 
\begin{multline}\label{E:I_1-mu-2}
|\frac {\tilde\lambda_1-\tilde\lambda_2}{4\pi i}\iint_{\tilde s\le 2}\widetilde\gamma_1({\tilde\zeta})\frac{ \varphi _{f^\sharp}(\frac {\tilde s}{  \tilde\sigma },\beta,X)
-\varphi _{f^\sharp}(\frac {\tilde r_2}{ \tilde\sigma },\alpha_2,X)
 }{({\tilde\zeta}-\tilde\lambda_1)({\tilde\zeta}-\tilde\lambda_2)} d\bar{\tilde\zeta}\wedge d{\tilde\zeta}| 
 \\
 \le C\epsilon_0|f^\sharp|_{C^\mu_{\tilde\sigma}(D_{\kappa_1,\corr{\frac{1}{\tilde\sigma}}})}|\tilde\lambda_1-\tilde\lambda_2|^\mu,\quad\textit{for $|\tilde\lambda_j-\kappa_1|\le 1,\ j=1,2$. }
\end{multline}  

Plugging \eqref{E:I_1-mu-3}, \eqref{E:I_1-mu-4}, \eqref{E:I_1-mu-1}, and \eqref{E:I_1-mu-2} into \eqref{E:linear-norm-dec}, we obtain
\be 
\label{E:linear-norm-est}
|I_3(x, \lambda_1 )-I_3(x, \lambda_2)|\le C\epsilon_0 |f^\sharp|_{C^\mu_{\tilde\sigma}(D_{\kappa_1,\corr{\frac{1}{\tilde\sigma}}})} |\tilde\lambda_1-\tilde\lambda_2|^\mu 
\ee for $|\tilde\lambda_j-\kappa_1|\le 1,\ j=1,2$. 
Hence
\be \label{E:mean-vekua-3}
 | I_3|_{C^\mu_{\tilde\sigma}(D_{\kappa_1,\corr{\frac{1}{\tilde\sigma}}})}\le C\epsilon_0    | f^\sharp|_{C^\mu_{\tilde\sigma}(D_{\kappa_1,\corr{\frac{1}{\tilde\sigma}} })} .
 \ee

\end{proof}

\subsubsection{Estimates for $I_4$ and $I_5$}\label{SSS:FN-homo}   \hfill \\

We first deal with the generic cubic case $X_3\ne0$. Without loss of generality and for simplicity, we first  assume   
\be\label{E:assumption}
\begin{gathered}
\kappa_j=\kappa_1,\quad |\lambda-\kappa_1|\le \frac \delta 2,\quad  X_3> 0 , X_2\ge 0, X_1\ge 0.
\end{gathered}
\ee Degenerate cases will be proved at the end of  this subsection. 

For $ \widehat \sigma\in\{X_1,\sqrt{X_2},\sqrt[3]{X_3}\}$,  consider the  $\widehat\sigma$-scaled coordinates
\[
\zeta =\kappa_1+se^{i\beta}=\kappa_1+\frac{\widehat s }{\widehat \sigma}e^{i\beta}\quad\in D_{\kappa_1},
\] and the deformation 
\[
\widehat s\in\RR \longrightarrow \widehat se^{i\tau}\in\CC. 
\]
   Observe that
\[
\begin{split}
&\mathfrak {Re} ({-i\wp(\frac {\widehat se^{i\tau}}{\widehat \sigma},\beta,X)})=\frac{X_3}{\widehat \sigma ^3} \sin3\tau\sin3\beta\widehat s^3+  \frac{X_2}{\widehat \sigma^2} \sin2\tau\sin2\beta\widehat s^2 +\frac{X_1}{\widehat \sigma} \sin \tau\sin \beta\widehat s,\\
&\partial_{\widehat s}\wp(\frac{\widehat s}{\widehat \sigma},\beta,X)=3\frac{X_3}{\widehat \sigma^3}  \sin3\beta\widehat s^2+ 2 \frac{X_2}{\widehat \sigma^2}  \sin2\beta\widehat s  +\frac{X_1}{\widehat \sigma } \sin \beta.  
\end{split}
\] Hence as $|\tau|\ll 1$, non zero roots of $\partial_{\widehat s}\wp(\frac{\widehat s}{\widehat \sigma},\beta,X)$, the major obstruction of \eqref{E:deform-intro-1}, are approximated  by those of  $\mathfrak {Re} ({-i\wp(\frac {\widehat se^{i\tau_\dagger}}{\widehat \sigma},\beta,X)})/\widehat{s}$.

\begin{definition}\label{E:stationary} The stationary points are defined to be 
\be\label{E:determinant}
\widehat s_\pm = \frac{    - 1\pm\sqrt{1-\Delta}\, }{3\frac{X_3 }{\widehat\sigma X_2}\frac{\sin3\beta}{\sin2\beta}},\quad \Delta= 3\frac{ {X_1}{X_ 3}}{X_2^2}\frac{\sin\beta\sin3\beta}{\sin^22\beta},
\ee
  satisfying 
\[
\partial_{\widehat s}\wp(\frac{\widehat s_\pm}{\widehat \sigma},\beta,X)=3\frac{X_3}{\widehat \sigma^2} \sin3\beta \widehat s^2_\pm+2\frac{X_2}{\widehat \sigma^{2 }}\sin2\beta\widehat s_\pm+\frac{X_1}{\widehat \sigma }\sin\beta=0.
\]

Denote 
\be\label{E:Omega}
\ba{lll}
\Omega_1=\{0\le|\beta|\le\frac\pi 3 \},& \Omega_2=\{\frac {\pi} 3\le|\beta|\le\frac {\pi }2\},& \Omega_3=\{\frac {\pi} 2\le|\beta|\le\frac {2\pi }3\},\\
\Omega_4=\{\frac {2\pi} 3\le|\beta|\le\pi\},&  &
\ea
\ee one has Figure \ref{Fg:signature} for signatures of $sin (k\beta)$ on $\Omega_j$.  Moreover, according to the determinant $\Delta$, we have Table \ref{Tb:dynamic-FNH} for properties  of roots $\widehat s_\pm$. Notice that
\be\label{E:roots-distance}
  |\widehat s_ +-\widehat s_ -  
|=    |\frac{\sqrt{1-\Delta}}{3{ \frac{X_3}{\widehat\sigma X_2 } \frac{ \sin3\beta  }{  \sin2\beta} }}|.
\ee

{\begin{center}
\begin{figure} 
\begin{tikzpicture}[xscale=1.2 , yscale=.6]
\draw [  thick, ->](0,0) -- (3,0);
\draw [ thick, ->](0,0) -- (1.5,2.59) node [right]{$\frac{ \pi} 3$};
\draw [ thick, ->](0,0) -- (1.5,-2.59) node [right]{$-\frac{ \pi} 3$};
\draw [ thick,    ->](0,0) -- (0,3) ;
\draw [ thick,    ->](0,0) -- (0,-3) ;
\draw [ thick, ->](0,0) -- (-1.5,2.59) node [left]{$\frac{2\pi} 3$};
\draw [ thick, ->](0,0) -- (-1.5,-2.59) node [left]{$-\frac{2\pi} 3$};
\draw [ thick,   ->](0,0) -- (-3,0);
\node   at (3.5,1.6) {$(+,+,+)$};
\node   at (3.5,-1.6) {$(-,-,-)$};
\node  at (1,3.5) {$(+,+,-)$};\node  at (-1,3.5) {$(+,-,-)$};
\node  at (1,-3.5) {$(-,-,+)$};\node  at (-1,-3.5) {$(-,+,+)$};
\node   at (-3.5,1.6) {$(+,-,+)$};\node   at (-3.5,-1.6) {$(-,+,-)$};
\node [above] at (3.5,1.6) {$\Omega_1$};
\node [above] at (3.5,-1.6) {$\Omega_1$};
\node [above] at (1,3.5) {$\Omega_2$};
\node [above] at (1,-3.5) {$\Omega_2$};
\node [above] at (-1, 3.5) {$\Omega_3$};
\node [above] at (-1,-3.5) {$\Omega_3$};
\node [above] at (-3.5,1.6) {$\Omega_4$};
\node [above] at (-3.5,-1.6) {$\Omega_4$};
\end{tikzpicture}
\caption[Signatures]{\small Signatures of  $
( \sin\beta,\sin2\beta,\sin3\beta)$ for $X_1,X_2,X_3>0$}\label{Fg:signature}
\end{figure}
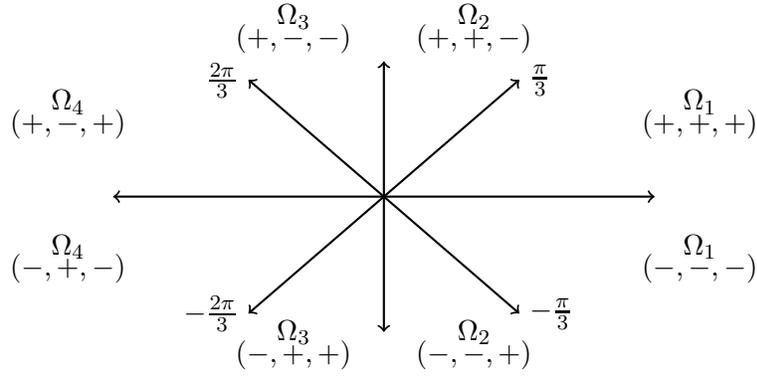
\end{center}}

\begin{table}
{\small\begin{center}
\vskip.1in
\begin{tabular}{|c|c|c|l|l|}
 \hline
 \bf{Type}&\bf{Subtype}&\bf{Range of $\Delta$}  &\bf{Properties of $\widehat s_\pm$}   &\bf{$\beta$-domain}\\
& &  && ({\bf order of } $\widehat s_\pm$)\\
   \hline\hline
 \multirow{2}*{$\mathfrak A(\beta,X)$}  & $\mathfrak A'$  &  {$  (2,\infty)$}   &\multirow{2}* {$\widehat s_\pm$ complex roots} & $\Omega_1(\frac{\widehat s_++\widehat s_-}{2}\le 0)$       \\
\cline{2-2}\cline{3-3}   
 & $\mathfrak A''$ & {$  (1,2)$}   &   & $ \Omega_4(\frac{\widehat s_++\widehat s_-}{2}\ge 0)$ \\
 \hline
 \multirow{2}*{$\mathfrak B(\beta,X)$}  & $\mathfrak B'$   &\multirow{2}* {$ { (\frac 12,1)}$}    &  {adjacent real roots w.  } & $ {\Omega_1 (\widehat s_-\le\widehat s_ +\le 0) }$       \\
\cline{2-2}\cline{5-5}  
 & $\mathfrak B''$ &   &   & $ {\Omega_4 (\widehat s_-\ge\widehat s_ +\ge 0)}$ \\
 \hline
 \multirow{2}*{$\mathfrak C(\beta,X)$}  & $\mathfrak C'$   &\multirow{2}* { $ (0, \frac { 1}2)$}    & {$\widehat s_ +=\frac{-\Delta}{6\frac{X _3}{\widehat \sigma X_2}\frac{\sin 3\beta}{\sin2\beta} } +\mbox{l.o.t.}, $} & $ {\Omega_1  (\widehat s_-\le\widehat s_+\le 0)}$    \\
\cline{2-2}\cline{5-5}  
 & $\mathfrak C''$ &   &$ {\widehat s_ -}= \frac{ -2 }{ 3\frac{X_3}{\widehat \sigma X_2}\frac{\sin3\beta}{\sin 2 \beta}  } +\mbox{l.o.t.}  $   & $ {\Omega_4   {(\widehat s_-\ge\widehat s_+\ge 0)}}$ \\
 \hline
 \multirow{2}*{$\mathfrak D(\beta,X)$}  &     &\multirow{2}* {$ (- \frac { 1}2,0)$}    &  {$\widehat s_ +=\frac{-\Delta}{6\frac{X _3}{\widehat \sigma X_2}\frac{\sin 3\beta}{\sin2\beta} } +\mbox{l.o.t.}, $   } & $ \Omega_2  (\widehat s_-\ge 0\ge\widehat s_ + )$     \\ 
 &   &   &$ {\widehat s_ -}= \frac{ -2 }{ 3\frac{X_3}{\widehat \sigma X_2}\frac{\sin3\beta}{\sin 2 \beta}  } +\mbox{l.o.t.}  $   & $ {\Omega_3   {(\widehat s_+\ge 0\ge \widehat s_ - )}}$ \\
 \hline
 \multirow{2}*{$\mathfrak E(\beta,X)$}  &     &\multirow{2}* {$  (-\infty,- \frac { 1}2)$}    &\multirow{2}* {$\widehat s_\pm$ real roots} & $ \Omega_2  (\widehat s_-\ge 0\ge\widehat s_ + )$      \\  
 &   &   &   & $ {\Omega_3  (\widehat s_+\ge 0\ge \widehat s_ - )}$ \\
 \hline
\end{tabular}
\end{center} }

\caption{\small Properties of $\widehat s_\pm$ and $\Delta$ for  Type $\mathfrak A,\cdots,\mathfrak E$ when  $X_1>0,X_2,X_3\ge 0$} 
\label{Tb:dynamic-FNH}
\end{table}

\end{definition}
In the following, we will  define essential stationary points $\widehat s_{j,\ast}$ and decompose  $[0,\widehat\sigma]$ into intervals $\mho_j$ around $\widehat s_{j,\ast}$. The deformation will be defined on $\mho_j $.  
\begin{definition}\label{D:deformation}   We define 
the essential stationary points $\widehat s_{j,\ast}=\widehat s_{j,\ast}(\widehat\sigma,\beta, X)$ by
\begin{align}
\widehat s_{0,\ast} = &  0, \nonumber\\
\widehat s_{1,\ast}= &  \left\{
{\ba{ll}  
 { \frac{\widehat s_{+ }+\widehat s_{- }}2\gtrless 0}  , &       Type\,\mathfrak A' \wedge(\widehat{\sigma}\in\{\sqrt{X_2},\sqrt[3]{X_3}\}), \,\mathfrak A'' , \\
  {  \inf\widehat s_{\pm }}>0 , &     Type\, \mathfrak B'',\mathfrak C'',\\
      \sup\widehat s_{\pm }>0 ,&    Type\,  \mathfrak D, \mathfrak  E,  \\
     - ,&  Type\,  \mathfrak A '\wedge(\widehat{\sigma}=X_1),\, \mathfrak B' ,\mathfrak C' ,  
 \ea}\right.  \ &\nonumber\\ 
  \widehat s_{2,\ast}= &  \left\{
{\ba{ll}  
  { \sup\widehat s_{\pm }}, &   \hskip.3in     Type\, \mathfrak B'',\mathfrak C'' ,\\
 -  ,&  \hskip.3in    \textit{others, }
 \ea}\right.\label{E:F3-ast} 
\end{align}   where $-$  means no definition.    Given $0<\epsilon_1<\frac{\pi}{2k}\ll 1$, define neighborhood $   \mho_j(\widehat\sigma,\beta,X )$ of essential critical points $\widehat s_{j,\ast}$ by 
\begin{align}
\mho_0 
=&\left\{
{\ba{l }  
{  [0, \frac 12]},    
\hskip3.15 in        Type\, \mathfrak A'\wedge(\widehat{\sigma}\in\{\sqrt{X_2},\sqrt[3]{X_3}\}),\mathfrak A'', \\
  {  [0, \frac 1{ 2\cos\epsilon_1 } \widehat s _{1,\ast} ],}   \hskip2.65 in   Type\, \mathfrak B'',\mathfrak C'',   \mathfrak D ,  \mathfrak E,\\
    {[0, \widehat\sigma\delta]=[0,\frac 12]\cup[\frac 12,\widehat\sigma\delta]=\mho_{0,<}\cup\mho_{0,>}} ,  
\hskip.9 in   Type\,\mathfrak A '\wedge(\widehat{\sigma}=X_1), \, \mathfrak B',\mathfrak C', 
 \ea}\right.\label{E:f-beta-interval}\\
 \mho_1  
=&
\left\{
{\ba{l}  
{  [\frac 12,\widehat s_{1,\ast}]\cup[\widehat s_{1,\ast}, \widehat\sigma\delta] }\equiv\mho_{1,<}\cup\mho_{1,>}, 
\hskip .1in Type\, \{[\mathfrak A'\wedge(\widehat{\sigma}\in      \{\sqrt{X_2},\sqrt[3]{X_3}\})]\vee\mathfrak A''\}\wedge (\widehat{s}_{1,\ast}>0),\\
{  [\frac 12,  \widehat\sigma\delta] }\equiv \mho_{1,>}, 
\hskip 1.4in  Type\, \{[\mathfrak A'\wedge(\widehat{\sigma}\in\{\sqrt{X_2},      \sqrt[3]{X_3}\})]\vee \mathfrak A''\}\wedge(\widehat{s}_{1,\ast}<0),\\
  { [ (1-  \frac 1{2\cos\epsilon_1 } ) \widehat s_{1,\ast}, \widehat s_{1,\ast} ]\cup [ \widehat s_{1,\ast}, \widehat s_{1,\ast}+ \frac {\widehat s_{2,\ast}-\widehat s_{1,\ast}  }{2\cos\epsilon_1  }  ]\equiv\mho_{1,<}\cup\mho_{1,>},}      
 \hskip .2in  Type\, \mathfrak B'',\mathfrak C'',\\
 { [ (1-  \frac 1{2\cos\epsilon_1 } )\widehat s_{1,\ast}, \widehat s_{1,\ast} ]\cup [  \widehat s_{1,\ast},\widehat\sigma\delta   ] \equiv\mho_{1,<}\cup\mho_{1,>},}      
\hskip 1in  Type\,  \mathfrak D ,  \mathfrak E,\\
  {\phi,}      
\hskip 3. in  Type\,\mathfrak A '\wedge(\widehat{\sigma}=X_1),\,  \mathfrak B',\mathfrak C', 
 \ea}\right.\nonumber\\
 \mho_2  
=&
\left\{
{\ba{l} 
  { [  \widehat s_{2,\ast}-\frac {\widehat s_{2,\ast}-\widehat s_{2,\ast}  }{2 \cos\epsilon_1 } , \widehat s_{2,\ast} ]\cup [  \widehat s_{2,\ast},\widehat\sigma\delta ] \equiv\mho_{2,<}\cup\mho_{2,>},}   \hskip.2in Type\, \mathfrak B'',\mathfrak C'',\\
 \phi      \hskip.2in   \textit{otherwise},
 \ea}\right.
 \nonumber
\end{align}

Write
\be\label{E:FNH-lambda-sigma}
\begin{gathered}
\lambda=\kappa_1+\frac{\widehat re^{i\alpha} }{\widehat \sigma}=\kappa_1+\frac{ \widehat s_{j,\ast}e^{i\beta}+\widehat r_je^{i\alpha_j}}{\widehat \sigma}  ,\\
 \widehat r_j=\widehat r_j(\widehat\sigma,\beta,X,\lambda),\ \alpha_j=\alpha_j(\widehat\sigma,\beta,X,\lambda),\,j=0,1,2.
\end{gathered}
\ee 
We define the deformation   defined by
\be \label{E:f3-deform-s}
\begin{gathered}
\zeta=\kappa_1+se^{i\beta}=\kappa_1+\frac{\widehat s e^{i\beta}}{\widehat\sigma}\\
 \widehat  s \mapsto \xi_j\equiv  \widehat s_{j,\ast} +\widehat s_je^{i\tau_j}  ,   \\
  \widehat s\equiv\widehat s_{j,\ast}\pm\widehat s_j\in \mho_j ,\,  |\tau_j|\lessgtr\frac\pi 2,\,\widehat s_j\ge 0,\, j=0,1,2, 
\end{gathered} 
\ee
with 
\begin{center}
\begin{tabular}{lll}
    \multirow{12}{*}{$   
 \left\{
{\ba{l}  
{  \mp\pi\lessgtr\tau_1\lessgtr \mp\pi\pm\epsilon_1 ,} \\
{ \mp\epsilon_1\lessgtr \tau_1\lessgtr 0  ,} \\
  {   \mp\pi\lessgtr \tau_1\lessgtr \mp\pi\pm\epsilon_1 ,} \\
   {   \pm\epsilon_1\gtrless\tau_1\gtrless 0 , }\\
   {  \pm\pi\gtrless\tau_1\gtrless\pm\pi\mp\epsilon_1,}\\
   {  \mp\epsilon_1\lessgtr\tau_1\lessgtr 0,}\\   
   { \mp\pi\lessgtr\tau_1\lessgtr \mp\pi\pm\frac{\epsilon_1 }4 ,} \\
   {  \mp\frac{\epsilon_1}4\lessgtr\tau_1\lessgtr 0 ,} \\
 {  \mp\pi\lessgtr \tau_1\lessgtr \mp\pi\pm\frac{\epsilon_1 }4,} \\
  {      \pm\frac{\epsilon_1 }4 \gtrless\tau_1\gtrless 0  ,} \\
   { \pm\pi\gtrless\tau_1\gtrless \pm\pi\mp\frac{\epsilon_1}4,}\\
   { \mp\frac{\epsilon_1}4\lessgtr\tau_1\lessgtr 0,}
 \ea}\right.$}
    &$\textit{for }\sin 3\beta\gtrless 0,\ \ ||\alpha_1-\beta|- \pi|\le \frac{\epsilon_1}2$,&$ \widehat s\in\mho_{1<}\ne\phi,\,  Type\, \mathfrak A ,$\\
    &$\textit{for }\sin 3\beta\gtrless0,\ \ |\alpha_1-\beta|\le \frac{\epsilon_1}2$,&$\widehat s\in\mho_{1>}\ne\phi, \,  Type\, \mathfrak A ,$\\ 
    &$\textit{for }\sin 3\beta\gtrless0,\ \ ||\alpha_1-\beta|- \pi|\le \frac{\epsilon_1}2$,&$\widehat s\in\mho_{1<}, $ $Type\, \mathfrak B'',\mathfrak C''$,\\
    &$\textit{for }\sin 3\beta\gtrless0,\ \ |\alpha_1-\beta|\le \frac{\epsilon_1}2$,&$ \widehat s\in\mho_{1>} ,$ $Type\, \mathfrak B'',\mathfrak C''$,\\
    &$\textit{for }\sin 3\beta\gtrless0,\ \ ||\alpha_1-\beta|- \pi|\le \frac{\epsilon_1}2$,&$ \widehat s\in\mho_{1<} ,$ $Type\, \mathfrak D,\mathfrak E$,\\
     &$\textit{for }\sin 3\beta\gtrless 0,\ \ |\alpha_1-\beta|\le \frac{\epsilon_1}2$,&$ \widehat s\in\mho_{1>} ,$ $Type\, \mathfrak D,\mathfrak E,\mathfrak B',\mathfrak C'$,\\
    &$ \textit{for }\sin 3\beta\gtrless 0,\ \ ||\alpha_1-\beta|- \pi|\ge \frac{\epsilon_1}2$,&$ \widehat s\in\mho_{1<}\ne\phi,\,  Type\, \mathfrak A ,$ \\
    &$ \textit{for }\sin 3\beta\gtrless 0,\ \ |\alpha_1-\beta|\ge \frac{\epsilon_1}2$,&$\widehat s\in\mho_{1>}\ne\phi,\,  Type\, \mathfrak A ,$ \\
    &$\textit{for }\sin 3\beta\gtrless 0,\ \ ||\alpha_1-\beta|- \pi|\ge \frac{\epsilon_1}2$,&$ \widehat s\in\mho_{1<}, $ $Type\, \mathfrak B'',\mathfrak C''$, \\
    &$\textit{for }\sin 3\beta\gtrless 0,\ \ |\alpha_1-\beta|\ge \frac{\epsilon_1}2$,&$\widehat s\in\mho_{1>}, $ $Type\, \mathfrak B'',\mathfrak C''$, \\
    &$\textit{for }\sin 3\beta\gtrless 0,\ \ ||\alpha_1-\beta|- \pi|\ge \frac{\epsilon_1}2$,&$ \widehat s\in\mho_{1<} ,$ $Type\, \mathfrak D,\mathfrak E$,\\
     &$\textit{for }\sin 3\beta\gtrless 0,\ \ |\alpha_1-\beta|\ge \frac{\epsilon_1}2$,&$ \widehat s\in\mho_{1>} ,$ $Type\, \mathfrak D, \mathfrak E,\mathfrak B',\mathfrak C'$,
\end{tabular}
\end{center}

\begin{center}
\begin{tabular}{lll}      
\multirow{8}{*}{$   
 \left\{
{\ba{l}  
  {\tau_0\equiv 0  ,} \\
  { \tau_0\equiv 0  ,} \\
   { \pm\epsilon_1\gtrless\tau_0\gtrless 0  ,} \\
  {\pm\frac{\epsilon_1 }4\gtrless\tau_0\gtrless 0  ,} \\
  { \mp\epsilon_1\lessgtr\tau_0\lessgtr 0  ,} \\
  {\mp\frac{\epsilon_1 }4\lessgtr\tau_0\lessgtr 0  ,}\\
  { \mp\epsilon_1\lessgtr\tau_0\lessgtr 0  ,} \\
  {\mp\frac{\epsilon_1 }4\lessgtr\tau_0\lessgtr 0  ,} 
  \ea}\right.$}&$\textit{for }\sin 3\beta\gtrless 0,\ \ |\alpha_0-\beta|\le \frac{\epsilon_1}2$,&$\widehat s\in  [0,\frac 12]\subset \mho_0$, $Type\ \mathfrak A,\mathfrak B',\mathfrak C'$,\\
    &$ \textit{for }\sin 3\beta\gtrless 0,\ \ |\alpha_0-\beta|\ge \frac{\epsilon_1}2$,&$\widehat s\in  [0,\frac 12]\subset \mho_0$, $Type\ \mathfrak A,\mathfrak B',\mathfrak C'$,\\
    &$\textit{for }\sin 3\beta\gtrless 0,\ \ |\alpha_0-\beta|\le \frac{\epsilon_1}2$,&$  \widehat s\in\mho_0$, $Type\ \mathfrak D,\mathfrak E$,\\
    &$ \textit{for }\sin 3\beta\gtrless 0,\ \ |\alpha_0-\beta|\ge \frac{\epsilon_1}2$,&$\widehat s\in\mho_0 $,  $Type\ \mathfrak D, \mathfrak E$,\\
    &$\textit{for }\sin 3\beta\gtrless 0,\ \ |\alpha_0-\beta|\le \frac{\epsilon_1}2$,&$\widehat s\in   \mho_{0,>}$, $Type\ \mathfrak A',\mathfrak B',\mathfrak C'$,\\
    &$ \textit{for }\sin 3\beta\gtrless 0,\ \ |\alpha_0-\beta|\ge \frac{\epsilon_1}2$,&$\widehat s\in   \mho_{0,>}$, $Type\ \mathfrak A',\mathfrak B',\mathfrak C'$,\\
    &$\textit{for }\sin 3\beta\gtrless 0,\ \ |\alpha_0-\beta|\le \frac{\epsilon_1}2$,&$  \widehat s\in\mho_0$, $Type\, \mathfrak B'',\mathfrak C''$,\\
    &$ \textit{for }\sin 3\beta\gtrless 0,\ \ |\alpha_0-\beta|\ge \frac{\epsilon_1}2$,&$\widehat s\in\mho_0$, $Type\, \mathfrak B'',\mathfrak C''$,\\
    \multirow{4}{*}{$   
 \left\{
{\ba{l}  
  {    \pm\pi \gtrless\tau_2\gtrless \pm\pi\mp\epsilon_1,}   \\
   {  \mp\epsilon_1\lessgtr\tau_2\lessgtr 0  , }  \\
  {    \pm\pi \gtrless\tau_2\gtrless \pm\pi\mp\frac{\epsilon_1 }4 ,}   \\
  { \mp\frac{\epsilon_1 }4\lessgtr\tau_2\lessgtr 0 ,}  
 \ea}\right.$}&$\textit{for }\sin 3\beta\gtrless 0,\ \ ||\alpha_2-\beta|- \pi|\le \frac{\epsilon_1}2$,&$  \widehat s\in\mho_{2<}$, $Type\,  \mathfrak B'', \mathfrak  C''$,\\
    &$\textit{for } \sin 3\beta\gtrless 0,\ \ |\alpha_2-\beta|\le \frac{\epsilon_1}2$,&$ \widehat s\in\mho_{2>},$ $Type\,  \mathfrak B'', \mathfrak  C''$,\\
    &$\textit{for }\sin 3\beta\gtrless 0,\ \ ||\alpha_2-\beta|- \pi|\ge \frac{\epsilon_1}2$,&$ \widehat s\in\mho_{2<} $, $Type\,  \mathfrak B'', \mathfrak  C''$,\\
    &$\textit{for }\sin 3\beta\gtrless 0,\ \ |\alpha_2-\beta|\ge \frac{\epsilon_1}2$,&$\widehat s\in\mho_{2>} $, $Type\,  \mathfrak B'', \mathfrak  C''$,
\end{tabular}
\end{center}
 and
 \begin{center}
\begin{tabular}{lll} 
\multirow{8}{*}{$ \tau_{0,\dagger}=  
 \left\{
{\ba{l}  
  { 0  ,} \\
  {0  ,} \\
  { \pm\epsilon_1   ,} \\
  {\pm\frac{\epsilon_1 }4   ,}\\
  { \mp\epsilon_1   ,} \\
  {\mp\frac{\epsilon_1 }4   ,}\\
  { \mp\epsilon_1   ,} \\
  {\mp\frac{\epsilon_1 }4   ,} 
  \ea}\right.$}&$\textit{for }\sin 3\beta\gtrless 0,\ \ |\alpha_0-\beta|\le \frac{\epsilon_1}2$,&$\widehat s\in  [0,\frac 12]\subset \mho_{0 }$, $Type\ \mathfrak A ,\mathfrak B',\mathfrak C'$,\\
    &$ \textit{for }\sin 3\beta\gtrless 0,\ \ |\alpha_0-\beta|\ge \frac{\epsilon_1}2$,&$\widehat s\in  [0,\frac 12]\subset \mho_{0 }$, $Type\ \mathfrak A ,\mathfrak B',\mathfrak C'$,\\
    &$\textit{for }\sin 3\beta\gtrless 0,\ \ |\alpha_0-\beta|\le \frac{\epsilon_1}2$,&$  \widehat s\in\mho_0,\,Type\, \mathfrak D, \mathfrak E $, \\
    &$ \textit{for }\sin 3\beta\gtrless 0,\ \ |\alpha_0-\beta|\ge \frac{\epsilon_1}2$,&$\widehat s\in\mho_0 ,\,Type\, \mathfrak D, \mathfrak E$,\\
    &$\textit{for }\sin 3\beta\gtrless 0,\ \ |\alpha_0-\beta|\le \frac{\epsilon_1}2$,&$\widehat s\in    \mho_{0,> }$, $Type\ \mathfrak A' ,\mathfrak B',\mathfrak C'$,\\
    &$ \textit{for }\sin 3\beta\gtrless 0,\ \ |\alpha_0-\beta|\ge \frac{\epsilon_1}2$,&$\widehat s\in    \mho_{0,> }$, $Type\ \mathfrak A' ,\mathfrak B',\mathfrak C'$,\\
    &$\textit{for }\sin 3\beta\gtrless 0,\ \ |\alpha_0-\beta|\le \frac{\epsilon_1}2$,&$  \widehat s\in\mho_0 $,  $Type\, \mathfrak B'',\mathfrak C''$,\\
    &$ \textit{for }\sin 3\beta\gtrless 0,\ \ |\alpha_0-\beta|\ge \frac{\epsilon_1}2$,&$\widehat s\in\mho_0 $,  $Type\, \mathfrak B'',\mathfrak C''$,
    \end{tabular}
\end{center}

\begin{center}
\begin{tabular}{lll}
 \multirow{12}{*}{$  \tau_{1,\dagger}=   
 \left\{
{\ba{l}  
{    \mp\pi\pm\epsilon_1 ,} \\
{ \mp\epsilon_1   ,} \\
  {    \mp\pi\pm\epsilon_1 ,} \\
   {   \pm\epsilon_1 , }\\
   {\pm\pi\mp\epsilon_1 ,}\\
   {  \mp\epsilon_1 ,}\\   
   {   \mp\pi\pm\frac{\epsilon_1 }4 ,} \\
   {  \mp\frac{\epsilon_1}4 ,} \\
 {  \mp\pi\pm\frac{\epsilon_1 }4,} \\
  { \pm\frac{\epsilon_1 }4   ,} \\
   {\pm\pi\mp\frac{\epsilon_1}4 ,}\\
   { \mp\frac{\epsilon_1}4 ,}
 \ea}\right.$}
    &$\textit{for }\sin 3\beta\gtrless 0,\ \ ||\alpha_1-\beta|- \pi|\le \frac{\epsilon_1}2$,&$ \widehat s\in\mho_{1<}\ne\phi,\,  Type\, \mathfrak A ,$\\
    &$\textit{for }\sin 3\beta\gtrless0,\ \ |\alpha_1-\beta|\le \frac{\epsilon_1}2$,&$\widehat s\in\mho_{1>}\ne\phi, \,  Type\, \mathfrak A ,$\\
    &$\textit{for }\sin 3\beta\gtrless0,\ \ ||\alpha_1-\beta|- \pi|\le \frac{\epsilon_1}2$,&$\widehat s\in\mho_{1<}, $ $Type\, \mathfrak B'', \mathfrak C''$,\\
    &$\textit{for }\sin 3\beta\gtrless0,\ \ |\alpha_1-\beta|\le \frac{\epsilon_1}2$,&$ \widehat s\in\mho_{1>} ,$ $Type\, \mathfrak B'', \mathfrak C''$,\\
    &$\textit{for }\sin 3\beta\gtrless0,\ \ ||\alpha_1-\beta|- \pi|\le \frac{\epsilon_1}2$,&$ \widehat s\in\mho_{1<} ,$ $Type\, \mathfrak D, \mathfrak E$,\\
     &$\textit{for }\sin 3\beta\gtrless 0,\ \ |\alpha_1-\beta|\le \frac{\epsilon_1}2$,&$ \widehat s\in\mho_{1>} ,$ $Type\, \mathfrak D, \mathfrak E$,\\
    &$ \textit{for }\sin 3\beta\gtrless 0,\ \ ||\alpha_1-\beta|- \pi|\ge \frac{\epsilon_1}2$,&$ \widehat s\in\mho_{1<},\  Type\, \mathfrak A ,$ \\
    &$ \textit{for }\sin 3\beta\gtrless 0,\ \ |\alpha_1-\beta|\ge \frac{\epsilon_1}2$,&$\widehat s\in\mho_{1>},\  Type\, \mathfrak A ,$ \\
    &$\textit{for }\sin 3\beta\gtrless 0,\ \ ||\alpha_1-\beta|- \pi|\ge \frac{\epsilon_1}2$,&$ \widehat s\in\mho_{1<}, $ $Type\, \mathfrak B'', \mathfrak C''$, \\
    &$\textit{for }\sin 3\beta\gtrless 0,\ \ |\alpha_1-\beta|\ge \frac{\epsilon_1}2$,&$\widehat s\in\mho_{1>}, $ $Type\, \mathfrak B'', \mathfrak C''$, \\
    &$\textit{for }\sin 3\beta\gtrless 0,\ \ ||\alpha_1-\beta|- \pi|\ge \frac{\epsilon_1}2$,&$ \widehat s\in\mho_{1<} ,$ $Type\, \mathfrak D, \mathfrak E$,\\
     &$\textit{for }\sin 3\beta\gtrless 0,\ \ |\alpha_1-\beta|\ge \frac{\epsilon_1}2$,&$ \widehat s\in\mho_{1>} ,$ $Type\, \mathfrak D, \mathfrak E$,  
\\     
 \multirow{4}{*}{$ \tau_{2,\dagger}=    
 \left\{
{\ba{l}  
  {\pm\pi\mp\epsilon_1   ,}   \\
   {\mp\epsilon_1  , }  \\
  {\pm\pi\mp\frac{\epsilon_1 }4    ,}   \\
  {\mp\frac{\epsilon_1 }4  ,}  
 \ea}\right.$}&$\textit{for }\sin 3\beta\gtrless 0,\ \ ||\alpha_2-\beta|- \pi|\le \frac{\epsilon_1}2$,&$  \widehat s\in\mho_{2<},$ $Type\, \mathfrak B'', \mathfrak C''$,\\
    &$\textit{for } \sin 3\beta\gtrless 0,\ \ |\alpha_2-\beta|\le \frac{\epsilon_1}2$,&$ \widehat s\in\mho_{2>},$ $Type\, \mathfrak B'', \mathfrak C''$,\\
    &$\textit{for }\sin 3\beta\gtrless 0,\ \ ||\alpha_2-\beta|- \pi|\ge \frac{\epsilon_1}2$,&$ \widehat s\in\mho_{2<},$ $Type\, \mathfrak B'', \mathfrak C''$,\\
    &$\textit{for }\sin 3\beta\gtrless 0,\ \ |\alpha_2-\beta|\ge \frac{\epsilon_1}2$,&$\widehat s\in\mho_{2>}, $ $Type\, \mathfrak B'', \mathfrak C''$,    
    \end{tabular}
\end{center}

\end{definition}

The following lemmas  show   that the deformation fulfills goals of \eqref{E:deform-intro-1} partially. 
\begin{table}
{\begin{center}
\vskip.1in
\begin{tabular}{|l|l|}
   \hline
\bf{Case}& \bf{Type $\mathfrak A$} \cr
   \hline\hline
   $\widehat s \in\mho_0$&$\blacktriangleright \mathfrak{Re}(-i\wp(\frac{\widehat s e^{i\tau }}{{\widehat\sigma}},\beta,X)) $\\
   &$\color{blue}=\frac{X_3}{\widehat\sigma^3} \sin3\tau\sin3\beta\widehat s(\widehat s-  \frac{3\sin2\tau }{ 2\sin3\tau }\widehat s_{1,\ast})^2 +\frac{X_1}{\widehat\sigma}\sin\tau\sin\beta\widehat s(1- \frac{3\sin^22\tau}{4\sin\tau\sin3\tau}\frac{1}{\Delta })$ \\
 \hline
 $\widehat s \in\mho_1$&$\blacktriangleright \mathfrak{Re}(-i\wp(\frac{\widehat s_1 e^{i \tau_1 } +\widehat s_{1,\ast}}{{\widehat\sigma}},\beta,X)) $\\
 &$\color{blue}=  \frac{X_3}{\widehat\sigma^3}\sin3\tau_1\sin3\beta \widehat s_1^3 +\frac{X_1}{\widehat\sigma}\sin\tau_1\sin\beta\widehat s_1(1-\frac{1}{\Delta })  $   \\
   \hline
 \hline
\bf{Case} &\bf{Type $\mathfrak B$,   $ \mathfrak C$} \cr
\hline\hline  
$\widehat s\in\mho_0$& $\blacktriangleright\mathfrak{Re}(-i\wp(\frac{\widehat s e^{i\tau}}{{\widehat\sigma}},\beta,X))$\\  
& $\color{blue}=\frac{X_3}{\widehat\sigma^3}\sin3\tau \sin3\beta \widehat s   (\widehat s-\frac{    - 1+\sqrt{1-\frac 43 \frac{\sin\tau \sin3\tau  }{\sin^22\tau }\Delta}  }{ 3\frac{X_3 }{ \widehat\sigma X_2}\frac{ \sin3\beta}{ \sin2\beta}\frac{2\sin3\tau}{3\sin2\tau}})    (\widehat s- \frac{    - 1-\sqrt{1-\frac 43 \frac{\sin\tau \sin3\tau }{\sin^22\tau }\Delta}  }{ 3\frac{X_3 }{ \widehat\sigma X_2}\frac{ \sin3\beta}{ \sin2\beta}\frac{2\sin3\tau}{3\sin2\tau}})$ \\
\hline
$ \widehat s\in\mho_1$ & $\blacktriangleright\mathfrak{Re}(-i\wp(\frac{\widehat s_1 e^{i\tau_1 } +\widehat s_{1,\ast}}{{\widehat\sigma}},\beta,X))$\\
$\color{black} $&$\color{blue} = \frac{X_3}{\widehat\sigma^3} \sin3\tau_1 \sin3\beta  \widehat s_1^2  (  \widehat s_1+   \frac {\sqrt{1-\Delta}}{3\frac{X_3}{\widehat\sigma X_2}\frac{ \sin3\beta}{ \sin 2\beta}\frac{\sin 3\tau_1}{   3 \sin2\tau_1}}  ) $\\
\hline
$\widehat s\in\mho_2$ & $\blacktriangleright\mathfrak{Re}(-i\wp(\frac{\widehat s_2 e^{i\tau_2 } +\widehat s_{2,\ast}}{{\widehat\sigma}},\beta,X))$\\
$\color{black}  $&$\color{blue} = \frac{X_3}{\widehat\sigma^3} \sin3\tau_2 \sin3\beta  \widehat s_2^2  (  \widehat s_2- \frac {\sqrt{1-\Delta}}{3\frac{X_3}{\widehat\sigma X_2}\frac{ \sin3\beta}{ \sin 2\beta}\frac{\sin 3\tau_2}{   3 \sin2\tau_2}}  ) $\\
\hline\hline
\bf{Case} &\bf{Type $\mathfrak D$, $\mathfrak E$} \cr
\hline\hline  
$\widehat s\in\mho_0$& $\blacktriangleright\mathfrak{Re}(-i\wp(\frac{\widehat s e^{i\tau}}{{\widehat\sigma}},\beta,X))$\\  
& $\color{blue}=\frac{X_3}{\widehat\sigma^3}\sin3\tau \sin3\beta \widehat s   (\widehat s-\frac{    - 1+\sqrt{1-\frac 43 \frac{\sin\tau \sin3\tau  }{\sin^22\tau }\Delta}  }{ 3\frac{X_3 }{ \widehat\sigma X_2}\frac{ \sin3\beta}{ \sin2\beta}\frac{2\sin3\tau}{3\sin2\tau}})    (\widehat s- \frac{    - 1-\sqrt{1-\frac 43 \frac{\sin\tau \sin3\tau }{\sin^22\tau }\Delta}  }{ 3\frac{X_3 }{ \widehat\sigma X_2}\frac{ \sin3\beta}{ \sin2\beta}\frac{2\sin3\tau}{3\sin2\tau}})$ \\
\hline
$\widehat s\in\mho_1$ & $\blacktriangleright\mathfrak{Re}(-i\wp(\frac{\widehat s_1 e^{i\tau_1 } +\widehat s_{1,\ast}}{{\widehat\sigma}},\beta,X))$\\
$\widehat s_{1,\ast}= \widehat s_\pm $ &$\color{blue}= \frac{X_3}{\widehat\sigma^3} \sin3\tau_1 \sin3\beta  \widehat s_1^2  (  \widehat s_1\pm  \frac {\sqrt{1-\Delta}}{3\frac{X_3}{\widehat\sigma X_2}\frac{ \sin3\beta}{ \sin 2\beta}\frac{\sin 3\tau_1}{   3 \sin2\tau_1}}  ) $  
  \cr
  
  \hline
  
\end{tabular}
\end{center} }
\caption{\small $\mathfrak{Re}(-i\wp(\frac{\widehat s}{\widehat\sigma},\beta,X))$ for  deformation defined by Definition \ref{D:deformation}} 
\label{Tb:f3-deformation}
\end{table}
\begin{lemma}\label{L:ideas} 
Define the deformation $\widehat s\mapsto \xi_j=\widehat s_je^{i\tau_{j}}+\widehat s_{j,\ast}$ on $\mho_j$ by Definition \ref{D:deformation}. Therefore, for $j=0,1,2$,  
\begin{align}
{|\widehat s_j e^{i\tau_{j,\dagger}}-\widehat r_je^{i(\alpha_j-\beta)}| }  \ge & \frac 1C \max\{\widehat r_j,\widehat s_j\},\quad\textit{ if $\,\tau_{j,\dagger}\ne 0$},
 \label{E:f3-denominator}
\end{align} and   
\be\mathfrak{Re}(-i\wp(\frac{\xi_j}{\widehat\sigma},\beta,X))\le   0.
 \label{E:definitive}
 \ee
\end{lemma}
\begin{proof}\hfill \\

\begin{itemize}
\item [$\blacktriangleright$] {\bf Proof of \eqref{E:f3-denominator} :}  According to the definition of $\tau_{j,\dagger}$ in Definition \ref{D:deformation}, if $\tau_{j,\dagger}\ne 0$ then
\be \label{E:angle}
|\alpha_j-\beta-\tau_{j,\dagger}|\ge  \frac{\epsilon_1} {4},\quad j=0,1,2.
\ee As a result,  \eqref{E:f3-denominator} is justified.

\item [$\blacktriangleright$] {\bf Proof of \eqref{E:definitive} :} 
In view of Definition \ref{D:deformation}, Table \ref{Tb:f3-deformation}, and Figure \ref{Fg:signature}, 
\begin{itemize}
\item for Type $\mathfrak A $:
\begin{itemize}
\item On $\mho_1$,   both terms of $\mathfrak{Re}(-i\wp(\frac{\xi_j}{\widehat\sigma},\beta,X))$ are of  the same signatures. Along with Definition \ref{D:deformation}  and Figure \ref{Fg:signature}, yields \eqref{E:definitive}.

\item On $\mho_{0}$, according to the assumption, we need only to consider $Type\, \mathfrak A'$ (i.e., $ \Delta\ge 2$),  $\widehat{\sigma}=X_1$. Along with Definition \ref{D:deformation}  and Figure \ref{Fg:signature},      
\begin{align}
& 3\frac{X_3}{X_1^3}\sin3\beta(\widehat s-\frac32\frac{\sin2\tau}{\sin3\tau } \frac{\widehat s_{+}+\widehat s_{-}}{2} )^2 \label{E:a-1-est-2} \\
+&\frac{X_1}{X_1 }\sin \beta(1- \frac{3\sin^22\tau}{4\sin\tau\sin3\tau }  \frac 1\Delta) \nonumber\\
\le  &-| \frac{X_1}{X_1 }\sin \beta(1- \frac{3\sin^22\tau}{4\sin\tau\sin3\tau }  \frac 1\Delta)|\nonumber\\
\le &- C_0|\sin\beta|,\nonumber
\end{align}
 Hence we prove   \eqref{E:definitive}.
\end{itemize}

\item For Type $\mathfrak B'',\mathfrak C''$:
\begin{itemize}
\item on $\mho_0$, from $\widehat s\le \frac 1{2\cos\epsilon_1}  \widehat s_{1,\ast}$ and $\epsilon_1\ll 1$, we prove \eqref{E:definitive} by means of
\be\label{E:bc-0}
\hskip1in (\widehat s-\frac{    - 1+\sqrt{1-\frac 43 \frac{\sin\tau \sin3\tau  }{\sin^22\tau }\Delta}  }{ 3\frac{X_3 }{ \widehat\sigma X_2}\frac{ \sin3\beta}{ \sin2\beta}\frac{2\sin3\tau}{3\sin2\tau}})    (\widehat s- \frac{    - 1-\sqrt{1-\frac 43 \frac{\sin\tau \sin3\tau }{\sin^22\tau }\Delta}  }{ 3\frac{X_3 }{ \widehat\sigma X_2}\frac{ \sin3\beta}{ \sin2\beta}\frac{2\sin3\tau}{3\sin2\tau}}) 
\ge   \frac 18\widehat s^2_{1,\ast}. 
\ee  

\item On $\mho_1$, from Figure \ref{Fg:signature}, $\frac{\sin3\beta}{\sin2\beta}\frac{\sin3\tau_1}{\sin2\tau_1} \gtrless 0$ on $\mho_{1,\lessgtr}$. It reduces to proving \eqref{E:definitive} on $\mho_{1,>}$. From the definition of $\mho_{1,>}$ and \eqref{E:roots-distance}, we have 
\be\label{E:bc-1-new}
\widehat s_1\le  \frac {\widehat s_{2,\ast}-\widehat s_{1,\ast}  }{2  \cos\epsilon_1}= \frac{1}{\cos\epsilon_1} |\frac{\sqrt{1-\Delta}}{3{ \frac{X_3}{\widehat\sigma X_2 } \frac{ \sin3\beta  }{  \sin2\beta} }}|.
\ee Therefore, 
\be\label{E:bc-1}
\hskip1.5in \widehat s_1+   \frac {\sqrt{1-\Delta}}{3\frac{X_3}{\widehat\sigma X_2}\frac{ \sin3\beta}{ \sin 2\beta}\frac{\sin 3\tau_1}{   3 \sin2\tau_1}} \le  \frac 12\frac {\sqrt{1-\Delta}}{3\frac{X_3}{\widehat\sigma X_2}\frac{ \sin3\beta}{ \sin 2\beta}\frac{\sin 3\tau_1}{   3 \sin2\tau_1}} \le 0, \quad \widehat s_1\in\mho_{1,>}.
\ee As a result,  \eqref{E:definitive} follows. 

\item On $\mho_2$, from Figure \ref{Fg:signature}, $-\frac{\sin3\beta}{\sin2\beta}\frac{\sin3\tau_2}{\sin2\tau_2} \lessgtr 0$ on $\mho_{2,\lessgtr}$. It reduces to proving \eqref{E:definitive} on $\mho_{2,<}$. From the definition of $\mho_{2,<}$, we have 
\be\label{E:bc-2-new}
\widehat s_2\le  \frac {\widehat s_{2,\ast}-\widehat s_{1,\ast}  }{2  \cos\epsilon_1} = \frac{1}{\cos\epsilon_1}|\frac{\sqrt{1-\Delta}}{3{ \frac{X_3}{\widehat\sigma X_2 } \frac{ \sin3\beta  }{  \sin2\beta} }}|.
\ee Therefore, 
\be\label{E:bc-2}
\hskip1.35in \widehat s_2-   \frac {\sqrt{1-\Delta}}{3\frac{X_3}{\widehat\sigma X_2}\frac{ \sin3\beta}{ \sin 2\beta}\frac{\sin 3\tau_2}{   3 \sin2\tau_2}}\le  -\frac 12\frac {\sqrt{1-\Delta}}{3\frac{X_3}{\widehat\sigma X_2}\frac{ \sin3\beta}{ \sin 2\beta}\frac{\sin 3\tau_2}{   3 \sin2\tau_2}} \le 0  ,\quad \widehat s_2\in \mho_{2,<}. 
\ee Hence  \eqref{E:definitive} is proved. 
\end{itemize} 

\item For Type $ \mathfrak B', \mathfrak C'$: it is sufficient to consider \eqref{E:definitive} for   $\widehat{s}\in \mho_{0,>}$. Hence \eqref{E:definitive} is proved by $\widehat s_\pm\le 0$ and $\widehat s>\frac 12$. 

\item For Type $\mathfrak D,\mathfrak E$: 

\begin{itemize}
\item on $\mho_0$,  we prove \eqref{E:definitive} by means of $\widehat s\le \frac 1{2\cos\epsilon_1}  \widehat s_{1,\ast}$, $\epsilon_1\ll 1$,   and
\be\label{E:de-0}
\begin{split}
\hskip1. in &(\widehat s-\frac{    - 1+\sqrt{1-\frac 43 \frac{\sin\tau \sin3\tau  }{\sin^22\tau }\Delta}  }{ 3\frac{X_3 }{ \widehat\sigma X_2}\frac{ \sin3\beta}{ \sin2\beta}\frac{2\sin3\tau}{3\sin2\tau}})    (\widehat s- \frac{    - 1-\sqrt{1-\frac 43 \frac{\sin\tau \sin3\tau }{\sin^22\tau }\Delta}  }{ 3\frac{X_3 }{ \widehat\sigma X_2}\frac{ \sin3\beta}{ \sin2\beta}\frac{2\sin3\tau}{3\sin2\tau}}) 
\\
\hskip1. in \le &  -\frac 14\widehat s\widehat s_{1,\ast}\le  -\frac 12\widehat s^2.
\end{split}
\ee 

\item on $\mho_1$, if $\widehat s_{1,\ast}=\widehat s_+$, by means of  Table \ref{Tb:dynamic-FNH} and Figure \ref{Fg:signature}, $\frac{\sin3\beta}{\sin2\beta}\frac{\sin3\tau_1}{\sin2\tau_1} \lessgtr 0$ on $\mho_{1,\lessgtr}$. It reduces to proving \eqref{E:definitive} on $\mho_{1,<}$. From the definition of $\mho_{1,<}$, \eqref{E:roots-distance}, and two roots are opposite signs, we have 
\be\label{E:de-1}
\begin{split}
& \widehat s_1+ \frac {\sqrt{1-\Delta}}{3\frac{X_3}{\widehat\sigma X_2}\frac{ \sin3\beta}{ \sin 2\beta}\frac{\sin 3\tau_1}{   3 \sin2\tau_1}} \le +\frac 12  \frac {\sqrt{1-\Delta}}{3\frac{X_3}{\widehat\sigma X_2}\frac{ \sin3\beta}{ \sin 2\beta}\frac{\sin 3\tau_1}{   3 \sin2\tau_1}}  \le -\frac 14\widehat s _1 ,
\end{split}
\ee and  \eqref{E:definitive} is proved.

On $\mho_1$, if $\widehat s_{1,\ast}=\widehat s_-$, by means of  Table \ref{Tb:dynamic-FNH} and Figure \ref{Fg:signature}, $-\frac{\sin3\beta}{\sin2\beta}\frac{\sin3\tau_1}{\sin2\tau_1} \lessgtr 0$ on $\mho_{1,\lessgtr}$. It reduces to proving \eqref{E:definitive} on $\mho_{1,<}$. From the definition of $\mho_{1,<}$, \eqref{E:roots-distance}, and two roots are opposite signs, we have 
\be\label{E:de-1-s-}
\begin{split}
& \widehat s_1- \frac {\sqrt{1-\Delta}}{3\frac{X_3}{\widehat\sigma X_2}\frac{ \sin3\beta}{ \sin 2\beta}\frac{\sin 3\tau_1}{   3 \sin2\tau_1}} \le -\frac 12  \frac {\sqrt{1-\Delta}}{3\frac{X_3}{\widehat\sigma X_2}\frac{ \sin3\beta}{ \sin 2\beta}\frac{\sin 3\tau_1}{   3 \sin2\tau_1}}  \le -\frac 14\widehat s _1 ,
\end{split}
\ee and  \eqref{E:definitive} is proved. 
\end{itemize}

\end{itemize}

\end{itemize} 
\end{proof}

\begin{lemma}\label{L:p-estimate}  Let $\tilde\sigma=   \max\{ 1,  X_1 ,  \sqrt{ X_2 },  \sqrt[3]{ X_3 }\}$ and $X_j\ge 0$ defined by \eqref{E:phase}.
\begin{itemize}
\item [(i)] For $\,Type\,\mathfrak B'',\,\mathfrak C'',\,\mathfrak D,\,\mathfrak E$ on $\mho_j$, $j=0,1,2$, and $\,Type\,\mathfrak B',\,\mathfrak C'$ on $\mho_{0,>}$, one has
\begin{align}
\mathfrak{Re}(-i\wp(\frac{ \widehat s_je^{i\tau_{j,\dagger}}+ \widehat  s_{j,\ast}}{\widehat  \sigma},\beta,X)) \le &
     -\frac {| \sin 3\beta|}C \widehat  s_j^{3},\hskip.72in \widehat\sigma=\sqrt[3]{X_3};\label{E:f3-deform-wp}\\
     \mathfrak{Re}(-i\wp(\frac{\widehat   s_je^{i\tau_{j,\dagger}}+ \widehat  s_{j,\ast}}{\widehat  \sigma},\beta,X)) \le& 
     -\frac {| \sin 2\beta|}C \min\{ \widehat s_j^{3},  \widehat s_j^{2}\}, \ \widehat\sigma=\sqrt{X_2}.\label{E:f2-deform-wp}
\end{align}

\item [(ii)] For $\,Type\,\mathfrak A$ on $\mho_1$, 
\begin{align}
\mathfrak{Re}(-i\wp(\frac{ \widehat s_je^{i\tau_{j,\dagger}}+ \widehat  s_{j,\ast}}{\widehat  \sigma},\beta,X)) \le &
     -\frac {| \sin 3\beta|}C \widehat  s_j^{3},\hskip.72in \widehat\sigma=\sqrt[3]{X_3};\label{E:f3-deform-wp-1}\\
     \mathfrak{Re}(-i\wp(\frac{\widehat   s_je^{i\tau_{j,\dagger}}+ \widehat  s_{j,\ast}}{\widehat  \sigma},\beta,X)) \le& 
     -\frac {|  \sin 2\beta|}C   \widehat s_j^{3}, \hskip.72in   \widehat\sigma=\sqrt{X_2}=\tilde\sigma.\label{E:f2-deform-wp-1}
\end{align}
\item [(iii)] For 
\begin{itemize}
\item  $Type\,\mathfrak B'',\,\mathfrak C'',\,\mathfrak D,\,\mathfrak E$ on $\mho_{0}$;
\item $Type\,\mathfrak A',\,\mathfrak B',\,\mathfrak C' $  on $\mho_{0,>}$,
\end{itemize}  one has
\begin{align}  
\mathfrak{Re}(-i\wp(\frac{\widehat   s_0e^{i\tau_{0,\dagger}}}{\widehat  \sigma},\beta,X))\le &- \frac 1C| \sin\beta|  \widehat   s, \hskip1 in \widehat\sigma=  {X_1}=\tilde\sigma.     \label{E:all-0}  
\end{align}
\end{itemize}

\end{lemma}
\begin{proof}  \hfill \\
\begin{itemize}
\item [$\blacktriangleright$] {\bf Proof of \eqref{E:f3-deform-wp} and \eqref{E:f3-deform-wp-1}:}
In view of $\widehat\sigma=\sqrt[3]{X_3}$, Definition \ref{D:deformation}, Table \ref{Tb:f3-deformation}, and Figure \ref{Fg:signature}, 
\begin{itemize}
\item for Type $\mathfrak A $ on $\mho_1$, from Figure \ref{Fg:signature},   both terms of $\mathfrak{Re}(-i\wp(\frac{\xi_1}{\sqrt[3]{X_3}},\beta,X))$ are of  the same signatures. Therefore,
\begin{align*}
&{Re}(-i\wp(\frac{\xi_1}{ \sqrt[3]{X_3}},\beta,X))\\
=& \sin3\tau_1\sin3\beta \widehat s_1^3 +\frac{X_1}{\sqrt[3]{X_3}}\sin\tau_1\sin\beta\widehat s_1(1-\frac{1}{\Delta })\\
\le &-|\sin3\tau_1\sin3\beta| \widehat s_1^3
\end{align*} and  \eqref{E:f3-deform-wp-1} is proved.
\item For Type $\mathfrak B'',\mathfrak C''$: 
\begin{itemize}
\item On $\mho_0$, by means of \eqref{E:bc-0},   we prove \eqref{E:f3-deform-wp}. 

\item On $\mho_1$, from Figure \ref{Fg:signature}, $\frac{\sin3\beta}{\sin2\beta}\frac{\sin3\tau_1}{\sin2\tau_1} \gtrless 0$ on $\mho_{1,\lessgtr}$. It reduces to proving \eqref{E:f3-deform-wp} on $\mho_{1,>}$. Then \eqref{E:f3-deform-wp} follows from \eqref{E:bc-1-new}, \eqref{E:bc-1} and the definition  of interval  $\mho_{1,>}$.
\item On $\mho_2$, from Figure \ref{Fg:signature}, $-\frac{\sin3\beta}{\sin2\beta}\frac{\sin3\tau_2}{\sin2\tau_2} \lessgtr 0$ on $\mho_{2,\lessgtr}$. It reduces to proving \eqref{E:f3-deform-wp} on $\mho_{2,<}$. Then \eqref{E:f3-deform-wp} follows from \eqref{E:bc-2-new}, \eqref{E:bc-2} and the definition  of interval  $\mho_{2,<}$.
\end{itemize}

\item For Type $ \mathfrak B', \mathfrak C'$ on $\mho_{0,>}$, \eqref{E:f3-deform-wp} is proved since $\widehat s_\pm\le 0$ and $\widehat s>\frac 12$.

\item For Type $\mathfrak D,\mathfrak E$: 
\begin{itemize}
\item on $\mho_0$, by means of \eqref{E:de-0},   we prove \eqref{E:f3-deform-wp}. 
\item On $\mho_1$, from Figure \ref{Fg:signature}, $\widehat s_{1,\ast}=\widehat s_\pm$, $\pm\frac{\sin3\beta}{\sin2\beta}\frac{\sin3\tau_1}{\sin2\tau_1} \lessgtr 0$ on $\mho_{1,\lessgtr}$. It reduces to proving \eqref{E:f3-deform-wp} on $\mho_{1,<}$. Then \eqref{E:f3-deform-wp} follows from \eqref{E:de-1}, \eqref{E:de-1-s-}, and the definition  of interval  $\mho_{1,<}$.
\end{itemize}

\end{itemize}
\item [$\blacktriangleright$] {\bf Proof of \eqref{E:f2-deform-wp}:} In view of $\widehat\sigma=\sqrt{X_2}$, Definition \ref{D:deformation}, Table \ref{Tb:f3-deformation}, and Figure \ref{Fg:signature}, 
\begin{itemize}
\item For Type $\mathfrak B''$ on $\mho_j$, $j=0,1,2$:  Notice  \eqref{E:determinant} implies 
\be\label{E:f2-coef}
3\frac{X_3}{X_2^{3/2}}\sin3\beta=\frac{X_2^{1/2}}{X_1 } \frac{\sin2\beta}{\sin\beta}\Delta \sin2\beta .
\ee Together with Table \ref{Tb:dynamic-FNH},  $\Delta>\frac 12$,  and similar argument as that for  $\widehat{\sigma}=\sqrt[3]{X_3}$, yields 
\be\label{E:f2-12-a-1}
\begin{gathered}
 \mathfrak{Re}(-i\wp(\frac{\widehat s_je^{i\tau_{j,\dagger}}+\widehat s_{j,\ast}}{ \widehat \sigma },\beta,X)) \le -\frac 1C| \sin 2\beta|\widehat s_j^{3}.
  \end{gathered}
 \ee 
 
\item For Type $   \mathfrak B'$ on $\mho_{0,>}$: Using \eqref{E:f2-coef}, Table \ref{Tb:dynamic-FNH},  $\Delta>\frac 12$,  $\widehat{s}_\pm\le 0$,  and   $\widehat{s}>\frac 12$, yields 
\be\label{E:f2-12-a-2}
\begin{gathered}
 \mathfrak{Re}(-i\wp(\frac{\widehat s_je^{i\tau_{j,\dagger}}+\widehat s_{j,\ast}}{ \widehat \sigma },\beta,X)) \le -\frac 1C| \sin 2\beta|\widehat s_j^{3}.
  \end{gathered}
 \ee

\item For Type $ \mathfrak C'', \mathfrak D, \mathfrak E$ on $\mho_j$, $j=0,1,2$:   
Via \eqref{E:f2-coef}, similar argument as that for   $\widehat{\sigma}=\sqrt[3]{X_3}$, 
\begin{align}
| \widehat s_j+(-1)^{j+1}\frac {\sqrt{1-\Delta}}{3\frac{X_3}{\widehat\sigma X_2}\frac{ \sin3\beta}{ \sin 2\beta}\frac {\sin 3\tau_j}{   2 \sin2\tau_j} }    |\ge\frac 1C|  \frac{ 1}{ \frac{X _3}{X_2^{3/2}}\frac{\sin3\beta}{\sin 2 \beta}  }|,& \  \textit{on $\mho_j$, $j=1,2$, for $Type\ \mathfrak C''$},\nonumber\\
| \widehat s_1\pm\frac {\sqrt{1-\Delta}}{3\frac{X_3}{\widehat\sigma X_2}\frac{ \sin3\beta}{ \sin 2\beta}\frac {\sin 3\tau_1}{   2 \sin2\tau_1} }    |\ge\frac 1C|  \frac{1}{ \frac{X _3}{X_2^{3/2}}\frac{\sin3\beta}{\sin 2 \beta}  }|,&  \ \textit{on $\mho_1, 0<\widehat s_\pm$  for $Type\ \mathfrak D,\mathfrak E $,}\nonumber \\
| \widehat s_{1,\ast}-\widehat s_{2,\ast}|\ge\frac 1C|  \frac{ 1 }{ \frac{X _3}{X_2^{3/2}}\frac{\sin3\beta}{\sin 2 \beta}  }| ,&  \ \textit{for $Type\ \mathfrak C,\mathfrak D, \mathfrak E $}.\label{E:2-0-1}
\end{align}
As a result,
\be\label{E:f2-12-cde}
\begin{gathered}
   \mathfrak{Re}(-i\wp(\frac{\widehat s_je^{i\tau_{j,\dagger}}+\widehat s_{j,\ast}}{ \widehat \sigma },\beta,X))\le-  \frac 1C|\sin  2\beta | \widehat s_j^2.
\end{gathered} 
\ee
\item For Type $ \mathfrak C' $ on $\mho_{0,>}$: 
Thanks to $\widehat s_\pm\le 0$, $\widehat{s}>\frac 12$, $0\le \Delta\le \frac 12$, \eqref{E:f2-coef}, and \eqref{E:2-0-1},  
\be\label{E:f2-12-cde-new}
\begin{gathered}
   \mathfrak{Re}(-i\wp(\frac{\widehat s_je^{i\tau_{j,\dagger}}+\widehat s_{j,\ast}}{ \widehat \sigma },\beta,X))\le-  \frac 1C|\sin  2\beta | \widehat s_j^2.
\end{gathered} 
\ee
\end{itemize} Hence \eqref{E:f2-deform-wp} is justified.

\item [$\blacktriangleright$] {\bf Proof of   \eqref{E:f2-deform-wp-1}:} In this case, 
$\widehat\sigma=\sqrt{X_2}=\max\{X_1,\sqrt{X_2},\sqrt[3]{X_3}\}$ for Type $ \mathfrak A $ on $\mho_1$.  Together with  \eqref{E:f2-coef}, Table \ref{Tb:dynamic-FNH},  $\Delta>\frac 12$,   yields 
\be\label{E:f2-12-a}
\begin{gathered}
 \mathfrak{Re}(-i\wp(\frac{\widehat s_je^{i\tau_{j,\dagger}}+\widehat s_{j,\ast}}{ \widehat \sigma },\beta,X)) \le -\frac 1C| \sin 2\beta|\widehat s_j^{3}.
  \end{gathered}
 \ee  Hence \eqref{E:f2-deform-wp-1} is justified.

 \item[$\blacktriangleright$] {\bf Proof of \eqref{E:all-0} :} In this case, 
$\widehat\sigma= X_1=\max\{X_1,\sqrt{X_2},\sqrt[3]{X_3}\}$,
\begin{itemize}
\item For Type $ \mathfrak A '$ on $\mho_0$, the estimates then follows from \eqref{E:all-0}.
\item For $\,Type\,\mathfrak B'',\,\mathfrak C'',\,\mathfrak D,\,\mathfrak E$  on $\mho_0$, and  for $\,Type\,\mathfrak B',\,\mathfrak C',$   on $\mho_{0,>}$,  we use 
\be\label{E:roots}  
3\frac{X_3}{X_1^3}\sin3\beta =\left(\frac{X_2}{X_1^2}\frac{\sin2\beta}{\sin\beta}\right)^2\sin\beta\Delta,\quad 3\frac{X_3}{X_1X_2 }\frac{\sin3\beta}{\sin2\beta}=\frac{X_2}{X_1^2}\frac{\sin2\beta}{\sin\beta}\Delta.  
\ee Hence  \eqref{E:all-0} follows from
\be\label{E:u-0-f1}
\qquad \left\{
{\ba{ll}   \eqref{E:f3-deform-s}, \eqref{E:roots},  \frac 12\le\Delta\le 1, \textit{ \eqref{E:determinant}, Table \ref{Tb:f3-deformation}}&    \textit{ $Type\,\mathfrak B    $},\\
\textit{\eqref{E:determinant}, Table \ref{Tb:f3-deformation}},  \eqref{E:roots},  |(\widehat s -\widehat s_+')(\widehat s -\widehat s_-')|\ge \frac 1C  \frac{ 1 }{ (\frac{X _2}{X_1^2}\frac{\sin2\beta}{\sin   \beta})^2|\Delta|  },&   \textit{ $Type\,\mathfrak C,\mathfrak D$},\\\textit{\eqref{E:determinant}, Table \ref{Tb:f3-deformation}},  \eqref{E:roots},  |(\widehat s -\widehat s_+')(\widehat s -\widehat s_-')|\ge \frac 1C  \frac{ 1 }{ (\frac{X _2}{X_1^2}\frac{\sin2\beta}{\sin   \beta})^2|\Delta|  },&    \textit{  $Type\,\mathfrak E$}.\ea}\right.   \ee    
\end{itemize}  
\end{itemize}

\end{proof}

\begin{definition}\label{D:dominant} Let $\tilde\sigma=   \max\{ 1,  X_1 ,  \sqrt{ X_2 },  \sqrt[3]{ X_3 }\}$   with $X_j$ defined by \eqref{E:phase}, and scaled coordinates
\[
\zeta =\kappa_1+\frac{\tilde s }{\tilde \sigma}e^{i\beta}\quad\in D_{\kappa_1} 
\]  by replacing $\widehat \sigma$, $\widehat s$, $\widehat s_j$, $\widehat s_{j,\ast}$, $\widehat \lambda$, $\widehat \lambda_{j,\ast}$ by $\tilde \sigma$, $\tilde s$, $\tilde s_j$, $\tilde s_{j,\ast}$, $\tilde \lambda$, $\tilde \lambda_{j,\ast}$ in Definition \ref{D:deformation}. We decompose $X_1,X_2\ge 0,X_3> 0$ into following three cases
\be\label{E:123-grad-assump}
\begin{array}{rl}
(F1)&   \ \tilde\sigma=   \ \ \    X_1  ,     \\
(F2)&   \  \tilde\sigma=   \sqrt{ X_2 } ,    \\
(F3)&   \  \tilde\sigma=  \sqrt[3]{ X_3 }.      \end{array}
 \ee

\end{definition}
\begin{proposition}\label{P:non-homogeneous-F3}
For Case $(F3)$, $(F2)$, and $f\in L^\infty (D_{\kappa_j})$ is $\tilde s$-holomorphic, 
\begin{align}
   |I_4| _{C^\mu_{\tilde\sigma}(D_{\kappa_j,\frac 1{\tilde\sigma}})}\le &  C\epsilon_0   |  f|_{L^\infty(D_{\kappa_j})},\label{E:32-CIE}\\
|I_5| _{L^\infty(D_{\kappa_j )}}\le &  C\epsilon_0   |  f|_{L^\infty(D_{\kappa_j})}.\label{E:32-CIE-1}   
 \end{align}    
 
 Moreover, let
 \be\label{E:I-5-2}
I_{\frac{5}{2}}= -\frac {\theta(\tilde r-\frac12)}{2\pi i}\iint_{  \tilde s< {\tilde\sigma} \delta} \frac{  \widetilde\gamma_j(\tilde s, \beta) e^{-i\wp(\frac {\tilde s}{{\tilde\sigma}},\beta,X)}f(\frac {\tilde s}{ {\tilde\sigma} },-\beta,X) }{{\tilde\zeta}-\tilde\lambda}d\overline{\tilde\zeta} \wedge d{\tilde\zeta}, 
\ee one has
\be\label{E:I-5-2-est}
\begin{split}
&\textit{$I_{\frac52}(\frac{\tilde s }{\tilde\sigma} ,\beta,X)$ is holomorphic in $\tilde s>\frac12$}.
\end{split}
\ee 

\end{proposition}
  
\begin{proof}  \hfill \\

\begin{itemize}
\item {\bf Estimates for $|I_4|_{C^\mu_{\tilde\sigma}(D_{\kappa_1,\frac 1{\tilde\sigma}})}$ :} 
 Using the $\tilde s$-holomorphic property of $f$, and a residue theorem,
\begin{align}
    I_{4} =& - \frac { \theta(1-\tilde r)}{2\pi i}\int_{-\pi}^\pi  d\beta[\partial_ {\beta}   \ln (1-\gamma |\beta| )]   
   \{(  \int_{   S _<}+\int_{  \Gamma  _{40} }) \frac{  e^{-i \wp(\frac{\xi_0 }{{\tilde\sigma}},\beta,X)}f(\frac{ \xi_0 }{\tilde\sigma} ,-\beta,X)}{ \tilde s_0e^{i\tau_0}  -  \tilde r_0 e^{ i (\alpha_0 -\beta) } } d{\xi_0 } \label{E:21-I-4-4-dec-1}\\
   +& \int_{\Gamma  _{41} }    \frac{  e^{-i \wp(\frac{\xi_1 }{{\tilde\sigma}},\beta,X)}f(\frac{ \xi_1 }{\tilde\sigma} ,-\beta,X)}{ \tilde s_1e^{i\tau_1}  -  \tilde r_1 e^{ i (\alpha_1 -\beta) } } d{\xi_1 }+\int_{\Gamma  _{42} }    \frac{  e^{-i \wp(\frac{\xi_2 }{{\tilde\sigma}},\beta,X)}f(\frac{ \xi_2 }{\tilde\sigma} ,-\beta,X)}{ \tilde s_1e^{i\tau_2}  -  \tilde r_1 e^{ i (\alpha_2 -\beta) } } d{\xi_2 }\nonumber\\
+&  \int_{   S  _>  }\frac{  e^{-i \wp(\frac{\xi_h }{{\tilde\sigma}},\beta,X)}f(\frac{ \xi_h }{\tilde\sigma} ,-\beta,X)}{ \tilde s_he^{i\tau_h}  -  \tilde r_h e^{ i (\alpha_h -\beta) } } d{\xi_h }  \}, \nonumber \end{align}with 
\begin{align}
   S    _{<}(\tilde\sigma,\beta,X,\lambda)=& \{\xi_0 : \tilde s=2      \},\label{21-I4-01-contour}  \\
 \Gamma   _{40} (\tilde\sigma,\beta,X,\lambda)= &
  \{\xi_0   : 
  \tilde s\in (2,\tilde\sigma\delta)\cap  \mho _0 ,\  \textit{$\tau_0 =\tau_{0, \dagger}$}    \}, 
  \nonumber\\
  \Gamma   _{41} (\tilde\sigma,\beta,X,\lambda)= &
 \{\xi_1   : 
  \tilde s\in (2,\tilde\sigma\delta)\cap  \mho _1 ,\  \textit{$\tau_1 =\tau_{1, \dagger}$}    \},  
  \nonumber\\
  \Gamma   _{42} (\tilde\sigma,\beta,X,\lambda)= &
 \{\xi_2   : 
  \tilde s\in (2,\tilde\sigma\delta)\cap  \mho _2 ,\  \textit{$\tau_2 =\tau_{2, \dagger}$}    \},  
  \nonumber\\
  S _>(\tilde\sigma,\beta,X,\lambda)  = & \{\xi_h   : h=\sup_{\mho_j\ne\phi}j,\, \tilde s=\tilde \sigma\delta   \},
    \nonumber 
\end{align}
and   $\xi_j  $, $\tau_j$,  $\tau_{j, \dagger}$, $\mho_j=\mho_j(\tilde\sigma,\beta,X)$  defined by Definition \ref{D:deformation}. 

In view of \eqref{E:assumption}, $\tilde r<1$, \eqref{E:f3-denominator}, \eqref{E:definitive}, and \eqref{21-I4-01-contour},%
\begin{multline}\label{E:homo-4}
| \frac { \theta(1-\tilde r)}{2\pi i}  \int_{-\pi}^\pi  d\beta[\partial_ {\beta}   \ln (1-\gamma |\beta| )]  \\
\times(\int_{S _<}+\int_{S  _>  })\frac{  e^{-i \wp(\frac{\tilde s}{\tilde\sigma}e^{i\tau},\beta,X) }f(\frac{\tilde s}{\tilde\sigma}e^{i\tau},-\beta,X)}{\tilde se^{i\tau}-  \tilde r e^{ i (\alpha-\beta) }} d\tilde s e^{i\tau}|_{C^\mu_{\tilde\sigma}(D_{\kappa_1,\frac 1{\tilde\sigma}})}
\le C\epsilon_0|f|_{L^\infty(D_{\kappa_1})}.
\end{multline}

Applying \eqref{E:f3-deform-wp} and \eqref{E:f3-deform-wp-1}  for Case $(F3)$, or \eqref{E:f2-deform-wp} and \eqref{E:f2-deform-wp-1} for $(F2)$,  \eqref{E:f3-denominator}, $\tilde r<1$, \eqref{E:definitive},   \eqref{21-I4-01-contour}, and improper integrals,
\begin{align}
&\sum_{j=0}^2| \frac { \theta(1-\tilde r)}{2\pi i}  \int_{-\pi}^\pi  d\beta[\partial_ {\beta}   \ln (1-\gamma |\beta| )] \label{E:homo-4-gamma}\\ 
 \times& \int_{\Gamma _{4j}}\frac{  e^{-i \wp(\frac{\tilde s_j}{\tilde\sigma}e^{i\tau_j},\beta,X) }f(\frac{\tilde s_j}{\tilde\sigma}e^{i\tau_j},-\beta,X)}{\tilde s_je^{i\tau_j}-  \tilde r_j e^{ i (\alpha_j-\beta) }} d\tilde s e^{i\tau_j}|_{C^\mu_{\tilde\sigma}(D_{\kappa_1,\frac 1{\tilde\sigma}})}\nonumber\\
\le &C \sum_{n= 1 }^2\sum_{j=0}^2|\frac { \theta(1-\tilde r)}{2\pi i}  \int_{-\pi}^\pi  d\beta[\partial_ {\beta}   \ln (1-\gamma |\beta| )]e^{-i(n-1)\beta}\nonumber\\\times&
 \int_{\Gamma _{4j}} \frac{ e^{-i \wp(\frac{\tilde s_j}{\tilde\sigma}e^{i\tau_j},\beta,X)}f(\frac{\tilde s_j}{\tilde\sigma}e^{i\tau_j},-\beta,X)} {( \tilde s_je^{i\tau_j}-  \tilde r_j e^{ i (\alpha_j-\beta) })^n  }d\tilde s_j e^{i\tau_j}|_{L^\infty (D_{\kappa_1})} \nonumber\\
\le &
\left\{
{\ba{ll} 
C \epsilon_0|f|_{L^\infty(D_{\kappa_1})} \int_{-\pi}^\pi  d\beta 
 \int  e^{-\frac 1C\tilde s^3_j|\sin 3\tau_\dagger\sin 3\beta|}d\tilde s _j &\textit{if }\tilde\sigma=\sqrt[3]{X_3},\\
 C\epsilon_0|f|_{L^\infty(D_{\kappa_1})}  \int_{-\pi}^\pi   d\beta  
 \int  e^{-\frac 1C\tilde s^2_j|\sin 2\tau_\dagger\sin 2\beta|}d\tilde s _j &\textit{if }\tilde\sigma=\sqrt {X_2},\ea}\right.   
  \nonumber\\     
 \le& 
 \left\{
{\ba{ll}   C\epsilon_0|f|_{L^\infty(D_{\kappa_1})}  \int_{-\pi}^\pi   d\beta  \frac{1}{\sqrt[3]{|\sin 3\beta|}} 
 \int_0^\infty  e^{-  t^3 |\sin  3\tau_\dagger |   } dt &\textit{if }\tilde\sigma=\sqrt[3]{X_3},\\
 C\epsilon_0|f|_{L^\infty(D_{\kappa_1})}  \int_{-\pi}^\pi   d\beta  \frac{1}{\sqrt {|\sin 2\beta|}} 
 \int_0^\infty  e^{-  t^2 |\sin  2\tau_\dagger |   } dt &\textit{if }\tilde\sigma=\sqrt {X_2},\ea}\right.   
  \nonumber\\ 
  \le &C\epsilon_0|f|_{L^\infty(D_{\kappa_1})}.\nonumber
\end{align} 

Combining \eqref{E:21-I-4-4-dec-1}, \eqref{E:homo-4},   and \eqref{E:homo-4-gamma}, we derive
\be\label{E:i-I_4-k}
\begin{split}
&|I_4|_{C^\mu_{\tilde\sigma}(D_{\kappa_1,\frac 1{\tilde\sigma}})}\le C\epsilon_0 |f|_{L^\infty(D_{\kappa_1})}.
\end{split}
\ee

\item {\bf Estimates for $|I_5|_{L^\infty(D_{\kappa_1})}$
:} Using the $\tilde s$-holomorphic property of $f$, and the residue theorem,
\begin{align}
  I_{5} 
= &- \frac { \theta( \tilde r-1)}{2\pi i}  \int_{-\pi}^\pi  d\beta[\partial_ {\beta}   \ln (1-\gamma |\beta| )]  \{\int_{  \Gamma  _{50} } \frac{ e^{-i \wp(\frac{\xi_0 }{{\tilde\sigma}},\beta,X)}f(\frac{ \xi_0 }{\tilde\sigma} ,-\beta,X)}{\tilde s_0 e^{i\tau_0} -  \tilde r_0 e^{ i (\alpha_0 -\beta) }} d{\xi_0 }  \label{E:21-I-5-dec-1}\\
+& \int_{\Gamma  _{51} }     \frac{  e^{-i \wp(\frac{\xi_1 }{{\tilde\sigma}},\beta,X)}f(\frac{ \xi_1 }{\tilde\sigma} ,-\beta,X)}{\tilde s_1 e^{i\tau_1}-  \tilde r_1 e^{ i (\alpha_1 -\beta) }} d{\xi_1 } + \int_{\Gamma  _{52} }     \frac{  e^{-i \wp(\frac{\xi_2 }{{\tilde\sigma}},\beta,X)}f(\frac{ \xi_2 }{\tilde\sigma} ,-\beta,X)}{\tilde s_2 e^{i\tau_2}-  \tilde r_2e^{ i (\alpha_2 -\beta) }} d{\xi_2 }\nonumber\\
+& \int_{   S  _>  }  \frac{  e^{-i \wp(\frac{\xi_h }{{\tilde\sigma}},\beta,X)}f(\frac{ \xi_h }{\tilde\sigma} ,-\beta,X)}{\tilde s_h e^{i\tau_h}  -  \tilde r_h e^{ i (\alpha_h -\beta) }} d{\xi_h }\}\nonumber\\
-&    \theta(\tilde r-1)\theta(\tilde r_1-\frac 14)\theta(\tilde r_2-\frac 14) \int_{\beta\in\mathfrak\Delta(\lambda)}   d\beta[\partial_ {\beta}   \ln (1-\gamma |\beta| )] \sgn(\beta) \nonumber\\
\times& e^{-i \wp(\frac{\tilde r e^{ i (\alpha-\beta) }}{\tilde \sigma} , \beta,X)  }f(\frac {\tilde r e^{i(\alpha-\beta)}}{ {\tilde\sigma} },-\beta ,X)  , \nonumber 
\end{align}where 
{\be\label{E:alpha-beta}
 \mathfrak\Delta (\lambda) \equiv \{\beta: \ 
\left\{
{\ba{lll}  
  |\alpha_0-\beta|<\frac{\epsilon_1}2, &(\alpha_0-\beta)\beta<0,&\tilde r\in\mho_0,  \\
  ||\alpha_1-\beta|-\pi|<\frac{\epsilon_1}2,&(\alpha_1-\beta)\beta<0, &\tilde r\in\mho_{1<}, \\
  |\alpha_1-\beta|<\frac{\epsilon_1}2,&(\alpha_1-\beta)\beta<0, &\tilde r\in\mho_{1>},  \\
  ||\alpha_2-\beta|-\pi|<\frac{\epsilon_1}2,&(\alpha_2-\beta)\beta<0, &\tilde r\in\mho_{2<}, \\
  |\alpha_2-\beta|<\frac{\epsilon_1}2,&(\alpha_2-\beta)\beta<0, &\tilde r\in\mho_{2>}  
  \ea}\right.  \} ,
\ee} and  $  S _> $ defined by \eqref{21-I4-01-contour},   $  \Gamma  _{5j}= \Gamma  _{5j}(\beta,X,\lambda), \ j=0,1,2$ defined by
\be\label{E:21-I-5-2-new}
\begin{split}
\Gamma_{50}=& 
   \{\xi_0 : 
    \tilde s\in  \mho_0 ,\, \textit{$\tau_0=\tau_{0,\dagger}$} \}\cup S_{50},    
  \\
\Gamma _{51}  =&  \Gamma_{51,out}\cup S_{51}\cup\Gamma_{51,in},    
  \\
\Gamma _{52}  =&  \Gamma_{52,out}\cup S_{52}\cup\Gamma_{52,in},  \end{split}
\ee with
\begin{center}
\begin{tabular}{ll} 
\multirow{2}{*}{$\quad\ S_{50}=  
 \left\{
{\ba{l}
\{\xi_1 : 
      \tilde s_0 =\frac 12^+  \} \\  
   {\phi,}   
  \ea}\right.$} &$Type\,\mathfrak A '\wedge(\widehat{\sigma}=X_1), \, \mathfrak B',\mathfrak C'$,  \\
    &otherwise,\\

    \multirow{3}{*}{$\ \Gamma_{51,in}=  
 \left\{
{\ba{l}  
    {\phi,} \\
  \{\xi_1 : 
    \tilde s\in   \mho_1 ,\, \tau_1=0 \textit{ on }\mho_{1>},   \\
  \hskip.87in\tau_1=\pi \textit{ on }\mho_{1<}, \, \tilde s_1<1/2\},  
  \ea}\right.$} &$ \tilde r_1 >\frac 14$,  \\
    &$ \tilde r_1 <\frac 14$, \\ &\\
\multirow{2}{*}{$\Gamma_{51,out}=  
 \left\{
{\ba{l}  
  {\{\xi_1 : 
    \tilde s\in  \mho_1 ,\, \textit{$\tau_1=\tau_{1,\dagger}$} \},} \\
  {\{\xi_1 : 
    \tilde s\in  \mho_1,\, \textit{$\tau_1=\tau_{1,\dagger}$},\,  \tilde s_1 >1/2 \},}   
  \ea}\right.$}& $  \tilde r_1 >\frac 14$,  \\
    &$ \tilde r_1 <\frac 14$, 
\\ 
    
\multirow{2}{*}{$\quad\ S_{51}=  
 \left\{
{\ba{l}  
   {\phi,} \\
 \{\xi_1 : 
      \tilde s_1 =1/2 \}   
  \ea}\right.$} &$ \tilde r_1 >\frac 14$,  \\
    &$ \tilde r_1 <\frac 14$,\\

    \multirow{3}{*}{$\ \Gamma_{52,in}=  
 \left\{
{\ba{l}  
    {\phi,} \\
  \{\xi_2 : 
    \tilde s\in   \mho_2 ,\, \tau_2=0 \textit{ on }\mho_{2>},   \\
  \hskip.87in\tau_2=\pi \textit{ on }\mho_{2<}, \, \tilde s_2<1/2\},  
  \ea}\right.$} &$ \tilde r_2 >\frac 14$,  \\
    &$ \tilde r_2 <\frac 14$, \\ &\\
 \multirow{2}{*}{$\Gamma_{52,out}=  
 \left\{
{\ba{l}  
  {\{\xi_2 : 
    \tilde s\in  \mho_2 ,\, \textit{$\tau_2=\tau_{2,\dagger}$} \},} \\
  {\{\xi_2 : 
    \tilde s\in  \mho_2,\, \textit{$\tau_2=\tau_{2,\dagger}$},\,  \tilde s_2 >1/2 \},}   
  \ea}\right.$}& $  \tilde r_2 >\frac 14$,  \\
    &$ \tilde r_2 <\frac 14$,\\    
\multirow{2}{*}{$\quad\ S_{52}=  
 \left\{
{\ba{l}  
   {\phi,} \\
 \{\xi_2 : 
      \tilde s_2 =1/2 \}   
  \ea}\right.$} &$ \tilde r_2 >\frac 14$,  \\
    &$ \tilde r_2 <\frac 14$,
 \end{tabular}
\end{center}
and $\alpha_j$,  $\xi_j  $,   $\tau_{ j,\dagger}$, $\mho_j=\mho_j(\tilde\sigma,\beta,X)$  defined by  Definition \ref{D:renormal}.

Using Lemma \ref{L:vekua}, $\tilde r>1$,   and argument as that for $I_4$,
\begin{align}
|I_5|_{L^\infty(D_{\kappa_1})} \le C\epsilon_0|f|_{L^\infty(D_{\kappa_1})}+\sum_{j=1}^2|  \frac { \theta( \tilde r-1)}{2\pi i}  \int_{-\pi}^\pi  d\beta[\partial_ {\beta}   \ln (1-\gamma |\beta| )]\label{E:21-i-5-reduce} \\
\times   \int_{\Gamma  _{5j,in}  }     \frac{  e^{-i \wp(\frac{\xi_j }{{\tilde\sigma}},\beta,X)}f(\frac{ \xi_j }{\tilde\sigma} ,-\beta,X)}{\tilde s_j e^{i\tau_j}-  \tilde r_j e^{ i (\alpha_j -\beta) }} d{\xi_j } |_{L^\infty(D_{\kappa_1})}.\nonumber
\end{align}

Namely, we have to pay extra attention when both $\tilde \zeta$ and $\tilde\lambda$ is close to one of the essential stationary point $\tilde s_{j,\ast}(\tilde\sigma,\beta, X)$, say $\tilde s_{1,\ast}(\tilde\sigma,\beta, X)$ without loss of generality because the other case can be done by analogy.  In this situation,  $\Gamma_{52,in}=\phi$.
  For the estimates on $\Gamma_{51,in}$, in view of 
  \eqref{E:21-I-5-2-new}, we have $|\kappa_1-\tilde\zeta|\ge 1/4$ for  $   \tilde s\in\Gamma_{51,in}  $. Hence  $I_5$ for $\tilde s\in\Gamma_{51,in}$ is no longer a singular integral and we can apply Lemma \ref{L:vekua}. Namely,
\begin{align}
&|\theta(\tilde r-1)\int_{-\pi}^\pi d\beta [\partial_ {\beta}   \ln (1-\gamma |\beta| )] \int_{\Gamma _{51,in }}\frac{  e^{-i \wp(\frac{\xi_1 }{{\tilde\sigma}},\beta,X)}f(\frac{ \xi_1 }{\tilde\sigma} ,-\beta,X)}{\tilde s_1 e^{i\tau_1}-  \tilde r_1 e^{ i (\alpha_1 -\beta) }} d{\xi_1 }|_{L^\infty(D_{\kappa_1})}\label{E:residue-beta} \\
\le &C|\iint_{ \tilde s\in\Gamma_{51,in} }\frac{ \widetilde \gamma_1(\tilde s, \beta) e^{-i\wp(\frac {\tilde s}{{\tilde\sigma}},\beta,X)} f(\frac {\tilde s}{{\tilde\sigma}},-\beta,X)  }{{\tilde\zeta}-\tilde \lambda}d\overline{\tilde\zeta} \wedge d{\tilde\zeta}|_{L^\infty(D_{\kappa_1})}\nonumber\\
\le &    C\epsilon_0 |f|_{L^\infty(D_{\kappa_1})} . \nonumber
\end{align} 

Consequently, 
\be\label{E:I-50-est}
\begin{split}
&|I_{5}|_{L^\infty(D_{\kappa_1})}\le C\epsilon_0|f|_{L^\infty(D_{\kappa_1})}.
\end{split}
\ee

Moreover, via the same argument as for $I_5$, one can prove $I_{\frac52}(\frac{\tilde s }{\tilde\sigma} ,\beta,X)$ is holomorphic in $\tilde s>\frac12$  but may not have $L^\infty$-estimates for the derivatives.
\end{itemize}

\end{proof}

Estimates for Case $(F1)$ are more involved. Let's point out the difficulties and outline our strategy first for clarity.
\begin{itemize}
\item  On $\mho_j$ for $j\ge 1$, in view of Lemma \ref{L:p-estimate}, there are no good uniform estimates (means the coefficients for the polynomial $\tilde s$ only depends on functions of $\beta$)   for $\mathfrak{Re}(-i\wp(\frac{ \tilde s_je^{i\tau_{j,\dagger}}+ \tilde s_{j,\ast}}{X_1}$, $\beta,X)) $. 

We shall take advantage of the scaling invariant properties of the Hilbert transform and estimates    \eqref{E:f2-deform-wp}  and \eqref{E:f3-deform-wp-1}, where $\sqrt{X_2}$ and $\sqrt[3]{X_3}$ are not necessarily equal to $\tilde\sigma=X_1$, to derive estimates on $\Gamma_{4j}$ or $\Gamma_{5j}$, $j\ge 1$.   In view of argument of \eqref{E:residue-beta}, one additional important issue  is   that,  to apply  Lemma \ref{L:vekua}, the Cauchy integral is no longer a singular integral near renormalized essential critical points (see Definition \ref{D:renormal}). This is done in Lemma \ref{L:root-scaling}. 
 
\item On $\mho_0$, the above scaling argument won't work since $0$ is a singular point.

We have to take advantage of estimate \eqref{E:all-0}! If we apply    argument as for Case $(F3)$ or $(F2)$ in Proposition \ref{P:non-homogeneous-F3} directly, then the difficulty encountered is the Jacobian 
\[
\frac 1{|\sin\beta|}d\beta,
\]instead of $\frac 1{\sqrt[3]{|\sin3\beta|}}d\beta$ or $\frac 1{\sqrt {|\sin2\beta|}}d\beta$.  To overcome this difficulty when $|\beta|\ll 1$, 
 we shall squeeze out extra $|\sin\beta|$-decay of integrals on $\Gamma_{40}$ or $\Gamma_{50}$ via  deformation of a finer decomposition. This is done via $\mathfrak J_1$-$\mathfrak J_5$ decomposition \eqref{E:homo-i-41}-\eqref{E:homo-i-45}.
\end{itemize}

\begin{definition}\label{D:renormal} For Case $(F1)$, introduce new scaled $\sigma_j$-parameters on $\mho_j(\beta,X)$, $j=0,1,2,$
\be\label{E:affine-coord-new}
\begin{split}
&\left\{
{\ba{ll}  
{\sigma_0={\tilde\sigma},\  {  \sigma_1}=  \sigma_2=
 \sqrt[3]{|X_3|}},&  \textit{for }\ Type\ \mathfrak A, \mathfrak B,\mathfrak E,\\
 {\sigma_0={\tilde\sigma},\  {  \sigma_1}=  \sigma_2=
 \sqrt[2]{|X_2|}},&  \textit{for }\ Type\ \mathfrak C,\mathfrak D,
 \ea}\right. 
 \end{split}
 \ee and scaled $\sigma_j$-coordinates on $\mho_j(\tilde\sigma, \beta,X)$
 \be\label{E:12-lambda-sigma}
\begin{gathered}
\lambda=\kappa_1+\frac{\tilde r }{\tilde\sigma}e^{i\alpha}=\kappa_1+\frac{   s_{j,\ast}e^{i\beta}+ r_je^{i\alpha_j}}{  \sigma_j}  ,\\
   s_{j,\ast}\equiv \tilde s_{j,\ast}\frac{\sigma_j}{\tilde\sigma},\ r_j=\tilde r_j \frac{  \sigma_j}{{ \tilde \sigma }} \ge 0,\  \alpha_j = \alpha_j(\beta,X,\lambda), \\  
   \frac{\tilde  s}{ \tilde  \sigma } \mapsto  \frac{ \vartheta_j}{\sigma_j}\equiv \frac{s_{j,\ast}+s_je^{i\tau_j}}{\sigma_j}  ,  \\ 
\tilde s\equiv (s_{j,\ast}\pm  s_j)\frac{\tilde\sigma}{\sigma_j}  
  \in \mho_j(\tilde\sigma, \beta,X) 
\end{gathered}
\ee  where $\mho_j(\tilde\sigma, \beta,X)$, $\tilde s_{j,\ast}$, and $\tilde r_j$ are defined by Definition \ref{D:dominant}.
 
\end{definition}
\begin{lemma}\label{L:root-scaling}
For Case $(F1)$, introduce the new scaled $\sigma_j$-coordinates defined by Definition \ref{D:renormal}, 
\be\label{E:F1-s-1-inf}
\inf_\beta s_{1,\ast} = c_0,\quad 0<c_0<1. 
\ee 

\end{lemma}

\begin{proof}\hfill \\
\begin{itemize}
\item [$\blacktriangleright$] {\bf Proof  for $ Type\ \mathfrak A'', \mathfrak B'',\mathfrak E$ :} from Table \ref{Tb:dynamic-FNH}, \eqref{E:affine-coord-new}, 
and
\[
\begin{gathered}
|s_+s_-|=|\frac{X_1}{X_3^{1/3}}\frac{\sin\beta}{\sin3\beta}|\ge \frac 13,\quad\
s_\pm=\frac{    - 1\pm\sqrt{1-\Delta}\, }{3\frac{X_3^{2/3}}{ X_2}\frac{\sin3\beta}{\sin2\beta}},
\end{gathered}
\]we derive $|s_+|\sim |s_-| $ and then \eqref{E:F1-s-1-inf} for $ Type\ \mathfrak A'', \mathfrak B'',\mathfrak E$.

\item [$\blacktriangleright$] {\bf Proof  for $ Type\ \mathfrak C'', \mathfrak D$ :} from \eqref{E:determinant}, Table \ref{Tb:dynamic-FNH}, and \eqref{E:affine-coord-new}, 
\begin{align*}
|  s_ +|=&|\frac{-\Delta}{6\frac{X _3}{\tilde \sigma X_2}\frac{\sin 3\beta}{\sin2\beta} } +\mbox{l.o.t.}|=|\frac{-3\frac{ {X_1}{X_ 3}}{X_2^2}\frac{\sin\beta\sin3\beta}{\sin^22\beta}}{6\frac{X _3 }{ X_2^{3/2}}\frac{\sin 3\beta}{\sin2\beta} } +\mbox{l.o.t.}|\ge\frac 1C,\\
|  s_ -|=&|\frac{-2}{3\frac{X _3}{\tilde \sigma X_2}\frac{\sin 3\beta}{\sin2\beta} } +\mbox{l.o.t.}|\ge\frac 1C|  s_ +|\ge\frac 1C.
\end{align*} Hence \eqref{E:F1-s-1-inf} is proved for Type $\mathfrak C''$ and $\mathfrak D$.
\end{itemize}

\end{proof}
  
\begin{proposition}\label{P:non-homogeneous-F1}
For Case $(F1)$, and $f\in L^\infty (D_{\kappa_j})$ is $\tilde s$-holomorphic, 
\begin{align}
   |I_4| _{C^\mu_{\tilde\sigma}(D_{\kappa_j,\frac 1{\tilde\sigma}})}\le &  C\epsilon_0   |   f|_{L^\infty(D_{\kappa_j})},\label{E:f1-CIE}\\
|I_5| _{L^\infty(D_{\kappa_j )}}\le &  C\epsilon_0   |  f|_{L^\infty(D_{\kappa_j})},\label{E:f1-CIE-1}   
 \end{align}and
 \be\label{E:I-5-2-est-f1}
\begin{split}
&\textit{$I_{\frac52}(\frac{\tilde s }{\tilde\sigma} ,\beta,X)$ is holomorphic in $\tilde s>\frac12$} 
\end{split}
\ee with $I_{\frac{5}{2}}$ defined by \eqref{E:I-5-2}.
\end{proposition}  

\begin{proof}

Thanks to Lemma \ref{L:root-scaling}, we can following the same argument \eqref{E:21-I-4-4-dec-1}-\eqref{E:homo-4},  \eqref{E:21-I-5-dec-1}-\eqref{E:21-I-5-2-new} to obtain 
\begin{align}
|I_{4}|_{C^\mu_{\tilde\sigma}(D_{\kappa_1,\frac 1{\tilde\sigma}})} 
\le &C\epsilon_0 |f|_{L^\infty(D_{\kappa_1})} 
     + \sum_{j=0}^2|I_{4j}|_{L^\infty (D_{\kappa_1})} ,\label{E:renormal-1}\\
     |I_5|_{L^\infty(D_{\kappa_1})} \le & C\epsilon_0|f|_{L^\infty(D_{\kappa_1})}+\sum_{j=0}^2|I_{5j}|_{L^\infty (D_{\kappa_1})},\label{E:renormal-5-new}
     \end{align}
 with
\begin{align}
    I_{4j}=&  
   \sum_{h= 1,2}  \frac { \theta(1-\tilde r)}{2\pi i}  \int_{-\pi}^\pi  d\beta[\partial_ {\beta}   \ln (1-\gamma |\beta| )]e^{-i(n-1)\beta}
 \int_{\Gamma _{4j}} \frac{ e^{-i \wp(\frac{\tilde s_j}{\tilde\sigma}e^{i\tau_j},\beta,X)}f(\frac{\tilde s_j}{\tilde\sigma}e^{i\tau_j},-\beta,X)} {( \tilde s_je^{i\tau_j}-  \tilde r_j e^{ i (\alpha_j-\beta) })^h  }d\tilde s_j e^{i\tau_j},   \nonumber \\
 I_{5j}=&\frac { \theta( \tilde r-1)}{2\pi i}  \int_{-\pi}^\pi  d\beta[\partial_ {\beta}   \ln (1-\gamma |\beta| )] \int_{\Gamma  _{5j,out}  }    \frac{ e^{-i \wp(\frac{\tilde s_j}{\tilde \sigma}e^{i\tau_j},\beta,X)}f(\frac{\tilde s_j}{\tilde \sigma}e^{i\tau_j},-\beta,X)} {    \tilde s_je^{i\tau_j}- \tilde r_j e^{ i (\alpha_j-\beta) }   }d \tilde s_j e^{i\tau_j}.\nonumber
 \end{align} 
 
\begin{itemize}
\item [$\blacktriangleright$] {\bf Proof  for $I_{41}, I_{42},I_{51},I_{52}$ :} 
 For $j\ge 1$, using the scaling invariant of the Hilbert transform and
 \begin{itemize}
 \item for $\Gamma _{4j}\cap \{s_j<1\}$: applying Lemma \ref{L:root-scaling},  Lemma \ref{L:vekua}; 
 \item for $\Gamma _{4j}\cap \{s_j>1\}$: thanks to the scaling invariant of the Hilbert transform, applying  \eqref{E:f3-deform-wp-1} on $\mho_1,\,\mho_2 $ for $Type\,\mathfrak A'',\,\mathfrak B''\,\mathfrak E$,  and  \eqref{E:f2-deform-wp-1} on $\mho_1,\,\mho_2 $ for $Type\,\mathfrak C'',\,\mathfrak D$, 
 \end{itemize} we have
\begin{align}
&|I_{4j}|_{C^1(D_{\kappa_1,\frac{1}{\tilde\sigma}})}\label{E:41}\\
\le & \sum_{h=0}^1\{|\frac { \theta(1-\tilde r)}{2\pi i}  \int_{-\pi}^\pi  d\beta[\partial_ {\beta}   \ln (1-\gamma |\beta| )] \nonumber\\\times&
 \int_{\Gamma _{4j}\cap \{s_j<1\}} \frac{ e^{-i \wp(\frac{  s_j}{ \sigma}e^{i\tau_j},\beta,X)}\widetilde \chi(\frac{  s_j}{ \sigma}e^{i\tau_j},-\beta,X)} {   (\tilde s_je^{i\tau_j}-  \tilde  r_j e^{ i (\alpha_j-\beta) }  )^h( s_je^{i\tau_j}-    r_j e^{ i (\alpha_j-\beta) }  ) }d  s_j e^{i\tau_j}|_{L^\infty (D_{\kappa_1})}\nonumber \\
+& \sum_{h=0}^1\{|\frac { \theta(1-\tilde r)}{2\pi i}  \int_{-\pi}^\pi  d\beta[\partial_ {\beta}   \ln (1-\gamma |\beta| )] \nonumber\\
\times&\int_{\Gamma _{4j}\cap \{s_j>1\}} \frac{ e^{-i \wp(\frac{  s_j}{ \sigma}e^{i\tau_j},\beta,X)}\widetilde \chi(\frac{  s_j}{ \sigma}e^{i\tau_j},-\beta,X)} {   (\tilde s_je^{i\tau_j}-  \tilde  r_j e^{ i (\alpha_j-\beta) }  )^h( s_je^{i\tau_j}-    r_j e^{ i (\alpha_j-\beta) }  ) }d  s_j e^{i\tau_j}|_{L^\infty(D_{\kappa_1})} \}\nonumber\\
\le &C\epsilon_0|\widetilde \chi|_{L^\infty(D_{\kappa_1})}.\nonumber
\end{align}

 In an entirely similar way, namely, for $j\ge 1$,
\begin{itemize}
 \item for $\Gamma _{5j}\cap \{s_j<1\}$: applying Lemma \ref{L:root-scaling},  Lemma \ref{L:vekua}; 
 \item for $\Gamma _{5j}\cap \{s_j>1\}$: thanks to the scaling invariant of the Hilbert transform, applying  \eqref{E:f3-deform-wp-1} on $\mho_1,\,\mho_2 $ for $Type\,\mathfrak A'',\,\mathfrak B''\,\mathfrak E$,  and  \eqref{E:f2-deform-wp-1} on $\mho_1,\,\mho_2 $ for $Type\,\mathfrak C'',\,\mathfrak D$, 
 \end{itemize} we have
\be\label{E:51}
\begin{split}
&|I_{5j}|_{L^\infty(D_{\kappa_1})}\le C\epsilon_0|f|_{L^\infty(D_{\kappa_1})},\\
&\textit{$I_{  \frac52,j}(\frac{\tilde s }{\tilde\sigma} ,\beta,X)$ are holomorphic in $\tilde s>\frac 12$}.
\end{split}
\ee Here
\begin{align*}
I_{\frac52,j}=&\frac { \theta( \tilde r-\frac 12)}{2\pi i}  \int_{-\pi}^\pi  d\beta[\partial_ {\beta}   \ln (1-\gamma |\beta| )] \int_{\Gamma  _{5j,out}  }    \frac{ e^{-i \wp(\frac{  s_j}{ \sigma}e^{i\tau_j},\beta,X)}f(\frac{  s_j}{ \sigma}e^{i\tau_j},-\beta,X)} {    s_je^{i\tau_j}-    r_j e^{ i (\alpha_j-\beta) }   }d  s_j e^{i\tau_j}.
\end{align*}

As a result, proof for $I_{4j}$, $I_{5j}$ for $n=1,\,j\ge 1$ is done.

\item [$\blacktriangleright$] {\bf Proof  for $I_{40 },I_{50}$ :} We shall use the estimate \eqref{E:all-0}. But we have to deal with the singularity $\sin\beta=0$ caused by the change of variables. We shall squeeze out extra $|\sin\beta|$-decay of integrals on $\Gamma_{40}$ or $\Gamma_{50}$ via a finer decomposition. Namely, decompose
\begin{align}
 I_{40}=&-\frac { \theta(1-\tilde r)}{2\pi i} \int_{0}^\pi  d\beta[\partial_ {\beta}   \ln (1-\gamma |\beta| )] \int_{\Gamma_{40}}   \mathfrak J_{<} ,\label{E:scal-4-linear}
\end{align}
\begin{align}
\mathfrak J_<=&
\left\{
{\ba{ll} 
 \mathfrak J_{1}+\mathfrak J_{2}+\mathfrak J_{3}+\mathfrak J_{4}+\mathfrak J_{5} ,  &\textit{if } 0<\beta<\epsilon_1/8,\\
 0,&\textit{if }  -\epsilon_1/8<\beta<0,\\
 \frac{  e^{-i \wp(\frac{\vartheta_0 }{{ \sigma_0}},\beta,X)}f(\frac{ \vartheta_0 }{ \sigma_0} ,-\beta,X)}{s_0 e^{i\tau_0} -    r_0 e^{ i (\alpha_0 -\beta) }} d{\vartheta_0 }, &\textit{if }  |\beta|>\epsilon_1/8; 
 \ea}\right. \label{E:NH-L}
\end{align} and $\mathfrak J_1,\cdots,\mathfrak J_5$ defined by 
 \begin{align}
 \mathfrak J_1= & \theta(\frac 1{|\sin\beta|}-|\tilde s-\tilde r|)   \frac{ {[e^{-i\wp(\frac{\tilde se^{i\tau }}{X_1 } ,\beta,X)  }-1]  f(\frac{\tilde se^{i\tau }}{X_1 } ,-\beta,X) }}{\tilde s e^{i\tau}- \tilde r  e^{ i (\alpha-\beta ) }}d\tilde se^{i\tau },\label{E:homo-i-41} \\
\mathfrak J_2= & \theta(\frac 1{|\sin\beta|}-|\tilde s-\tilde r|)   \frac{ {[ 1-e^{-i\wp(\frac{\tilde se^{-i\tau }}{X_1 } ,-\beta,X)  } ]  f(\frac{\tilde se^{-i\tau }}{X_1 } ,-\beta,X ) }}{\tilde se^{-i\tau } - \tilde r  e^{ i (\alpha-\beta ) }}d\tilde se^{-i\tau },\label{E:homo-i-42} \\
\mathfrak J_3=& \theta(\frac 1{|\sin\beta|}-|\tilde s-\tilde r|) e^{-i\wp(\frac{\tilde se^{-i\tau }}{X_1 } ,-\beta,X)  }  f(\frac{\tilde se^{-i\tau }}{X_1 } ,-\beta,X)  \label{E:homo-i-43} \\
\times&  [\frac{1  }{\tilde se^{-i\tau } - \tilde r  e^{ i (\alpha-\beta ) }}-\frac{ 1  }{\tilde se^{-i\tau } - \tilde r  e^{ i (\alpha+\beta ) }}]d\tilde se^{-i\tau },\nonumber\\
\mathfrak J_4= &\theta(\frac 1{|\sin\beta|}-|\tilde s-\tilde r|)   e^{-i\wp(\frac{\tilde se^{-i\tau }}{X_1 } ,-\beta,X)  }\frac{ f(\frac{\tilde se^{-i\tau }}{X_1 } ,-\beta,X)- f(\frac{\tilde se^{-i\tau }}{X_1 } ,+\beta,X)   }{\tilde se^{-i\tau } - \tilde r  e^{ i (\alpha+\beta ) }} d\tilde se^{-i\tau },\label{E:homo-i-44}\\     
\mathfrak J_5=& \theta(|\tilde s-\tilde r|-\frac 1{|\sin\beta|}) ( e^{-i\wp(\frac{\tilde se^{i\tau }}{X_1 } ,\beta,X)  } \frac{ {f(\frac{\tilde se^{i\tau } }{X_1 } ,-\beta,X) }}{\tilde se^{i\tau } -\tilde r  e^{ i (\alpha-\beta ) }}d\tilde se^{i\tau } \label{E:homo-i-45} \\
  - &e^{-i\wp(\frac{\tilde se^{-i\tau }}{X_1 } ,-\beta,X)  }  \frac{ {f(\frac{\tilde se^{-i\tau } }{X_1 } ,+\beta,X) }}{\tilde se^{-i\tau } - \tilde r  e^{ i (\alpha+\beta ) }}d\tilde se^{-i\tau }),\nonumber
\end{align} with 
 $\tau $  defined by  
Definition \ref{D:deformation} for $\beta\in[0,\pi]$. 

Let's explain our strategy before providing detailed estimates for $\mathfrak J_1,\cdots,\mathfrak J_5$.

For $\mathfrak J_1$-$\mathfrak J_4$, one can squeeze out extra $|\sin\beta|$-decay of nominators of integrands by taking advantage of difference terms.   
Through the change of variables   \be\label{E:second-scale}
\tilde s\mapsto t=\tilde s|\sin\beta|,  
\ee the cut off function $\theta(\frac 1{|\sin\beta|}-|\tilde s-\tilde r|)$ makes the $t$-domain of $\mathfrak J_1$-$\mathfrak J_4$  compact, thus estimates can be derived. For $\mathfrak J_5$,   we use the scaling invariant of the Hilbert transform to cancel the Jacobian $\frac 1{|\sin\beta|}$, the cut off function $\theta(|\tilde s-\tilde r|-\frac 1{|\sin\beta|})$ making the scaled Hilbert transform $\frac{  1 }  {  t-\tilde r|\sin\beta|\,  }$  bounded,    and the nominator is $t$-exponentially decaying. Hence estimates can be derived.

Precisely, 
\begin{itemize}
\item [$\bullet$]   $\mathfrak J_1,\mathfrak J_2$:  From   the mean value theorem  and  \eqref{E:all-0},
{\begin{align} 
  &|e^{-i\wp(\frac{\tilde se^{\pm i\tau }}{X_1 } ,\pm \beta,X)  } -1|\le C  |  \sin\beta|,\quad\textit{ for $\tilde s\in \mho_0(\beta,X)$, $|\beta|<\frac{\epsilon_1}8$}  .\label{E:homo-i-41-42-aux} 
\end{align}}
\item  [$\bullet$] $\mathfrak J_3$:   From $|\beta|<\epsilon_1/8$,   
\begin{align}
 & |\frac1{ \tilde se^{-i\tau_\dagger }-  \tilde r e^{ i (\alpha-\beta) }} -\frac1{ \tilde se^{-i\tau_\dagger }-  \tilde r e^{ i (\alpha+\beta) }}|\label{E:L-H-53-aux}\\
=&  |\frac{\tilde r e^{i\alpha}2\sin\beta}{(\tilde se^{-i\tau_\dagger }-  \tilde r e^{ i (\alpha-\beta) })(\tilde se^{-i\tau_\dagger }-  \tilde r e^{ i (\alpha+\beta) })}| \le C|\sin\beta| .\nonumber
\end{align}

\item [$\bullet$]   $\mathfrak J_4$:  From $\tilde s$-holomorphic properties  of $f$,  
\begin{align} 
    &|\frac{   f(\frac{\tilde s}{\tilde \sigma}e^{-i\tau_\dagger },-\beta,X)-f(\frac{\tilde s}{\tilde \sigma}e^{-i\tau_\dagger },+\beta,X) }{ \tilde se^{-i\tau_\dagger }-  \tilde r e^{ i (\alpha+\beta) }  }|\label{E:homo-i-44-aux}  \\
\le  &  C|f|_{L^\infty(D_{\kappa_1})}\frac{\tilde s| \sin\beta| }{ |\tilde se^{-i\tau_\dagger }-  \tilde r e^{ i (\alpha+\beta) } | } 
\le    C|f|_{L^\infty(D_{\kappa_1})}| \sin\beta|. \nonumber
\end{align}

\end{itemize} 

Applying  \eqref{E:homo-i-41}-\eqref{E:homo-i-44}, \eqref{E:homo-i-41-42-aux}-\eqref{E:homo-i-44-aux}, and the change of variables \eqref{E:second-scale},  for $j=1,\cdots, 4$,
\begin{align}
& |-\frac { \theta(1-\tilde r)}{2\pi i}  \int_{0}^\pi  d\beta[\partial_ {\beta}   \ln (1-\gamma |\beta| )]   +  \int_{  \Gamma _{40}  }  \mathfrak J_j |_{L^\infty(D_{\kappa_1 })} \label{E:H-I-41}\\
\le & C|f|_{L^\infty(D_{\kappa_1})}|\int_{0}^\pi  d\beta[\partial_ {\beta}   \ln (1-\gamma |\beta| )]     \int_{ 2|sin\beta|  } ^{X_1\delta|\sin\beta|} \theta(1-  \left|t-(\tilde r|\sin\beta|)\, \right|)   dt   |_{L^\infty(D_{\kappa_1})}\nonumber\\
 \le &C\epsilon_0 |f|_{L^\infty(D_{\kappa_1})}.\nonumber
\end{align}

On the other hand, using   \eqref{E:homo-i-45}, $\epsilon_1>0$,  the rescaling \eqref{E:second-scale}, the  scaling invariant property of the Hilbert transform,
    one obtains
\begin{align}
 &|-\frac { \theta(1-\tilde r)}{2\pi i}  \int_{0}^\pi  d\beta[\partial_ {\beta}   \ln (1-\gamma |\beta| )]  \int_{  \Gamma _{40}  }  \mathfrak J_5|_{L^\infty(D_{\kappa_1})} \label{E:H-I-45-new} \\
 \le &   C\epsilon_0|f|_{L^\infty(D_{\kappa_1})}\int_{ 2|sin\beta|  } ^{X_1\delta|\sin\beta|}   \theta(| t-\tilde r|\sin\beta|\,|-1)   \frac{  e^{-t\sin\frac{\epsilon_1}4 } }{ |t-\tilde r|\sin\beta|\,| } dt|_{L^\infty(D_{\kappa_1})}  \nonumber\\
   \le &  C\epsilon_0 |f|_{L^\infty(D_{\kappa_1})}. \nonumber
\end{align}

From \eqref{E:H-I-41} and \eqref{E:H-I-45-new},
\be\label{E:i-I_4}
\begin{split}
&|I_{40}|_{L^\infty(D_{\kappa_1})}\le C\epsilon_0 |f|_{L^\infty(D_{\kappa_1})}.
\end{split}
\ee

In an entirely similar way,
\be\label{E:i-I_5}
\begin{split}
&|I_{50}|_{L^\infty(D_{\kappa_1})}\le C\epsilon_0 |f|_{L^\infty(D_{\kappa_1})},\\
&\textit{$I_{  \frac 52,0}(\frac{\tilde s }{\tilde\sigma} ,\beta,X)$ is holomorphic in $\tilde s>\frac 12$}.
\end{split}
\ee Here
\begin{align*}
I_{\frac52,0}=&\frac { \theta( \tilde r-\frac 12)}{2\pi i}  \int_{-\pi}^\pi  d\beta[\partial_ {\beta}   \ln (1-\gamma |\beta| )] \int_{\Gamma  _{50}  }    \frac{ e^{-i \wp(\frac{ \tilde s_j}{ \sigma}e^{i\tau_j},\beta,X)}\chi(\frac{ \tilde  s_j}{ \sigma}e^{i\tau_j},-\beta,X)} {   \tilde  s_je^{i\tau_j}-  \tilde   r_j e^{ i (\alpha_j-\beta) }   }d \tilde  s_j e^{i\tau_j}.
\end{align*}

Combining with   \eqref{E:renormal-1}-\eqref{E:51},   we complete the proof of this proposition.


\end{itemize}

\end{proof}

\begin{proposition}\label{P:degenerate}
For   $X_3=0$, and $f\in L^\infty (D_{\kappa_j})$ is $\tilde s$-holomorphic, 
\begin{align}
   |I_4| _{C^\mu_{\tilde\sigma}(D_{\kappa_j,\frac 1{\tilde\sigma}})}\le &  C\epsilon_0   |   f|_{L^\infty(D_{\kappa_j})},\label{E:f1-CIE-deg}\\
|I_5| _{L^\infty(D_{\kappa_j )}}\le &  C\epsilon_0   |  f|_{L^\infty(D_{\kappa_j})},\label{E:f1-CIE-1-deg}   
 \end{align}and
 \be\label{E:I-5-2-est-f1-deg}
\begin{split}
&\textit{$I_{\frac52}(\frac{\tilde s }{\tilde\sigma} ,\beta,X)$ is holomorphic in $\tilde s>\frac12$} 
\end{split}
\ee with $I_{\frac{5}{2}}$ defined by \eqref{E:I-5-2}.
\end{proposition}
\begin{proof}
For degenerate  cases, in view of Proposition \ref{P:non-homogeneous-F3}, \ref{P:non-homogeneous-F1}, it suffices to discuss the following cases
\begin{itemize}
\item Quadratic Cases  : $X_3=0,\, X_2>0$;
\item Linear Cases : $X_3=X_2=0,\, X_1>0$.
\end{itemize} It is easy to see, linear cases can be tackled by using the same method of Type $\mathfrak A'$ of Case $(F1)$; quadratic cases can be proved by adapting the argument of Type $\mathfrak B',\mathfrak D$ of Case $(F1)$, $(F2)$. 

For illustration, we provide more details for the quadratic cases and $ X_1 \ge \sqrt{X_2}> 0 $ and $\tilde\sigma=X_1$. 
In this case, the stationary point is $\tilde s_\ast =-\frac{  X_1^2}{2X_2}\frac{\sin\beta}{\sin2\beta}$ which satisfies
\[
\partial_{\widehat s}\wp(\frac{\tilde s_\ast}{X_1},\beta,X)=2\frac{X_2}{X_1^{2 }}\sin2\beta\tilde s_\ast+ \sin\beta=0.
\] The essential stationary points $\tilde s_{j,\ast}=  \tilde s_{j,\ast}(\beta,X)$, $j=0,1$ are defined by
\be\label{E:h-critical-points-21}
\tilde s_{0,\ast}\equiv 0;\quad \tilde s_{1,\ast}\equiv 
\left\{
{\ba{ll}  
 - &  {\tilde s_\ast<0},\\
  {\tilde s_\ast  },&   {\tilde s_\ast>0}, 
 \ea}\right.
\ee  where $-$ means no definition is needed.   
 Given $\epsilon_1<\frac{\pi}{2k}\ll 1$, define $ \mho_j(\beta,X )$  by
\begin{align}
\mho_0  
\equiv&
\left\{
{\ba{ll}  
  {[0, \widehat\sigma\delta]=[0,\frac 12]\cup[\frac 12,\widehat\sigma\delta]=\mho_{0,<}\cup\mho_{0,>}} ,  
&   \hskip.57in    \tilde s_\ast<0, \\
  {[0, \frac 1{2\cos\epsilon_1}\tilde s_{1,\ast} ],} &  \hskip.57in  \tilde s_\ast>0,\\
 \ea}\right.\label{E:beta-interval-21}\\
 \mho_1  
\equiv&
\left\{
{\ba{ll}  
 {\phi,}  &   \,   \tilde s_\ast<0,\\ 
 {[(1-\frac 1{2\cos\epsilon_1})\tilde s_{1,\ast} ,\tilde s_{1,\ast} ]\cup[\tilde s_{1,\ast},\tilde \sigma\delta ]\equiv \mho_{1<}\cup\mho_{1  >},} &  \,  \tilde s_\ast>0,
 \ea}\right.   \nonumber
\end{align} 
Write
\be\label{E:lambda-sigma}
\begin{gathered}
\lambda=\kappa_1+\frac{\tilde r e^{i\alpha}}{\tilde \sigma}=\kappa_1+\frac{ \tilde s_{j,\ast}e^{i\beta}+\tilde r_je^{i\alpha_j}}{\tilde \sigma}  ,\\ 
 \tilde r_j=\tilde r_j(\beta,X,\lambda),\ \alpha_j=\alpha_j(\beta,X,\lambda),\ j=0,1.
\end{gathered}
\ee

Due to  Table \ref{Tb:DNH-phase} and Figure \ref{Fg:signature},  we define
the deformation on $\mho^\sharp_j$:
\be\label{E:21-I-4-0}
\begin{gathered}
\tilde  s \mapsto \xi_j \equiv \tilde s_{j,\ast}+\tilde s_je^{i\tau_j},  \\ \tilde s\equiv \tilde s_{j,\ast}\pm\tilde s_j \in \mho_j \  \textit{if }\ |\tau_j|\lessgtr \frac\pi 2 ,\    \tilde s_j \ge 0,  \end{gathered}
\ee 
where
\begin{table}
{\begin{center}
\vskip.1in
\begin{tabular}{|l|l|}
   \hline
\bf{Case} &  $X_1\ge\sqrt{X_2},\ X_3=0$  \cr
   \hline\hline
 $ \tilde s\in\mho_0$ 
 &$\blacktriangleright \mathfrak{Re}(-i\wp(\frac{\xi_0}{\tilde\sigma},\beta,X))= \frac{X_2}{X_1^2} {\sin2\tau_0\sin2\beta\tilde s(\tilde s-\frac{2\sin\tau_0 }{\sin2\tau_0 }\tilde s_{1,\ast}) } $     \cr  
 \hline   
    $\tilde s\in\mho_1$ 
& $\blacktriangleright \mathfrak{Re}(-i\wp(\frac{  \xi_1}{\tilde\sigma},\beta,X))= \frac{X_2}{X_1^2}\sin2\tau_1\sin2\beta\tilde  s_1^2$  \cr 
   \hline
  
\end{tabular}
\end{center} }
\caption{\small Deformation for $\mathfrak{Re}(-i\wp(\frac{\tilde s}{\tilde\sigma},\beta,X))$ for quadratic cases $X_1\ge\sqrt{X_2},\ X_3=0$ } 
\label{Tb:DNH-phase}
\end{table}
\begin{center}
\begin{tabular}{lll} 
\multirow{6}{*}{$   
 \left\{
{\ba{l} 
 {  \tau_0 \equiv 0  ,} \\
  {  \tau_0 \equiv 0  ,} \\
  { \mp\epsilon_1\lessgtr\tau_0 \lessgtr 0  ,} \\
  {\mp\frac{\epsilon_1 }4\lessgtr\tau_0 \lessgtr 0  ,} \\
  { \pm\epsilon_1\gtrless\tau_0 \gtrless 0  ,} \\
  {\pm\frac{\epsilon_1 }4\gtrless\tau_0 \gtrless 0  ,} 
  \ea}\right.$}&$\textit{for }\sin k\beta\gtrless 0,\ \ |\alpha-\beta|\le \frac{\epsilon_1}2$, &$\ \tilde s\in\mho _{0,<},\tilde s_\ast<0,$\\
    &$ \textit{for }\sin k\beta\gtrless 0,\ \ |\alpha-\beta|\ge \frac{\epsilon_1}2$, &$\ \tilde s\in\mho _{0,<},\tilde s_\ast<0,$\\
    &$\textit{for }\sin k\beta\gtrless 0,\ \ |\alpha-\beta|\le \frac{\epsilon_1}2$, &$\ \tilde s\in\mho _{0,>},\tilde s_\ast<0,$\\
    &$ \textit{for }\sin k\beta\gtrless 0,\ \ |\alpha-\beta|\ge \frac{\epsilon_1}2$, &$\ \tilde s\in\mho _{0,>},\tilde s_\ast<0,$\\ 
    &$\textit{for }\sin k\beta\gtrless 0,\ \ |\alpha-\beta|\le \frac{\epsilon_1}2$, &$\ \tilde s\in\mho _0,\tilde s_\ast>0,$\\
    &$ \textit{for }\sin k\beta\gtrless 0,\ \ |\alpha-\beta|\ge \frac{\epsilon_1}2$, &$\ \tilde s\in\mho _0,\tilde s_\ast>0,$\\   
    \multirow{4}{*}{$   
 \left\{
{\ba{l}  
  { \pm\epsilon_1\gtrless\tau_1 \gtrless 0  ,} \\
  {\pm\frac{\epsilon_1 }4\gtrless\tau_1 \gtrless 0  ,} \\  
 { \pm\pi\gtrless\tau_1 \gtrless \pm\pi\mp\epsilon_1  ,} \\
  {\pm\pi\gtrless\tau_1 \gtrless \pm\pi\mp\frac{\epsilon_1}4,} 
    \ea}\right.$}
  &$\textit{for }\sin k\beta\lessgtr 0,\ \ |\alpha_1-\beta|\le \frac{\epsilon_1}2$, &$\ \tilde s\in\mho _{1>},$\\
    &$ \textit{for }\sin k\beta\lessgtr 0,\ \ |\alpha_1-\beta|\ge \frac{\epsilon_1}2$,  &$\ \tilde s\in\mho _{1>},$\\ 
    &$\textit{for }\sin k\beta\gtrless 0,\ \ ||\alpha_1-\beta|-\pi|\le \frac{\epsilon_1}2$, &$\ \tilde s\in\mho _{1<},$\\
    &$ \textit{for }\sin k\beta\gtrless 0,\ \ ||\alpha_1-\beta|-\pi|\ge \frac{\epsilon_1}2$,  &$\ \tilde s\in\mho _{1<},$
 \\  
 \multirow{6}{*}{$   
 \tau_{0,\dagger}
= \left\{
{\ba{l}  
  { 0   ,} \\
  {0   ,} \\
  { \mp\epsilon_1   ,} \\
  {\mp\frac{\epsilon_1 }4   ,} \\
  { \pm\epsilon_1   ,} \\
  {\pm\frac{\epsilon_1 }4   ,} 
  \ea}\right.$}&$\textit{for }\sin k\beta\gtrless 0,\ \ |\alpha-\beta|\le \frac{\epsilon_1}2$, &$\ \tilde s\in\mho _{0,<},\tilde s_\ast<0,$\\
    &$ \textit{for }\sin k\beta\gtrless 0,\ \ |\alpha-\beta|\ge \frac{\epsilon_1}2$, &$\ \tilde s\in\mho _{0,<},\tilde s_\ast<0,$ \\ 
    &$\textit{for }\sin k\beta\gtrless 0,\ \ |\alpha-\beta|\le \frac{\epsilon_1}2$, &$\ \tilde s\in\mho _{0,>},\tilde s_\ast<0,$\\
    &$ \textit{for }\sin k\beta\gtrless 0,\ \ |\alpha-\beta|\ge \frac{\epsilon_1}2$, &$\ \tilde s\in\mho _{0,>},\tilde s_\ast<0,$ \\
    &$\textit{for }\sin k\beta\gtrless 0,\ \ |\alpha-\beta|\le \frac{\epsilon_1}2$, &$\ \tilde s\in\mho _0,\tilde s_\ast>0,$\\
    &$ \textit{for }\sin k\beta\gtrless 0,\ \ |\alpha-\beta|\ge \frac{\epsilon_1}2$, &$\ \tilde s\in\mho _0,\tilde s_\ast>0,$ \\ 
 \multirow{4}{*}{$\tau_{1,\dagger}\equiv  
 \left\{
{\ba{l} 
{\pm\pi\mp\epsilon_1,}  \\
{\pm\pi\mp\frac{\epsilon_1 }4,} \\ 
     {\mp\epsilon_1,}  \\
       {\mp\frac{\epsilon_1 }4,}   
 \ea}\right.$}
    &$\textit{for } \sin k\beta\gtrless0,\ \ ||\alpha_1-\beta |-\pi|\le \frac{\epsilon_1}2$, &$ \tilde s\in\mho  _{1<},$\\
    &$\textit{for }\sin k\beta\gtrless0,\ \ \left||\alpha_1-\beta |-\pi\right|\ge \frac{\epsilon_1}2$, &$ \tilde s\in\mho  _{1<},$\\
    
    &$\textit{for } \sin k\beta\gtrless0,\ \ |\alpha_1-\beta|\le \frac{\epsilon_1}2$, &$ \tilde s\in\mho  _{1>},$\\
    &$\textit{for }\sin k\beta\gtrless0,\ \ |\alpha_1-\beta|\ge \frac{\epsilon_1}2$, &$ \tilde s\in\mho  _{1>} $.
 \end{tabular}
\end{center}

In view of Table \ref{Tb:DNH-phase}, in particular, we will use the scaling invariant property of the Hilbert transform to derive $L^\infty$-estimates on $\mho_1$.  To this aim, we introduce
\be\label{E:affine-coord-deg}
\begin{split}
&{\sigma_0={\tilde\sigma},\  {  \sigma_1} =
 \sqrt[2]{ X_2 }},  
 \end{split}
 \ee and scaled $\sigma_j$-coordinates on $\mho_j(\tilde\sigma, \beta,X)$
\begin{gather*}
\lambda=\kappa_1+\frac{\tilde r }{\tilde\sigma}e^{i\alpha}=\kappa_1+\frac{   s_{j,\ast}e^{i\beta}+ r_je^{i\alpha_j}}{  \sigma_j}  ,\\
s_{0,\ast}=0,\  s_{1,\ast}\equiv \tilde s_{ \ast}\frac{\sigma_1}{\tilde\sigma}=-\frac{  X_1^2}{2X_2}\frac{\sin\beta}{\sin2\beta}\times \frac{\sqrt{X_2}}{X_1}\ge\frac 14,\\ 
r_j=\tilde r_j \frac{  \sigma_j}{{ \tilde \sigma }} \ge 0,\  \alpha_j = \alpha_j(\beta,X,\lambda), \\  
   \frac{\tilde  s}{ \tilde  \sigma } \mapsto  \frac{ \vartheta_j}{\sigma_j}\equiv \frac{s_{j,\ast}+s_je^{i\tau_j}}{\sigma_j}  ,  \\ 
\tilde s\equiv (s_{j,\ast}\pm  s_j)\frac{\tilde\sigma}{\sigma_j}  
  \in \mho_j(\tilde\sigma, \beta,X). 
\end{gather*}
Hence   following  argument as in the proof of Proposition \ref{P:non-homogeneous-F1}, one can complete the proof.
   
\end{proof}

\subsubsection{Proof of Theorem \ref{T:basic}}\label{SS:CIO-k-2} 
  
\begin{proof}
For simplicity and WLOG,   we only give a proof  assuming $l=0$, 
 and reduce the proof to estimating principal parts by Lemma \ref{L:vekua}.  
Define   ${\tilde\sigma}$ and rescaled coordinates   by Definition \ref{D:phase},  \eqref{E:scale-homo}, and decompose the principal part of the CIO  by \eqref{E:scale-mu}-\eqref{E:scal-5}. Thanks to 
 \be\label{E:kappa-scale}
|\partial_{\tilde s}\wp(\frac{\tilde s}{ {\tilde\sigma}}  ,\beta,X  ) |\le C, \qquad \forall\  \tilde s  <2,
\ee
estimates on compact domains  $I_1$-$I_3$ can be derived via the same approach of the proof of Proposition \ref{P:l-c-mu}. As for $I_4$ and $I_5$, estimates has been derived by Proposition \ref{P:non-homogeneous-F3}, \ref{P:non-homogeneous-F1}, and \ref{P:degenerate}. Consequently, Theorem \ref{T:basic} is justified for $n=0,1$. 

To complete the proof, it is sufficient to prove
\begin{itemize}
\item [$\blacktriangleright$] {\bf Proof of $F(3)$ and $(F2)$ when $n>1$ and $f=\widetilde\chi$ :}  Let
\begin{align}
\widetilde\chi^{[n]}=&(\mathcal C  T E_{\kappa_1})^n \widetilde \chi= \mathcal C  T E_{\kappa_1}  \widetilde\chi^{[n-1]},\quad \widetilde\chi^{[0]}=\widetilde\chi,\label{E:CT-n}\\
I_1^{[n]}=&-\frac {\theta(1-\tilde r)}{2\pi i}\iint_{ \tilde s<2} \frac{  \widetilde \gamma_j(\tilde s, \beta)  \widetilde\chi^{[n-1], \flat}(\frac {\tilde s}{ {\tilde\sigma} },-\beta,X) }{\tilde\zeta-\tilde \lambda}d\overline{\tilde\zeta} \wedge d\tilde\zeta  ,\nonumber\\
I_2^{[n]}=&-\frac {\theta(1-\tilde r)}{2\pi i}\iint_{ \tilde s<2} \frac{  \widetilde \gamma_j(\tilde s, \beta)[ e^{-i\wp(\frac {\tilde s}{ {\tilde\sigma} },\beta,X)}-1]\widetilde\chi^{[n-1],\flat}(\frac {\tilde s}{ {\tilde\sigma} },-\beta,X) }{\tilde\zeta-\tilde \lambda}d\overline{\tilde\zeta} \wedge d\tilde\zeta ,\nonumber\\
I_3^{[n]}=&-\frac {\theta(1-\tilde r)}{2\pi i}\iint_{ \tilde s<2} \frac{  \widetilde \gamma_j(\tilde s, \beta) e^{-i\wp(\frac {\tilde s}{ {\tilde\sigma} },\beta,X)}\widetilde\chi^{[n-1], \sharp}(\frac {\tilde s}{ {\tilde\sigma} },-\beta,X) }{{\tilde\zeta}-\tilde \lambda}d\overline{\tilde\zeta} \wedge d{\tilde\zeta} ,\nonumber\\
I_4^{[n]}=&-\frac {\theta(1-\tilde r)}{2\pi i}\iint_{2< \tilde s< {\tilde\sigma} \delta} \frac{ \widetilde \gamma_1(\tilde s, \beta) e^{ -i\wp(\frac { \tilde s}{{\tilde\sigma}},\beta,X)}\widetilde\chi^{[n-1]}(\frac {\tilde s}{{\tilde\sigma}},-\beta,X)  }{{\tilde\zeta}-\tilde \lambda}d\overline{\tilde\zeta} \wedge d{\tilde\zeta},\nonumber\\
I_5^{[n]}=&-\frac {\theta( \tilde r-1)}{2\pi i}\iint_{ \tilde s< {\tilde\sigma} \delta} \frac{ \widetilde \gamma_1(\tilde s, \beta) e^{ -i\wp(\frac { \tilde s}{{\tilde\sigma}},\beta,X)}\widetilde\chi^{[n-1]}(\frac {\tilde s}{{\tilde\sigma}},-\beta,X)  }{{\tilde\zeta}-\tilde \lambda}d\overline{\tilde\zeta} \wedge d{\tilde\zeta}.\nonumber
\end{align} and 
\be\label{E:I-5-n-2}
I_{\frac 52}^{[n]}= -\frac {\theta( \tilde r-\frac 12)}{2\pi i}\iint_{ \tilde s< {\tilde\sigma} \delta} \frac{ \widetilde \gamma_1(\tilde s, \beta) e^{ -i\wp(\frac { \tilde s}{{\tilde\sigma}},\beta,X)}\widetilde\chi^{[n-1]}(\frac {\tilde s}{{\tilde\sigma}},-\beta,X)  }{{\tilde\zeta}-\tilde \lambda}d\overline{\tilde\zeta} \wedge d{\tilde\zeta}.
\ee

  Estimates for $|I^{[n]}_j|_{C^\mu_{\tilde\sigma}(D_{\kappa_1,\frac 1{\tilde\sigma}})}$, $1\le j\le 3$ follow directly from the argument of  $I_1$-$I_3$ and the induction hypothesis on $I_1^{[n-1]}$-$I_4^{[n-1]}$.

  
On the other hand,
\begin{align}
I_4^{[n]}=&-\frac {\theta(1-\tilde r)}{2\pi i}\iint_{2< \tilde s< {\tilde\sigma} \delta} \frac{ \widetilde \gamma_1(\tilde s, \beta) e^{ -i\wp(\frac { \tilde s}{{\tilde\sigma}},\beta,X)}I_5^{[n-1]}(\frac {\tilde s}{{\tilde\sigma}},-\beta,X)  }{{\tilde\zeta}-\tilde \lambda}d\overline{\tilde\zeta} \wedge d{\tilde\zeta}\label{E:i-4-n-dec}  
\end{align}   
Hence from \eqref{E:I-50-est}, \eqref{E:I-5-2-est},  \eqref{E:CT-n}, the induction hypothesis, the meromorphic property of the Cauchy kernel, and the residue theorem,
\begin{align}
  &| I_4^{[n]}|_{C^\mu_{\tilde\sigma}(D_{\kappa_1,\frac{1}{\tilde\sigma}})}\label{E:4-L-infty-mero} \\
  =& | 
  -\frac { \theta(1-\tilde r)}{2\pi i}  \int_{-\pi}^\pi  d\beta[\partial_ {\beta}   \ln (1-\gamma |\beta| )]\nonumber\\
 \times & (  \int_{  S _<}+ \sum_{j=0}^2 \int_{  \Gamma _{4j}  }+\int_{   S  _>  })\frac{  e^{-i \wp(\frac{\tilde s}{\tilde \sigma}e^{i\tau}, \beta,X)  }\widetilde\chi^{[n-1]}(\frac{\tilde s}{\tilde \sigma}e^{i\tau},- \beta,X)}{\tilde se^{i\tau}-  \tilde r e^{ i (\alpha-\beta) }} d\tilde s e^{i\tau}  |_{C^\mu_{\tilde\sigma}(D_{\kappa_1,\frac{1}{\tilde\sigma}})}  \nonumber\\
\le & C\epsilon_0 | \widetilde   \chi^{[n-1]}|_{L^\infty(D_{\kappa_1})}.\nonumber
 \end{align} 
 

Moreover, decompose
\begin{align}
I_5^{[n]}=&-\frac {\theta(\tilde r-1)}{2\pi i}\iint_{0< \tilde s< \frac 12} \frac{ \widetilde \gamma_1(\tilde s, \beta) e^{ -i\wp(\frac { \tilde s}{{\tilde\sigma}},\beta,X)}\widetilde\chi^{[n-1]}(\frac {\tilde s}{{\tilde\sigma}},-\beta,X)  }{{\tilde\zeta}-\tilde \lambda}d\overline{\tilde\zeta} \wedge d{\tilde\zeta}\label{E:i-5-n-dec}  \\
-& \frac {\theta(\tilde r-1)}{2\pi i}\iint_{\frac 12< \tilde s< {\tilde\sigma} \delta} \frac{ \widetilde \gamma_1(\tilde s, \beta) e^{ -i\wp(\frac { \tilde s}{{\tilde\sigma}},\beta,X)}\widetilde\chi^{[n-1]}(\frac {\tilde s}{{\tilde\sigma}},-\beta,X)  }{{\tilde\zeta}-\tilde \lambda}d\overline{\tilde\zeta} \wedge d{\tilde\zeta}\nonumber\\
\equiv&I_{5,<\frac 12}^{[n]}+I_{5,>\frac 12}^{[n]}.\nonumber
\end{align} 

Thanks to $\tilde r>1$, Lemma \ref{L:vekua}, and induction,
\be\label{E:I-5-21-est}
\begin{split}
&|I_{  5,<\frac 12}^{[n]}|_{L^\infty(D_{\kappa_1})}\le C\epsilon_0|    \widetilde\chi^{[n-1]}|_{L^\infty(D_{\kappa_1})},\\
&\textit{$I_{\frac 52,<\frac 12}^{[n]}(\frac{\tilde s }{\tilde\sigma} ,\beta,X)$ is holomorphic in $\tilde s>\frac 12$},
\end{split}
\ee with 
\be\label{E:I-51-n-2}
I_{\frac 52,<\frac 12}^{[n]}= -\frac {\theta( \tilde r-\frac 12)}{2\pi i}\iint_{ 0<\tilde s< \frac12} \frac{ \widetilde \gamma_1(\tilde s, \beta) e^{ -i\wp(\frac { \tilde s}{{\tilde\sigma}},\beta,X)}\widetilde\chi^{[n-1]}(\frac {\tilde s}{{\tilde\sigma}},-\beta,X)  }{{\tilde\zeta}-\tilde \lambda}d\overline{\tilde\zeta} \wedge d{\tilde\zeta}.
\ee

In a similar way as for $I_4^{[n]}$, 
\begin{align}
  &| I_{5,>\frac 12}^{[n]}|_{L^\infty(D_{\kappa_1})}\label{E:5-L-infty-mero} \\
   =&  | -\frac { \theta( \tilde r-1)}{2\pi i}  \int_{-\pi}^\pi  d\beta[\partial_ {\beta}   \ln (1-\gamma |\beta| )] \nonumber\\
 \times & (  \sum_{j=0}^2 \int_{  \Gamma _{5j}  }+\int_{   S  _>  })\frac{  e^{-i \wp(\frac{\tilde s}{\tilde \sigma}e^{i\tau}, \beta,X)  }\widetilde\chi^{[n-1]}(\frac{\tilde s}{\tilde \sigma}e^{i\tau},- \beta,X)}{\tilde se^{i\tau}-  \tilde r e^{ i (\alpha-\beta) }} d\tilde s e^{i\tau}    \nonumber\\
-&\frac {\theta( \tilde r-1)\theta( \tilde r_1-\frac 14)\theta( \tilde r_2-\frac 14)}{ 2\pi i } \nonumber\\
\times&\int_{\beta\in\mathfrak\Delta(\lambda)} d\beta  [\partial_ {\beta}   \ln (1-\gamma |\beta| )]   
e^{-i \wp(\frac{\tilde re^{i(\alpha-\beta )}}{\tilde \sigma} , \beta,X)  } \widetilde\chi^{[n-1]}(\tilde r,\alpha-2\beta,X)|_{L^\infty(D_{\kappa_1})} \nonumber  \\
\le &  C\epsilon_0 |\widetilde\chi^{[n-1]}|_{L^\infty(D_{\kappa_1})},\nonumber
 \end{align} where $\mathfrak\Delta(\lambda)$ is defined by \eqref{E:alpha-beta}, and
\be\label{E:H-CT2-4}
\begin{split}
&\textit{$I_{\frac 52,>\frac 12}^{[n]}(\frac{\tilde s }{\tilde\sigma} ,\beta,X)$ is holomorphic in $\tilde s>\frac 12$},
\end{split}
\ee with
\be\label{E:I-52-n-2}
I_{\frac 52,>\frac 12}^{[n]}= -\frac {\theta( \tilde r-\frac 12)}{2\pi i}\iint_{ \frac 12<\tilde s< {\tilde\sigma} \delta} \frac{ \widetilde \gamma_1(\tilde s, \beta) e^{ -i\wp(\frac { \tilde s}{{\tilde\sigma}},\beta,X)}\widetilde\chi^{[n-1]}(\frac {\tilde s}{{\tilde\sigma}},-\beta,X)  }{{\tilde\zeta}-\tilde \lambda}d\overline{\tilde\zeta} \wedge d{\tilde\zeta}.
\ee

Combining \eqref{E:I-5-21-est} and \eqref{E:H-CT2-4}, we prove
\be\label{E:5-n}
\begin{split}
& |I_{5}^{[n]}|_{L^\infty(D_{\kappa_1})}\le   C\epsilon_0 |    \widetilde\chi^{[n-1]}|_{L^\infty(D_{\kappa_1})},\\
&\textit{$I_{\frac 52}^{[n]}(\frac{\tilde s }{\tilde\sigma} ,\beta,X)$ is holomorphic in $\tilde s>\frac 12$} 
\end{split}
\ee and \eqref{E:32-CIE} is justified for $n\ge 2$ in case of $(F3)$ or $(F2)$.

\item [$\blacktriangleright$] {\bf Proof   for $(F1)$ when $n>1$ and $f=\widetilde \chi$ :} For for $n> 1$, we can   adapt the above argument  for Case $(F2)$ and $(F3)$ and techniques in Proposition \ref{P:non-homogeneous-F1}, the scaling invariant property of the Hilbert transform on $\mho_1,\mho_2$ and a finer decomposition on $\mho_0$, to derive estimates. The only exceptional situation is that,   we have used the $\tilde s$-holomorphic property  of $f$ to sqeeze out extra $|\sin\beta|$-decay for $\mathfrak J_4$ on $\mho_0$. To generalize the argument, it is sufficient to justify, for $\tilde s\ge\frac 12$, \begin{multline} \label{E:gamma}
\hskip.6in|\mathcal C_{\kappa_1+\frac{\tilde s}{\tilde \sigma}e^{-i(\tau_\dagger+\beta)}}TE_{\kappa_1}   g- \mathcal C_{\kappa_1+\frac{\tilde s}{\tilde \sigma}e^{-i(\tau_\dagger-\beta)}} TE_{\kappa_1} g|\\
\le  C\epsilon_0 |g|_{L^\infty(D_{\kappa_1})}\max\{\tilde s|\sin\beta|,|\tilde s\sin\beta |^\mu\} .\hskip.6in
\end{multline}

In the remaining part of the proof, we justify \eqref{E:gamma}.

  $\underline{\emph{(Proof for \eqref{E:gamma})}}:$  Decompose $ \mathcal CTE_{\kappa_1}g=\sum_{j=1}^5I_j$ via  \eqref{E:scale-mu}-\eqref{E:scal-5}  with $\theta(\tilde r-1)$ and $\theta(1-\tilde r)$ replaced by $\theta(\tilde s-\frac 12)$ and $\theta(\frac 12-\tilde s)$.
 Thanks to $\tilde s>\frac 12$, to prove \eqref{E:gamma},  it amounts to the investigation of $
 |  I_5(\frac{\tilde s}{\tilde \sigma}e^{-i\tau_\dagger},-\beta,X)- I_5(\frac{\tilde s}{\tilde \sigma}e^{-i\tau_\dagger},+\beta,X)  |$.

Write
\begin{align}
&  \theta(1/4 -| \tilde s\sin\beta|)[I_5(\frac{\tilde s}{\tilde \sigma}e^{-i\tau_\dagger},-\beta,X)- I_5(\frac{\tilde s}{\tilde \sigma}e^{-i\tau_\dagger},+\beta,X)] \label{E:ind-<}\\
= &  -\frac { \theta(\tilde s-1)\theta(1/4-| \tilde s\sin\beta|)}{2\pi i}  \int_{-\pi}^\pi  d\beta'[\partial_ {\beta'}   \ln (1-\gamma |\beta'| )] 
  \nonumber\\
  \times& \int_0^{\tilde \sigma\delta}\theta(1/2-|\tilde s'-\tilde se^{ i (-\tau_\dagger+\beta-\beta') }|)e^{-i \wp(\frac{\tilde s '}{{\tilde\sigma}},\beta',X)}  [\frac{g(\frac{\tilde s'}{\tilde\sigma} ,-\beta',X)}{  \tilde s '-  \tilde s e^{ i (-\tau_\dagger-\beta-\beta') } }-\frac{g(\frac{\tilde s'}{\tilde\sigma} ,-\beta',X)}{   \tilde s' -  \tilde s e^{ i (-\tau_\dagger+\beta-\beta') } }]  d\tilde s '\nonumber\\
  -& \frac { \theta(\tilde s-1)\theta(1/4-| \tilde s\sin\beta|)}{2\pi i}  \int_{-\pi}^\pi  d\beta'[\partial_ {\beta'}   \ln (1-\gamma |\beta'| )] 
  \nonumber\\
  \times&\int_0^{\tilde \sigma\delta}\theta(|\tilde s'-\tilde se^{ i (-\tau_\dagger+\beta-\beta') }|-1/2)e^{-i \wp(\frac{\tilde s' }{{\tilde\sigma}},\beta',X)}   \frac{-2i\tilde se^{i(-\tau_\dagger-\beta')}\sin\beta g(\frac{\tilde s'}{\tilde\sigma} ,-\beta',X)}{ (\tilde s' -  \tilde s e^{ i (-\tau_\dagger-\beta-\beta') })( \tilde s' -  \tilde s e^{ i (-\tau_\dagger+\beta-\beta') })}  d\tilde s'  \nonumber\\
\equiv &A_1+A_2.\nonumber  
\end{align}

In view of $|\tilde s\sin\beta|<1/4$,
\[
|\frac{\theta( |\tilde s'-\tilde se^{ i (-\tau_\dagger+\beta-\beta') }|-1/2)}{ (\tilde s' -  \tilde s e^{ i (-\tau_\dagger-\beta-\beta') })( \tilde s' -  \tilde s e^{ i (-\tau_\dagger+\beta-\beta') })}|\le \frac{C}{1+(\tilde s')^2},
\]  one obtains
\[
 | A_2|_{L^\infty(D_{\kappa_1})}\le C\epsilon_0 \tilde s|\sin\beta |\,|g|_{L^\infty(D_{\kappa_1})}.
 \]

From   $| \tilde s\sin\beta|<1/4$, and $|\tilde s'-\tilde se^{ i (-\tau_\dagger+\beta-\beta') }|<1/2$,  one can apply   Lemma \ref{L:vekua} to derive  
\[
  |A_1|_{L^\infty(D_{\kappa_1})}\le C\epsilon_0|g|_{L^\infty(D_{\kappa_1})}|\tilde s \sin\beta|^\mu.
\]

Plugging the above estimates into \eqref{E:ind-<}, we obtain
\begin{multline}\label{E:ind-<-est}
\hskip.5in \theta(1/4-| \tilde s\sin\beta|)|  I_5(\frac{\tilde s}{\tilde \sigma}e^{-i\tau_\dagger},-\beta,X)- I_5(\frac{\tilde s}{\tilde \sigma}e^{-i\tau_\dagger},+\beta,X)  |\\
\le C\epsilon_0 |g|_{L^\infty(D_{\kappa_1})} |\tilde s\sin\beta |^\mu .\hskip.5in
\end{multline}

On the other hand, write
\begin{align}
& \theta(| \tilde s\sin\beta|-1/4)( I_5(\frac{\tilde s}{\tilde \sigma}e^{-i\tau_\dagger},-\beta,X)- I_5(\frac{\tilde s}{\tilde \sigma}e^{-i\tau_\dagger},+\beta,X) )\label{E:ind->-dec}\\
=& 2i\tilde s\sin\beta\frac { \theta(\tilde s-1)\theta(| \tilde s\sin\beta|-1/4)}{2\pi i}  \int_{-\pi}^\pi  d\beta'[\partial_ {\beta'}   \ln (1-\gamma |\beta'| )]  \int_0^{\tilde\sigma\delta}e^{-i \wp(\frac{\tilde s' }{{\tilde\sigma}},\beta',X)}g(\frac{\tilde s'}{\tilde\sigma} ,-\beta',X)
 \nonumber\\
  \times&\theta(1/8-|\tilde s' -  \tilde s e^{ i (-\tau_\dagger-\beta-\beta') }| )    \frac{e^{i(-\tau_\dagger-\beta')}}{ (\tilde s' -  \tilde s e^{ i (-\tau_\dagger-\beta-\beta') })( \tilde s' -  \tilde s e^{ i (-\tau_\dagger+\beta-\beta') })}  d\tilde s'\nonumber\\
 + & 2i\tilde s\sin\beta\frac { \theta(\tilde s-1)\theta(| \tilde s\sin\beta|-1/4)}{2\pi i}  \int_{-\pi}^\pi  d\beta'[\partial_ {\beta'}   \ln (1-\gamma |\beta'| )]  \int_0^{\tilde\sigma\delta}e^{-i \wp(\frac{\tilde s' }{{\tilde\sigma}},\beta',X)}g(\frac{\tilde s'}{\tilde\sigma} ,-\beta',X)
 \nonumber\\
  \times&   \theta(1/8-|\tilde s' -  \tilde s e^{ i (-\tau_\dagger +\beta-\beta') }| )    \frac{e^{i(-\tau_\dagger-\beta')}}{ (\tilde s '-  \tilde s e^{ i (-\tau_\dagger-\beta-\beta') })( \tilde s' -  \tilde s e^{ i (-\tau_\dagger+\beta-\beta') })}  d\tilde s'\nonumber\\
  + & 2i\tilde s\sin\beta\frac { \theta(\tilde s-1)\theta(| \tilde s\sin\beta|-1/4)}{2\pi i}  \int_{-\pi}^\pi  d\beta'[\partial_ {\beta'}   \ln (1-\gamma |\beta'| )]  \int_0^{\tilde\sigma\delta}e^{-i \wp(\frac{\tilde s' }{{\tilde\sigma}},\beta',X)}g(\frac{\tilde s'}{\tilde\sigma} ,-\beta',X)
 \nonumber\\
  \times&  [1-\theta(1/8-|\tilde s '-  \tilde s e^{ i (-\tau_\dagger-\beta-\beta') }| )- \theta(1/8-|\tilde s' -  \tilde s e^{ i (-\tau_\dagger+\beta-\beta') }| ) ]\nonumber\\
  \times &   \frac{e^{i(-\tau_\dagger-\beta')}}{ (\tilde s' -  \tilde s e^{ i (-\tau_\dagger-\beta-\beta') })( \tilde s' -  \tilde s e^{ i (-\tau_\dagger+\beta-\beta') })}  d\tilde s'   \nonumber \\
  \equiv &A_1'+A_2'+A_3'.\nonumber
\end{align} 

Thanks to $\tilde s>\frac 12$, $| \tilde s\sin\beta|>1/4$, $|\tilde s' -  \tilde s e^{ i (-\tau_\dagger-\beta-\beta') }| <1/8$, and   Lemma \ref{L:vekua},   one has 
\be\label{E:ind->-12}
 |A'_1|\le C\epsilon_0 |g|_{L^\infty(D_{\kappa_1})} |\tilde s\sin\beta |  .
\ee

Analogously, 
\be\label{E:ind->-2}
 |A'_2|\le C\epsilon_0 |g|_{L^\infty(D_{\kappa_1})} |\tilde s\sin\beta |  .
\ee

Since 
\[
 |\frac{1-\theta(1/8-|\tilde s' -  \tilde s e^{ i (-\tau_\dagger-\beta-\beta') }| )- \theta(1/8-|\tilde s' -  \tilde s e^{ i (-\tau_\dagger+\beta-\beta') }| )}{ (\tilde s' -  \tilde s e^{ i (-\tau_\dagger-\beta-\beta') })( \tilde s' -  \tilde s e^{ i (-\tau_\dagger+\beta-\beta') })}|\le \frac{C}{1+(\tilde s')^2}  ,\] we have
\be\label{E:ind->-12-new}
 |A'_3|\le C\epsilon_0 |g|_{L^\infty(D_{\kappa_1})} |\tilde s\sin\beta |  .
\ee

Therefore, \eqref{E:gamma} is justified by \eqref{E:ind-<-est}-\eqref{E:ind->-12-new}. 

\item [$\blacktriangleright$] Thanks to Proposition \ref{P:degenerate}, for  degenerated quadratic and linear cases, estimates can be derived by adapting argument in previous steps.
\end{itemize}

\end{proof}

\vskip.2in
\subsection{Estimates for the Cauchy integral operator $\mathcal CT$}\label{SS:CIO}

After proving  Theorem \ref{T:basic}, estimates of the Cauchy integral operator reduce to studies of the CIO near $z_n$ and $\infty$.

  Due to the simple pole property of $\phi \in W$ at $z_n$, one can adapt the approach in Subsection \ref{SS:CIO-j} to derive estimates of the CIO near $z_n$. 
\begin{proposition}\label{P:hilbert} {\bf(Estimates near $z_n$)}  Suppose  $\mathcal S=(\{z_n\},\{\kappa_j\}, \mathcal D,s _c)$ is  a $d$-admissible   scattering data, $\phi(x,\lambda)=\frac{\phi_{z_n,\res}(x)}{\lambda-z_n}+\phi_{z_n,r}(x,\lambda) $, $\partial^l_x \phi_{z_n,\res},\, \partial^l_x \phi_{z_n,r}\in L^\infty$ for $0\le l_1+2l_2+3l_3\le d+5$. Then
\[\ba{r} 
    \sum_{0\le l_1+2l_2+3l_3\le d+5} |\partial^l_x\mathcal CTE_{z_n}\phi  |_ {L^\infty } +\sum_{j=1}^M \sum_{0\le l_1+2l_2+3l_3\le d+5} |\partial^l_x\mathcal CTE_{z_n}\phi  |_ {C^\mu(D_{\kappa_j}) } \\
\le   C\epsilon_0\sum_{0\le l_1+2l_2+3l_3\le d+5}(|\partial^l_x \phi_{z_n,\res}|_{L^\infty}+| \partial^l_x \phi_{z_n,r}|_{L^\infty ( D_{z_n})} ). 
\ea\] 
\end{proposition}

\begin{proof} Since proofs for $   \partial^l_{x}  \mathcal C  T E_{z_n}   \phi   $ are identical. We only prove $     \mathcal C  T E_{z_n} \phi    $ for simplicity. Set  $z=z_n$ in \eqref{E:scale-homo} and identify  $\phi(x,\zeta)=\phi(s,\beta,X)$ with $X$ defined by \eqref{E:phase}. 

Firstly, thanks to Lemma \ref{L:vekua}, one has
\be\label{E:c-mu}
   |(1-E_{z_n}) \mathcal CTE_{z_n}\phi  |_ {C^\mu  } \\
\le   C\epsilon_0 (|  \phi_{z_n,\res}|_{L^\infty}+|  \phi_{z_n,r}|_{L^\infty ( D_{z_n})} )
,\ee and it remains to consider $|\partial^l_x\mathcal CTE_{z_n}\phi  |_ {L^\infty ( D_{z_n})}$. To this aim, from \eqref{E:intro-s-c-N}, \eqref{E:cauchy-operator}, we   decompose
\begin{align*}
&  \mathcal C  TE_{z_n}\phi    \\
=&  -\frac { 1}{ 2\pi i}\iint_{ D_{z_n, {\tilde\sigma}\delta}} \frac{  \sgn(\beta) \hbar_n( \frac{\tilde s}{\tilde\sigma},\beta) e^{-i\wp(\frac {\tilde s}{ {\tilde\sigma} },\beta,X)}\phi_{z_n,\res}(x) }{ (\tilde\zeta-\tilde \lambda)(\overline{\tilde\zeta}-z_n)  }d\overline{\tilde\zeta} \wedge d\tilde\zeta +\mathcal C_\lambda E_{z_n}T\phi_{z_n,r} \\
\equiv & II_1+II_2+II_3+II_4+II_5, 
\end{align*}
where
\begin{align*}
II_1=&-\frac {\theta(1-\tilde r)}{2\pi i}\iint_{ \tilde s<2} \frac{  \sgn(\beta) \hbar_n(\frac{\tilde s}{\tilde\sigma}, \beta)  \phi_{z_n,\res}(x) }{(\tilde\zeta-\tilde \lambda)(\overline{\tilde\zeta}-z_n)}d\overline{\tilde\zeta} \wedge d\tilde\zeta  ,\\
II_2=&-\frac {\theta(1-\tilde r)}{2\pi i}\iint_{ \tilde s<2} \frac{  \sgn(\beta) \hbar_n(\frac{\tilde s}{\tilde\sigma}, \beta)[ e^{-i\wp(\frac {\tilde s}{ {\tilde\sigma} },\beta,X)}-1]\phi_{z_n,\res}( x) }{(\tilde\zeta-\tilde \lambda) (\overline{\tilde\zeta}-z_n) }d\overline{\tilde\zeta} \wedge d\tilde\zeta ,\\
II_3=&\ \mathcal C_\lambda E_{z_n}T\phi_{z_n,r} ,\\
II_4=&-\frac {\theta(1-\tilde r)}{2\pi i}\iint_{2< \tilde s< {\tilde\sigma} \delta} \frac{  \sgn(\beta) \hbar_n(\frac{\tilde s}{\tilde\sigma}, \beta)  e^{-i\wp(\frac {\tilde s}{ {\tilde\sigma} },\beta,X)} \phi_{z_n,\res}(x) }{(\tilde\zeta-\tilde \lambda) (\overline{\tilde\zeta}-z_n) }d\overline{\tilde\zeta} \wedge d{\tilde\zeta},\\
II_5=&-\frac {\theta(\tilde r-1)}{2\pi i}\iint_{  \tilde s< {\tilde\sigma} \delta} \frac{  \sgn(\beta) \hbar_n(\frac{\tilde s}{\tilde\sigma}, \beta)  e^{-i\wp(\frac {\tilde s}{ {\tilde\sigma} },\beta,X)} \phi_{z_n,\res}(x) }{(\tilde\zeta-\tilde \lambda) (\overline{\tilde\zeta}-z_n) }d\overline{\tilde\zeta} \wedge d{\tilde\zeta}. 
\end{align*}

From $|\hbar_{z_n}|_{C^1(D_{z_n})}<\epsilon_0$  (assured by the $d$-admissible condition),  the mean value theorem,  and the Hilbert transform theory \cite{Ga66},   
\[
\begin{split}
|II_1|_{L^\infty( D_{z_n})},\,|II_2|_{L^\infty( D_{z_n})}\le & C\epsilon_0 |\phi_{z_n,\res}|_{L^\infty} .
\end{split}
\] 

Writing $\hbar_n(\zeta)=\hbar_n(z_n)+ [\hbar_n(\zeta)-\hbar_n(z_n)]$ and  decompose $II_4=II_{41}+II_{42}$, $II_5=II_{51}+II_{52}$ with
\begin{align*}
II_{41}=&-\frac {\theta(1-\tilde r)}{2\pi i}\iint_{2< \tilde s< {\tilde\sigma} \delta} \frac{  \sgn(\beta) \hbar_{z_n}({z_n})  e^{-i\wp(\frac {\tilde s}{ {\tilde\sigma} },\beta,X)} \phi_{{z_n},\res}(x) }{(\tilde\zeta-\tilde \lambda)  (\overline{\tilde\zeta}-z_n)  }d\overline{\tilde\zeta} \wedge d{\tilde\zeta},\\
II_{51}=&-\frac {\theta(\tilde r-1)}{2\pi i}\iint_{  \tilde s< {\tilde\sigma} \delta} \frac{  \sgn(\beta) \hbar_{z_n}({z_n})  e^{-i\wp(\frac {\tilde s}{ {\tilde\sigma} },\beta,X)} \phi_{{z_n},\res}(x) }{(\tilde\zeta-\tilde \lambda)  (\overline{\tilde\zeta} -{z_n}) }d\overline{\tilde\zeta} \wedge d{\tilde\zeta}.
\end{align*}
  Thanks to $\tilde s$-meromorphic properties  and adapting argument for estimating $I_4$, $I_5$ in $\S$ \ref{SSS:FN-homo} for $ II_{41} $, $ II_{51} $, we obtain
\[
\begin{split}
|II_{41}|_{L^\infty(D_{z_n})},\,|II_{51}|_{ L^\infty(D_{z_n})}\le & C\epsilon_0|\phi_{{z_n},\res}|_{L^\infty}.
\end{split}
\]

For the remaining terms, by Lemma \ref{L:vekua} and $|\hbar_{z_n}|_{C^1(D_{z_n})}<\epsilon_0$  (assured by the $d$-admissible condition),  
\[
\begin{gathered}
|II_{42}|_{L^\infty(D_{z_n})},\,|II_{52}|_{ L^\infty(D_{z_n})}\le   C\epsilon_0|\phi_{{z_n},\res}|_{L^\infty},\\
|II_3|_{L^\infty(D_{z_n})} \le    C\epsilon_0 |\phi_{{z_n},r}|_{L^\infty(D_{z_n})}.
\end{gathered}
\] 
\end{proof}

Estimates of the CIO near $\infty$ has been done by Wickerhauser via Fourier analysis \cite{W87}. We sketch the proof as follows for convenience.
\begin{proposition}\label{P:wickerhauser} {\bf(Estimates near $\infty$)} If  $\mathcal S=(\{z_n\},\{\kappa_j\}, \mathcal D,s _c)$ is $d$-admissible, $(1-\sum_{n=1}^NE_{z_n}-\sum_{j=1}^ME_{\kappa_j} )\partial^l_x \phi \in L^\infty$ for $0\le l_1+2l_2+3l_3\le d+5$. Then
\[\ba{rl}
 &\sum_{0\le l_1+2l_2+3l_3\le d+5}|\partial_x^l\mathcal CT(1-\sum_{n=1}^NE_{z_n}-\sum_{j=1}^ME_{\kappa_j} )\phi|_{L^\infty }\\
 +&\sum_{0\le l_1+2l_2+3l_3\le d+5}\sum_{j=1}^M|\partial_x^l\mathcal CT(1-\sum_{n=1}^NE_{z_n}-\sum_{j=1}^ME_{\kappa_j} )\phi|_{C^\mu(D_{\kappa_j})}\\
   \le & C\epsilon_0\sum_{0\le l_1+2l_2+3l_3\le d+5}|(1-\sum_{n=1}^NE_{z_n}-\sum_{j=1}^ME_{\kappa_j} )\partial_x^l\phi|_{L^\infty }. 
 \ea\] 
\end{proposition}
\begin{proof} Thanks to $d$-admissibility,
\be\label{E:d-admissible-u}(1-\sum_{j=1}^ME_{{\kappa_j}}  )  \sum_{|l_1+2l_2+3l_3|\le {d+5}}|[\ |\overline\lambda-\lambda|^{l_1} +| \overline\lambda^2-\lambda^2|^{l_2}+| \overline\lambda^3-\lambda^3|^{l_3}\ ] s_c (\lambda)|   _{   L^\infty}<\infty.\ee  Together with   
$\partial_x^l\phi^{(k)}\in W$, $| E_{\kappa_j} \partial^l_{x}  \mathcal C  T E_{\kappa_j}   f | _{W}=| \sum_{l'+l''=l}E_{\kappa_j}   \mathcal C  (\partial^{l'}_{x}T )E_{\kappa_j}  (\partial^{l''}_{x} f) | _{W}$, proofs for $\partial_x^l\mathcal CT(1-\sum_{n=1}^NE_{z_n}-\sum_{j=1}^ME_{\kappa_j} )\phi$ are identical. We only prove $ \mathcal CT(1-\sum_{n=1}^NE_{z_n}-\sum_{j=1}^ME_{\kappa_j} )\phi$ for simplicity.

Via the coordinates change, 
\be\label{E:wicker-1}
\begin{gathered}
{ 2\pi i\xi= \overline\zeta-\zeta},\quad 2\pi i\eta=\overline\zeta^2-\zeta^2,\quad{ \frac{d\overline\zeta\wedge d\zeta}{\zeta-\lambda}=\frac{-2\pi^2\sgn(\xi)d\xi d\eta}{p_\lambda(\xi,\eta)}},\\
p_\lambda(\xi,\eta)=(2\pi\xi)^ 2-4\pi i\xi\lambda+2\pi i \eta , \quad \Omega_\lambda=\{(\xi,\eta)\in\mathbf R^2\ :\ |p_\lambda(\xi,\eta)|<1\},\\
\left|\frac 1{p_\lambda}\right|_{L^1(\Omega_\lambda, d\xi d\eta)}\le C, \quad  
\left|\frac 1{p_\lambda}\right|_{L^2(\Omega_\lambda^c,d\xi d\eta)}\le C.
\end{gathered}\ee
 \cite[Lemma 6.III]{W87}, we have
\begin{align}
 &|\mathcal CT(1-\sum_{n=1}^NE_{z_n}-\sum_{j=1}^ME_{\kappa_j})\phi|_{L^\infty(D_{\kappa_j}) } \nonumber\\
\le&C|(1-\sum_{n=1}^NE_{z_n}-\sum_{j=1}^ME_{\kappa_j})\phi|_{L^\infty}|\iint \frac{(1-\sum_{n=1}^NE_{z_n}-\sum_{j=1}^ME_{\kappa_j})|s_c(\zeta)| }{|p_\lambda (\xi,\eta)|}d\xi  d\eta|_{L^\infty(D_{\kappa_j})}\nonumber\\
\le&C|(1-\sum_{n=1}^NE_{z_n}-\sum_{j=1}^ME_{\kappa_j})\phi|_{L^\infty}
\{ |(1-\sum_{n=1}^NE_{z_n}-\sum_{j=1}^ME_{\kappa_j})s_c|_{L^2(d\xi d\eta)}\left|\frac 1{p_\lambda}\right|_{L^2(\Omega_\lambda^c,d\xi d\eta)}\nonumber\\
+&|(1-\sum_{n=1}^NE_{z_n}-\sum_{j=1}^ME_{\kappa_j})s_c|_{L^\infty(d\xi d\eta)}\left|\frac 1{p_\lambda}\right|_{L^1(\Omega_\lambda ,d\xi d\eta)}\}\nonumber\\
\le& C\epsilon_0|(1-\sum_{n=1}^NE_{z_n}-\sum_{j=1}^ME_{\kappa_j})\phi|_{L^\infty}\nonumber
\end{align}
 
Besides, $|\mathcal CT(1-\sum_{n=1}^NE_{z_n}-\sum_{j=1}^ME_{\kappa_j})\phi|_{C^\mu(D_{\kappa_j}) } \le C\epsilon_0|(1-\sum_{n=1}^NE_{z_n}-\sum_{j=1}^ME_{\kappa_j})\phi|_{L^\infty}$ can be proved applying Lemma \ref{L:vekua}, \eqref{E:d-admissible-u}, and 
$\partial_x^l\phi^{(k)}\in W$.

\end{proof}

We conclude Section \ref{S:CIO-0} by estimates from Theorem \ref{T:basic}, Proposition \ref{P:hilbert}, and \ref{P:wickerhauser}.
\begin{theorem}\label{T:CIO}  Suppose the scattering data  $\mathcal S=(\{z_n\},\{\kappa_j\}, \mathcal D,s _c)$ is $d$-admissible,  $ \partial^l_{x}\phi\in W$,  and $ \partial^l_{x}\phi$ are $\lambda$-holomorphic on $D_{\kappa_j}$   for $0\le l_1+2l_2+3l_3\le d+5$. Then
\[
\sum_{0\le l_1+2l_2+3l_3\le d+5}|\partial_x^l(\mathcal CT)^n \phi|_{W }\le (C\epsilon_0)^n\sum_{0\le l_1+2l_2+3l_3\le d+5}|\partial_x^l\phi|_{W}.  
\]
\end{theorem}

\section{Existence of the eigenfunction}\label{S:eigenfunction}

In this section, 
 we shall prove the unique solvability of the system of the CIE \eqref{E:intro-CIE} and the $\mathcal D$-symmetry   \eqref{E:intro-sym} by constructing the iteration sequence
\be\label{E:recursion-iteration-N}
\begin{split}
 \phi^{(k)}(x ,\lambda)  
\equiv&  1+\sum_{n=1}^N\frac{\phi^{(k)}_{z_n,\res} (x   )}{\lambda -z_n}  +\mathcal CT \phi^{(k-1)}(x, \lambda)  ,\ \ k>0,\\
\phi^{(0)}(x ,\lambda)\equiv &\widetilde\chi,
 \end{split}\ee
 such that $\partial_x^l\phi^{(k)}\in W$, $0\le l_1+2l_2+3l_3\le d+5$, 
 \be\label{E:recursion-iteration-D}
(e^{\kappa_1x_1+\kappa_1^2x_2+\kappa_1^3x_3}\phi^{(k)}(x,\kappa^+_1), \cdots,e^{\kappa_Mx_1+\kappa_M^2x_2+\kappa_M^3x_3}\phi^{(k)}(x,\kappa^+_M))\mathcal D=0,\ \ k>0,
\ee  and proving the  convergence of $\phi^{(k)}$ in $W$. 
 
In the following proposition, via the $d$-admissible condition, the $\mathcal D$-symmetry, Sato theory, and Theorem \ref{T:CIO}, we provide formula of the residues $\phi^{(k)}_{z_n,\res}$ and derive their estimates.
\begin{proposition}\label{P:N-alg-sym}
Suppose $ {\mathcal S}=(\{z_n\},\{\kappa_j\}, \mathcal D,s _c)$ is   $d$-admissible and $\phi^{(k)}$, $\phi^{(k)}_{z_n,\res}$ satisfy \eqref{E:recursion-iteration-N}, \eqref{E:recursion-iteration-D}. Then for $k>0$, 
\be \label{E:N-alg-sym}
 {\left(
\ba{c}
\phi^{(k)}_{z_1,\res}  \\
\vdots\\
\phi^{(k)}_{z_N,\res}  
\ea
\right)= -B^{-1}  \widetilde A
\left(
\ba{c}
1+\mathcal C_{{ \kappa_1^+}}T\phi^{(k-1)}\\
\vdots\\
\vdots\\
\vdots\\
\vdots\\
1+\mathcal C_{{ \kappa_M^+}}T\phi^{(k-1)}
\ea
\right),} \ee where
{\tiny\be\label{E:sym-BA}
\begin{split}
\widetilde A= &  \left(
\ba{cccccc}
\kappa_1^Ne^{\theta_1}&\cdots&0&\mathcal D_{N+1,1}e^{\theta_{N+1}}&\cdots&\mathcal D_{M,1}e^{\theta_M} \\
\vdots&\ddots&\vdots&\vdots&\ddots&\vdots \\
0&\cdots&\kappa_N^Ne^{\theta_N}&\mathcal D_{N+1,N}e^{\theta_{N+1}}&\cdots&\mathcal D_{M,N}e^{\theta_M} \ea 
\right)  ,\ \    
B=     \widetilde A 
   \left(
\ba{ccc} 
\frac 1{\kappa_1-z_1}&\cdots&\frac 1{\kappa_1-z_N}\\
\vdots&\ddots&\vdots\\
\vdots&\ddots&\vdots\\
\vdots&\ddots&\vdots\\
\vdots&\ddots&\vdots\\
\frac 1{\kappa_M-z_1}&\cdots&\frac 1{\kappa_M-z_N} 
\ea
\right), 
\end{split}
\ee} and $ e^{\theta_j}= e^{\kappa_j x_1+\kappa_j^2 x_2+\kappa_j^3 x_3}$.  

Moreover,    for $k>0$,
 \begin{align}
 \sum_{0\le l_1+2l_2+3l_3\le d+5}\left|\partial_x^l\phi^{(k)}_{z_n,\res}\right|_{L^\infty}\le   C(1+\epsilon_0\sum_{0\le l_1+2l_2+3l_3\le d+5}&\left|\partial_x^l\phi^{(k-1)} \right|_W)\label{E:k-est-new},\\
 \sum_{0\le l_1+2l_2+3l_3\le d+5}\left|\partial_x^l\left[\phi^{(k)}_{z_n,\res}-\phi^{(k-1)}_{z_n,\res}\right]\right|_{L^\infty}\le & (C\epsilon_0)^{k },  \label{E:k-difference-new}\\
 \sum_{0\le l_1+2l_2+3l_3\le d+5}\left|\partial_x^l\left[\phi^{(k)}_{z_n,\res}-\widetilde\chi_{z_n,\res}\right]\right|_{L^\infty}\le&  C\epsilon_0 .  \label{E:k-difference-0}
 \end{align}

\end{proposition} 
\begin{proof} Write the $\mathcal D$-symmetry and the evaluation at $ \kappa_j ^+ $ of $\phi^{(k)}$ as   a  linear system  for $M+N$ variables $\{\phi^{(k)}(x,\kappa_j^+),\phi^{(k)}_{z_n,\res}(x)\}$, 
{\small\begin{align}
&\left(
\ba{ccccccccc}
\kappa_1^Ne^{\theta_1}&\cdots&0&\mathcal D_{N+1,1}e^{\theta_{N+1}}&\cdots&\mathcal D_{M,1}e^{\theta_M}&0&\cdots&0\\
\vdots&\ddots&\vdots&\vdots&\ddots&\vdots&0&\ddots&0\\
0&\cdots&\kappa_N^Ne^{\theta_N}&\mathcal D_{N+1,N}e^{\theta_{N+1}}&\cdots&\mathcal D_{M,N}e^{\theta_M}&0&\cdots&0\\
-1&\cdots&0&0&\cdots&0&\frac 1{\kappa_1-z_1}&\cdots&\frac 1{\kappa_1-z_N}\\
\vdots&\ddots&\vdots&\vdots&\ddots&\vdots&\vdots&\ddots&\vdots\\
\vdots&\ddots&\vdots&\vdots&\ddots&\vdots&\vdots&\ddots&\vdots\\
\vdots&\ddots&\vdots&\vdots&\ddots&\vdots&\vdots&\ddots&\vdots\\
\vdots&\ddots&\vdots&\vdots&\ddots&\vdots&\vdots&\ddots&\vdots\\
0&\cdots&0&0&\cdots&-1&\frac 1{\kappa_M-z_1}&\cdots&\frac 1{\kappa_M-z_N}
\ea
\right) 
\left(
\ba{c}
\phi^{(k)}(x,\kappa_1^+)\\
\vdots\\
\vdots\\
\vdots\\
\vdots\\
\phi^{(k)}(x,\kappa_M^+)\\
\phi^{(k)}_{z_1,\res} (x)\\
\vdots\\
\phi^{(k)}_{z_N,\res} (x)
\ea
\right)\label{E:N-system}\\
&=
\left(
\ba{c}
0\\
\vdots\\
0\\
-1-\mathcal C_{\kappa_1^+}T\phi^{(k-1)}\\
\vdots\\
\vdots\\
\vdots\\
\vdots\\
-1-\mathcal C_{\kappa_M^+}T\phi^{(k-1)}
\ea
\right)  .\nonumber
 \end{align}} Solving $\phi^{(k)}(x,\kappa_j^+)$ in terms of $\phi^{(k)}_{z_n,\res}(x)$ and plugging the outcomes into \eqref{E:N-system} again yields 
{\small \be\label{E:bdd-res} 
 \begin{split}
&B\left(
\ba{c}
\phi^{(k)}_{z_1,\res}(x)\\
\vdots\\
\phi^{(k)}_{z_N,\res}(x)
\ea
\right)=-\widetilde A
\left(
\ba{c}
 1+\mathcal C_{\kappa_1^+}T\phi^{(k-1)}\\
\vdots\\
\vdots\\
\vdots\\
\vdots\\
 1+\mathcal C_{\kappa_M^+}T\phi^{(k-1)}
\ea
\right) ,
\end{split}\ee} with  $B$ and $\widetilde A $ defined by \eqref{E:sym-BA}.   By the $d$-admissible condition,    the   system \eqref{E:N-system} is just determined and is equivalent to \eqref{E:N-alg-sym}.
 
Proofs of \eqref{E:k-est-new}-\eqref{E:k-difference-0} rely heavily on Sato theory.    The key observation is that  the leading term of the right hand side of \eqref{E:N-alg-sym} can be realized by residue of some normalized Sato eigenfunction $\widetilde\chi'$. Through the direct computation of $\widetilde\chi'$ and matching with \eqref{E:line-grassmannian}, \eqref{E:chi-intro}, we shall prove boundedness of $B^{-1}A$ and then verify \eqref{E:k-est-new}-\eqref{E:k-difference-0} by multi linearity and the totally positive (TP) condition. 

Precisely, 
the $d$-admissible condition implies that 
\[
\begin{gathered}
A''= \textit{diag}\,(q_1, \cdots, q_M)^{-1} \times\mathcal D \times  \textit{diag}\,(q_1, \cdots, q_N) =\left(
\ba{ccc}
a''_{11}&\cdots &a''_{1N}\\
\vdots &\ddots &\vdots\\
a''_{N1}&\cdots &a''_{NN}\\
\vdots &\ddots &\vdots\\
a''_{M1}&\cdots &a''_{MN}
\ea
\right),\\
A'=A''\times 
\left(
\ba{ccc}
a''_{11}&\cdots &a''_{1N}\\
\vdots &\ddots &\vdots\\
a''_{N1}&\cdots &a''_{NN}
\ea
\right)^{-1},\qquad
A'\in  \mathrm {Gr}(N,M)_{\color{black}> 0},
\end{gathered}
\]   such that 
setting 
\be\label{E:sato-prime}
\widetilde \chi'(x,\lambda)=\widetilde \chi_{\{z_n\},\{\kappa_j\},A'}(x,\lambda) 
\ee as the normalized Sato eigenfunction with data $(\{z_n\},\{\kappa_j\},A',0)$, one has
\begin{gather} 
{  {  \widetilde\chi'}(x, \lambda) =1+\sum_{n=1}^N\frac{   \widetilde\chi'_{z_n, \res }(x   )}{\lambda -z_n }  ,} \label{E:sato-CIE} \\
(e^{\kappa_1x_1+\kappa_1^2x_2+\kappa_1^3x_3}\widetilde \chi'(x,\kappa _1),\cdots,e^{\kappa_Mx_1+\kappa_M^2x_2+\kappa_M^3x_3}\widetilde \chi'(x,\kappa _M))\mathcal D=0,\label{E:sato-D}
\end{gather}and, from the $d$-admissible condition,   $\forall k$,
\be\label{E:asy-chi'}
|\widetilde\chi'_{z_n, \res }(x   )-\widetilde\chi _{z_n, \res }(x   )|_{C^k}
\le C_k\epsilon_0.
\ee 

Moreover, using $\mathcal D$-symmetry and evaluating $\widetilde\chi'$ at $\kappa_j$ and   the above argument, yields 
\be\label{E:k-diff}
\quad{ \left(
\ba{c}
\widetilde\chi'_{z_1,\res} (x) \\
\vdots\\
\widetilde\chi'_{z_N,\res}(x)  
\ea
\right)= -  B^{-1}  \widetilde {  A}
\left(
\ba{c}
1 \\
\vdots\\  \vdots\\ \vdots\\ \vdots\\
1 
\ea
\right),}
\ee with $B$ and $\widetilde A $ defined by \eqref{E:sym-BA}. Let $E_j =e^{\theta_j}  = e^{\kappa_j x_1+\kappa_j^2 x_2+\kappa_j^3 x_3} 
 $ and write
\begin{align}
\widetilde {  A}= &  
\mathcal D^T \textit{diag}\,( E_1, \cdots,  E_M),\label{E:tilde-B-chi}\\
B
=&   
\mathcal D^T   \textit{diag}\,( E_1, \cdots,  E_M)\left(
\ba{ccc} 
 \frac 1{\kappa_1-z_1 }&\cdots&\frac{1}{\kappa_1-z_N }\\
\vdots&\ddots&\vdots\\
\vdots&\ddots&\vdots\\
\frac 1{\kappa_M-z_1 }&\cdots&\frac{1}{\kappa_M-z_N }
\ea
\right).\nonumber
\end{align} 

 From Sato theory,  \eqref{E:chi-intro},  \eqref{E:intro-sym-N-D-flat}, \eqref{E:k-diff}, \eqref{E:tilde-B-chi}, elementary row and column operations, and matching the   coefficients of $E_1\times\cdots\times E_N$,  
\begin{align}
  &B^{-1} = \frac{1}{\tau'(x)}
\left(
\ba{ccc} 
  b_{11}&\cdots &  b_{1N}\\
\vdots&\ddots&\vdots\\
  b_{N1}&\cdots &  b_{NN} 
\ea
\right),\nonumber
\\
&  {  b_{kl}= \sum_{ J(kl) =( j _{(kl),1},\cdots,j _{(kl),N-1})}\Delta _{J(kl)} E_{J(kl)}(x),\ \ 1\le j _{(kl),1}<\cdots<j _{(kl),N-1}\le M,}\label{E:determinant-B}\\
& \textit{$\tau'(x)$ is the tau function with data $\kappa_j$, $A'$,}\nonumber\\
 & |\Delta _{J(kl)}|=|\Delta _{J(kl)} (\{z_n\},\{\kappa_j\},A')|<C. \nonumber
\end{align}   
Moreover,   
 \begin{align}
   &\tau'(x) \widetilde\chi'_{z_h,\res}(x)\label{E:mathfrak-BA-0}\\
=&\textit{the $h$-row of}\ {\tiny \left(
\ba{ccc} 
 b_{11}&\cdots &  b_{1N}\\
\vdots&\ddots&\vdots\\
  b_{N1}&\cdots &  b_{NN} 
\ea
\right)
\left(
\ba{c}
\kappa_1^NE_1+\cdots+\mathcal D_{N+1,1}E_{N+1}+\cdots+\mathcal D_{M,1}E_M\\
\vdots\\  \vdots\\ \vdots\\ \vdots\\
\kappa_N^NE_{  1}+\cdots+\mathcal D_{N+1,N}E_{N+1}+\cdots+\mathcal D_{M,N}E_M 
\ea
\right)} \nonumber 
\\
=&(\kappa_1^NE_1+\cdots+\mathcal D_{N+1,1}E_{N+1}+\cdots+\mathcal D_{M,1}E_M)\sum_{|J(h1)|=N-1}\Delta _{J(h1)}  E_{J(h1)}(x)\nonumber\\
+&\cdots
+(\kappa_N^NE_{  1}+\cdots+\mathcal D_{N+1,N}E_{N+1}+\cdots+\mathcal D_{M,N}E_M) \sum_{|J(hN)|=N-1}\Delta _{J(hN)}  E_{J(hN)}(x)    
\nonumber\\
\equiv&(\widetilde{  a}_{11}E_1+\cdots+\widetilde{  a}_{1M}E_M)\sum_{|J(h1)|=N-1}\Delta _{J(h1)}  E_{J(h1)}(x)\nonumber\\
+&\cdots
+(\widetilde{ a}_{N1}E_1+\cdots+\widetilde{  a}_{NM}E_M) \sum_{|J(hN)|=N-1}\Delta _{J(hN)}  E_{J(hN)}(x),   
\nonumber
\end{align}and 
\be \label{E:chi-plucker}
\begin{split}
0=& \widetilde {  a}_{1k}E_k \ \sum_{k\in J(h1),\ |J(h1)|=N-1}\Delta _{J(h1)}  E_{J(h1)}(x) 
\ +\ \cdots \\
+&\widetilde {  a}_{Nk}E_k \sum_{k\in J(hN), |J(hN)|=N-1}\Delta _{J(hN)}  E_{J(hN)}(x)  .
\end{split} \ee  

Using \eqref{E:N-alg-sym}, \eqref{E:tilde-B-chi}-\eqref{E:chi-plucker}, and  multi linearity, 
\begin{align*}
&\tau'(x) \phi^{(k)}_{z_h,\res}(x) 
= \tau'(x) \widetilde\chi'_{z_h,\res}(x)\\
+&\textit{the $h$-row of } {\tiny \left(
\ba{ccc} 
 b_{11}&\cdots &  b_{1N}\\
\vdots&\ddots&\vdots\\
  b_{N1}&\cdots &  b_{NN} 
\ea
\right)
\left(
\ba{c}
\widetilde a_{11}E_1\mathcal C_{\kappa_1^+}T\phi^{(k-1)}+\cdots+ \widetilde a_{1M}E_M\mathcal C_{\kappa_M^+}T\phi^{(k-1)}\\
\vdots\\  \vdots\\ \vdots\\ \vdots\\
\widetilde a_{N1}E_1\mathcal C_{\kappa_1^+}T\phi^{(k-1)}+\cdots+ \widetilde a_{NM}E_M \mathcal C_{\kappa_M^+}T\phi^{(k-1)}
\ea
\right)}\\
=&\sum_{|J(h)|=N}\Delta _{J(h)} E_{J(h)}(x),\qquad 
\end{align*} with 
\begin{gather*}
 \sum_{0\le l_1+2l_2+3l_3\le d+5}|\partial_x^l \Delta _{J(h)}|
 <C(1+\sum_{j=1}^M \sum_{0\le l_1+2l_2+3l_3\le d+5}|\partial_x^l\mathcal C_{\kappa_j^+}T\phi^{(k-1)}|).
\end{gather*}
 Along with  the totally positive (TP) condition of $A' $,  yield 
\[
 \sum_{0\le l_1+2l_2+3l_3\le d+5}|\partial_x^l \phi^{(k)}_{z_n,\res}(x)|\le C(1+\sum_{j=1}^M \sum_{0\le l_1+2l_2+3l_3\le d+5}|\partial_x^l\mathcal C_{\kappa_j^+}T\phi^{(k-1)}|), 
\] and \eqref{E:k-est-new}-\eqref{E:k-difference-0} follow  from     and \eqref{E:asy-chi'}.

\end{proof}
Proof of Proposition \ref{P:N-alg-sym} explains that the $\mathcal D$-symmetry is a closed condition of the system \eqref{E:intro-sym}, \eqref{E:intro-CIE}. Besides, it justifies that, unlike the IST for the KdV equation and most $1+1$-dimensional integrable systems, residues and poles are no longer crucial to define the scattering data. Finally,   $L^\infty$ estimates of the  the residues, along with Theorem \ref{T:CIO},  will  induce   $L^\infty$-stability  of $  {\mathrm{Gr}(N, M)_{> 0}}$ KP solitons in the last section.

We are ready to prove \eqref{E:intro-CIE}-\eqref{E:intro-asymp} of Theorem \ref{T:intro-inverse}.   
\begin{theorem}\label{T:phi-k-N} Given a $d$-admissible scattering data $ {\mathcal S}=(\{z_n\},\{\kappa_j\}, \mathcal D,s _c)$, there exists uniquely $m$, $  \partial_x^l m\in W$, $0\le l_1+2l_2+3l_3\le d+5$,
   satisfying \eqref{E:intro-CIE}-\eqref{E:intro-asymp}. In particular, $m$ can be constructed 
   via the iteration sequence \eqref{E:recursion-iteration-N}, \eqref{E:N-alg-sym}, and, for $0\le l_1+2l_2+3l_3\le d+5$,
\be \label{E:m-dx-order}
\begin{gathered}
\left(\partial_{x }^l\phi^{(k)}  \right)_{z_n,\res} 
= \,\partial_{x }^l\phi^{(k)}_{z_n,\res}, \\ 
|\partial _{x }^l\phi^{(k)} -
 \partial _{x }^l\phi^{(k-1)}|_W \le (C\epsilon_0)^k,\quad  
\lim_{k\to\infty} \partial_{x }^l\phi^{(k)}= \,\partial_{x }^l m \   \in W .
 \end{gathered}
\ee
\end{theorem}
\begin{proof} Stipulating residues $\phi_{z_n,\res}^{(k)}(x)$ in \eqref{E:recursion-iteration-N} by \eqref{E:N-alg-sym}, we construct  the iteration sequence $\{\phi^{(k)}\}$ satisfying the $\mathcal D$-symmetry for $k>0$. Applying Theorem \ref{T:CIO} and Proposition \ref{P:N-alg-sym}, we have $\phi^{(k)}\in W$ and \eqref{E:intro-CIE}, \eqref{E:intro-sym}.

Suppose   $m_1,\,m_2\in W$ and satisfy \eqref{E:intro-CIE} and \eqref{E:intro-sym}. Hence
\be\label{E:iteration-induction}
\begin{split}
&m_1(x,\lambda)-m_2(x,\lambda)=\sum_{n=1}^N\frac{m_{1,z_n,\res}(x)-m_{2,z_n,\res}(x)}{\lambda-z_n}+\mathcal CT(m_1-m_2).
\end{split}
\ee
Following argument in the proof of Proposition \ref{P:N-alg-sym} and applying Theorem \ref{T:CIO}, 
\[\small
\left(
\ba{c}
m_{1,z_1,\res}-m_{2,z_1,\res}\\
\vdots\\
m_{1,z_N,\res}-m_{2,z_N,\res}
\ea
\right)
 =-B^{-1}  \widetilde A
\left(
\ba{c}
 \mathcal C_{{ \kappa_1^+}}T(m_1-m_2)\\
\vdots\\
\vdots\\
\vdots\\
\vdots\\
 \mathcal C_{{ \kappa_M^+}}T(m_1-m_2)
\ea
\right),
\] and
\begin{align*}
&| m_1-m_2|_{W}\le   C\epsilon_0 |m_1-m_2|_W.
\end{align*}Thanks to   $\epsilon_0\ll 1$,  $m_1(x,\lambda)\equiv m_2(x,\lambda)$.

To prove \eqref{E:intro-asymp}, it suffices to establish \eqref{E:m-dx-order}. 
Proofs for $   \partial^l_{x} \phi^{(k)}$, $0\le l_1+2l_2+3l_3\le d+5$, are identical. We only prove $\partial_{x_1}\phi^{(k)}$ for simplicity.

Taking $x_1$-derivatives on both sides of the iteration sequence \eqref{E:recursion-iteration-N}, formally
\begin{align}
  \partial_{x_1}\phi^{(k)}   
= & \sum_{n=1}^N\frac{\partial_{x_1}\phi^{(k)}_{z_n,\res}  }{\lambda -z_n}  +\mathcal CT(\overline\zeta-\zeta)  \phi^{(k-1)} +\mathcal CT \partial_{x_1}\phi^{(k-1)}  .\label{E:recursion-dx}
\end{align}

Using $\partial_{x_1}\phi^{(0)}\in W$, the $d$-admissible condition, and   Theorem \ref{T:CIO},  
\be\label{E:ct-0}
\begin{split}
  |\mathcal CT\partial_{x_1}\phi^{(0)}|_W\le  &C\epsilon_0 |\partial _{x_1}\phi^{(0)}|_W ,\quad
   |\mathcal CT(\overline\zeta-\zeta)  \phi^{(0)}|_W\le   C\epsilon_0 |\phi^{(0)}|_W . 
   \end{split}
\ee  

Adapting argument of Proposition \ref{P:N-alg-sym} and applying \eqref{E:ct-0}, we have
\[
|\partial_{x_1}\phi^{(1)}_{z_n,\res}|_{L^\infty}\le C(1+|\mathcal CT(\overline\zeta-\zeta)  \phi^{(0)}|_{L^\infty} +|\mathcal CT \partial_{x_1}\phi^{(0)} |_{L^\infty})<\infty.
\]
Consequently, \eqref{E:recursion-dx} is valid for $k'=1$,
\[
\partial_{x_1}\phi^{(1)}\in W, \quad \left(\partial_{x_1}\phi^{(1)}\right)_{z_n,\res}=\partial_{x_1}\phi^{(1)}_{z_n,\res},\] and 
\[
\begin{split}
&|\partial_{x_1}\phi^{(1)}-\partial_{x_1}\phi^{(0)}|_W 
\le   C( |\mathcal CT\partial_{x_1}\phi^{(0)}|_W+| \mathcal CT(\overline\zeta-\zeta)  \phi^{(0)}|_W)\le C\epsilon_0 (|\phi^{(0)}|_W+|\partial _{x_1}\phi^{(0)}|_W).
\end{split}
\]

Inductively, assume that   \eqref{E:recursion-dx} is valid for all $k'\le k$, 
\be\label{E:resi-iteration}
\begin{gathered}
\partial_{x_1}\phi^{(k')}\in W,\quad \left(\partial_{x_1}\phi^{(k')}\right)_{z_n,\res}=\partial_{x_1}\phi^{(k')}_{z_n,\res},\\
|\partial_{x_1}\phi^{(k')}-\partial_{x_1}\phi^{(k'-1)}|_W\le  (C\epsilon_0)^{k'}(|\phi^{(0)}|_W+|\partial _{x_1}\phi^{(0)}|_W).
\end{gathered}
\ee 
Together with  the $d$-admissible condition,   Theorem \ref{T:CIO}, and argument for proving Proposition \ref{P:N-alg-sym}, yields
\[
  |\mathcal CT\partial_{x_1}\phi^{(k)}|_W<\infty,\quad
   |\mathcal CT(\overline\zeta-\zeta)  \phi^{(k)}|_W<\infty , 
\]and
\begin{align*}
&|\partial_{x_1}\phi^{(k+1)}_{z_n,\res}-\partial_{x_1}\phi^{(k)}_{z_n,\res}|_{L^\infty}\le C|\mathcal CT(\overline\zeta-\zeta)(\phi^{(k)}- \phi^{(k-1)})|_W+|\mathcal CT(\partial_{x_1}\phi^{(k)}-\partial_{x_1}\phi^{(k-1)})|_W\\
\le &C\epsilon_0| \phi^{(k )}- \phi^{(k-1)} |_W+C\epsilon_0| \partial_{x_1}\phi^{(k)}-\partial_{x_1}\phi^{(k-1)} |_W\le  (C\epsilon_0)^{k+1}.
\end{align*}

As a result,   \eqref{E:recursion-dx} is valid for $k'=k+1$ and
\be\label{E:resi-iteration-1}
\begin{gathered}
\partial_{x_1}\phi^{(k+1)}\in W, \quad \left(\partial_{x_1}\phi^{(k+1)}\right)_{z_n,\res}=\partial_{x_1}\phi^{(k+1)}_{z_n,\res},\\
|\partial_{x_1}\phi^{(k+1)}-\partial_{x_1}\phi^{(k)}|_W\le   (C\epsilon_0)^{k+1}.
\end{gathered}\ee

Hence   \eqref{E:m-dx-order} is proved by induction.

\end{proof}



\section{The inverse scattering transform}\label{S:potential}
 
Provided the initial data $u_0$  decays sufficiently fast, applying Fourier theory and unique solvability of the CIE, Wickerhauser  derived  a Lax pair for $\Phi(x,\lambda)= e^{ \lambda  x_1+ \lambda ^2x_2}    m(x, \lambda)$, derive a representation formula for  the potential $u(x)$ of the Lax equation, and obtain Sobolev estimates of  $u(x)$ \cite{W87,W85}. Such a framework for perturbed $\mathrm{Gr}(N,M)_{>0}$ KP solitons is no longer efficient due to   singular structures of the continuous scattering data. We shall use the iteration procedure \eqref{E:recursion-iteration-N}, \eqref{E:N-alg-sym} to derive the Lax equation for $m(x,\lambda)$, construct the inverse scattering transform, and prove the continuity in this section. 

    Firstly, introduce the shorthand notation for the heat operator 
\be \label{E:op-de}
\begin{gathered}
-\partial_{x_2}+\partial_{x_1}^2+2 \lambda\partial_{x_1}=-\nabla_2+\nabla_1^2 ,\\
\nabla_1=\partial_{x_1}+\lambda,\   \nabla_2=\partial_{x_2}+\lambda^2 ,\end{gathered} 
\ee and the iteration sequence
\be\label{E:CIE-J}
\phi^{(k)}=  1+J\phi^{(k)}+\mathcal CT\phi^{(k-1)}, \qquad
    J  \, \phi^{(k)}  = \sum_{n=1}^N\frac{\phi^{(k)}_{z_n,\res}(x)}{\lambda-z_n}. 
\ee   
Therefore, formally
\be\label{E:cauchy}
\begin{split}
  (-\nabla_2+\nabla_1^2  )\phi^{(k)} 
=&  \left[-\nabla_2+\nabla_1^2,
  J  \right] \phi^{(k)}+\left[-\nabla_2+\nabla_1^2,
  \mathcal CT  \right] \phi^{(k-1)}\\
+  &   
  J(-\nabla_2+\nabla_1^2)\phi^{(k)}   +\mathcal CT(-\nabla_2+\nabla_1^2) \phi^{(k-1)}.
\end{split}
\ee
Note $
\left[\nabla_j,\ T\right]=0$ for $ j=1,\,2$, hence 
\begin{align}
& \left[-\nabla_2+\nabla_1^2 ,\mathcal CT \right]  \phi^{(k-1)}= \left[-\nabla_2+\nabla_1^2 ,\mathcal C \right] T \phi^{(k-1)} =2 \left[\lambda,\mathcal C \right] \partial_{x_1}\left(T  \phi^{(k-1)}\right)\label{E:rep-continuous}\\
=&-\frac { i}\pi\partial_{x_1}\iint  T \phi^{(k-1)} \ d\overline\zeta\wedge d\zeta   \nonumber \end{align}
holds, and, in view of \eqref{E:m-dx-order},
\begin{align}
&\left[-\nabla_2+\nabla_1^2 ,  J  \right]\phi^{(k)}\label{E:rep-discrete}\\
=&\left[-\partial_{x_2}+\partial_{x_1}^2+2 \lambda\partial_{x_1},  J  \right]  \phi^{(k)}\nonumber \\
=& (-\partial_{x_2}+\partial_{x_1}^2+2 \lambda\partial_{x_1})\sum_{n=1}^N\frac{\phi^{(k)}_{z_n,\res}(x)}{\lambda-z_n}-\sum_{n=1}^N\frac{\left((-\partial_{x_2}+\partial_{x_1}^2+2 \lambda\partial_{x_1})\phi^{ (k)}\right)_{z_n,\res}(x)}{\lambda-z_n}\nonumber\\
=& (-\partial_{x_2}+\partial_{x_1}^2+2 \lambda\partial_{x_1})\sum_{n=1}^N\frac{\phi^{(k)}_{z_n,\res}(x)}{\lambda-z_n}-\sum_{n=1}^N\frac{ (-\partial_{x_2}+\partial_{x_1}^2+2 \lambda\partial_{x_1})\psi^{ (k)} _{z_n,\res}(x)}{\lambda-z_n}\nonumber \\
 =&  2  \partial_{x_1} \sum_{n=1}^N\phi^{(k)}_{z_n,\res}(x)\nonumber. \nonumber
\end{align}  

 To carry out rigorous analysis of the above computation and investigate the convergence of \eqref{E:rep-continuous}-\eqref{E:rep-discrete} for $d$-admissible scattering data, we provide the following two lemmas.
 {\begin{lemma}\label{L:J} 
   Given a $d$-admissible scattering data $  {\mathcal S}=(\{z_n\},\{\kappa_j\}, \mathcal D,s _c)$, let $m$ be the unique solution to  \eqref{E:intro-CIE}, \eqref{E:intro-sym} which can be constructed via the iteration sequence \eqref{E:recursion-iteration-N}, \eqref{E:N-alg-sym}. Then  $\left[-\nabla_2+\nabla_1^2 ,  J  \right]\phi^{(k)}$ is  independent of $\lambda$ and bounded,  
$J   (-\nabla_2+\nabla_1^2  ) \phi^{(k)}\in W$ satisfying, as $k\to\infty$,
\begin{align}
  \left[-\nabla_2+\nabla_1^2 ,  J  \right]\phi^{(k)}\to \ &\left[-\nabla_2+\nabla_1^2 ,  J  \right]m=2\partial_{x_1}\sum_{n=1}^Nm_{z_n,\res}(x) ,\label{E:J-phi-infty }\\
 J   (-\nabla_2+\nabla_1^2  ) \phi^{(k)}\to \ & J   (-\nabla_2+\nabla_1^2  ) m\quad \textit{in $W$}.\label{E:j-phi-W}
\end{align}Moreover,
\be\label{E:iteration-u}
\sum_{0\le l_1+2l_2+3l_3\le d+4}\left|\partial_x^l\left(\left[-\nabla_2+\nabla_1^2 ,  J  \right]m+u_s\right)\right|_{L^\infty}\le C\epsilon_0 
.
\ee
\end{lemma}
\begin{proof} 

From \eqref{E:rep-discrete}  and Theorem \ref{T:phi-k-N}, we derive \eqref{E:J-phi-infty },  \eqref{E:j-phi-W}, and 
\begin{align}
 &\sum_{0\le l_1+2l_2+3l_3\le d+4}|\partial^l_x\left(\left[-\nabla_2+\nabla_1^2 ,  J  \right]m+u_s\right)|_{L^\infty } \label{E:u-rep-est-new}\\
\le &\sum_{0\le l_1+2l_2+3l_3\le d+4}|\partial^l_x\left(\left[-\nabla_2+\nabla_1^2 ,  J  \right]m-\left[-\nabla_2+\nabla_1^2 ,  J  \right]\phi^{(k)}\right)|_{L^\infty }\nonumber\\
+&\sum_{0\le l_1+2l_2+3l_3\le d+4}|\partial^l_x\left(\left[-\nabla_2+\nabla_1^2 ,  J  \right]\phi^{(k)}-\left[-\nabla_2+\nabla_1^2 ,  J  \right]\phi^{(0)}\right)|_{L^\infty }\nonumber\\
+&\sum_{0\le l_1+2l_2+3l_3\le d+4}|\partial^l_x\left(\left[-\nabla_2+\nabla_1^2 ,  J  \right]\widetilde\chi+u_s\right)|_{L^\infty }\le C \epsilon_0. \nonumber
\end{align}

\end{proof} }

\begin{lemma}\label{L:commutator}If $ {\mathcal S}=(\{z_n\},\{\kappa_j\}, \mathcal D,s _c)$ is $d$-admissible then $\left[-\nabla_2+\nabla_1^2 ,  \mathcal CT  \right]\phi^{(k-1)} $    is independent of $\lambda$, bounded, and
\begin{align}
\left[-\nabla_2+\nabla_1^2 ,  \mathcal CT  \right]\phi^{(k-1)}\to&
\  \left[-\nabla_2+\nabla_1^2 ,  \mathcal CT \right]m= \partial_{x_1}\iint  T m  \, d\overline\zeta\wedge d\zeta 
,\label{E:CT-u-est}\\
 \mathcal CT  (-\nabla_2+\nabla_1^2  ) \phi^{(k)}\to&\ \   \mathcal CT (-\nabla_2+\nabla_1^2  ) m\quad \textit{in $W$}, \label{E:CT-u-est-0} 
\end{align} and
\be
\sum_{0\le l_1+2l_2+3l_3\le d+4}\left|\partial_x^l \left[-\nabla_2+\nabla_1^2 ,  \mathcal CT \right]m \right|_{L^\infty}\le C\epsilon_0 .\label{E:CT-u-est-1}
\ee

\end{lemma}
\begin{proof} We  first consider \eqref{E:CT-u-est}.  In view of \eqref{E:rep-continuous},  formally     
\begin{equation}\label{E:2}
\begin{split}
{ \partial_{x_1}\iint  T \phi^{(k-1)}  \, d\overline\zeta\wedge d\zeta }
  =  \iint  T (\overline\zeta-\zeta)\phi^{(k-1)}  \, d\overline\zeta\wedge d\zeta +\iint  T \partial_{x_1}\phi^{(k-1)}  \, d\overline\zeta\wedge d\zeta. 
\end{split}
\end{equation}

Since that singularities of $T$  and $\phi^{(k)}\in W$  are $\kappa_j$ and $z_n$ respectively which are disjoint. Via Lemma \ref{L:vekua}, Theorem \ref{T:phi-k-N}, and the $d$-admissible condition,
\begin{align}
&|\iint  T (\overline\zeta-\zeta)(\phi^{(k-1)}-\phi^{(k-2)})  \, d\overline\zeta\wedge d\zeta|_{L^\infty}\label{E:2-1}\\
\le& (C\epsilon_0)^k+ |\iint (1-\sum_{n=1}^N E_{z_n}- \sum_{j=1}^M E_{\kappa_j})T (\overline\zeta-\zeta)(\phi^{(k-1)}-\phi^{(k-2)})  \, d\overline\zeta\wedge d\zeta\ |_{L^\infty}.\nonumber
\end{align}

The $d$-admissible condition assures 
\be\label{E:d-admissible-u-1}(1-\sum_{j=1}^ME_{{\kappa_j}}  )  \sum_{|l|\le {3}}\left|[\ |\overline\lambda-\lambda|^{l_1} +| \overline\lambda^2-\lambda^2|^{l_2}\ ] s_c (\lambda)\right|   _{   L^\infty}<\infty.\ee  Together with   
$\phi^{(k)}\in W$,   the change of variables \eqref{E:wicker-1}, 
   Theorem \ref{T:CIO}, and Theorem \ref{T:phi-k-N}, yields
\begin{align}
&|\iint (1-\sum_{n=1}^N E_{z_n}- \sum_{j=1}^M E_{\kappa_j})T (\overline\zeta-\zeta)(\phi^{(k-1)}-\phi^{(k-2)})  \, d\overline\zeta\wedge d\zeta \ |_{L^\infty}\label{E:2-1-infty}\\
\le&C \left|\iint (1-\sum_{n=1}^N E_{z_n}- \sum_{j=1}^M E_{\kappa_j})\left|s_c\times\xi \times (\phi^{(k-1)}-\phi^{(k-2)})  \right|\, \frac{d\xi d\eta}{ |\xi|}\right|_{L^\infty}\nonumber\\
\le&C (C\epsilon_0)^{k-1} |\iint (1-\sum_{n=1}^N E_{z_n}- \sum_{j=1}^M E_{\kappa_j}) s_c    \, d\xi d\eta\ |_{L^\infty}\nonumber\\
\le& C(C\epsilon_0)^k.\nonumber
\end{align}
Consequently,
\be\label{E:2-1-est}
|\iint  T (\overline\zeta-\zeta)(\phi^{(k-1)}-\phi^{(k-2)})  \, d\overline\zeta\wedge d\zeta\ |_{L^\infty}\le  (C\epsilon_0)^k.
\ee

Applying  \eqref{E:wicker-1}, \eqref{E:recursion-iteration-N}, \eqref{E:m-dx-order},
\be\label{E:d-admissible-u-new}  (1-\sum_{j=1}^ME_{{\kappa_j}}  )  \sum_{|l_1 |\le 1}\left|    |\overline\lambda-\lambda|^{l_1}   s_c \right|  _{     L^2(|\lambda_I| d\overline\lambda \wedge d\lambda) \cap L^\infty}  <\infty,
\ee assured by the $d$-admissible condition,  and  Theorem \ref{T:CIO},
\begin{align}
& | \iint  T [\partial_{x_1} \phi^{(k-1)} -\partial_{x_1} \phi^{(k-2)} ]  \, d\overline\zeta\wedge d\zeta\ |_{L^\infty}\label{E:2-2}\\
\le& |\iint  s_{c}(\zeta) e^{ (\overline\zeta-\zeta)x_1+(\overline\zeta^2-\zeta^2)x_2+(\overline\zeta^3-\zeta^3)x_3   }\sum_{n=1}^N\frac{\partial_{x_1}\phi^{(k-1)}_{z_n,\res}(x)-\partial_{x_1}\phi^{(k-2)}_{z_n,\res}(x)}{ \overline\zeta-z_n} \, d\overline\zeta\wedge d\zeta\  |_{L^\infty} \nonumber\\
 +&|\iint  s_{c}(\zeta) e^{ (\overline\zeta-\zeta)x_1+(\overline\zeta^2-\zeta^2)x_2+(\overline\zeta^3-\zeta^3)x_3   }\partial_{x_1}\mathcal CT (\phi^{(k-2)}-\phi^{(k-3)}  )\, d\overline\zeta\wedge d\zeta\  |_{L^\infty} \nonumber\\
\le& (C\epsilon_0)^k+ |\iint  s_{c}(\zeta) e^{ (\overline\zeta-\zeta)x_1+(\overline\zeta^2-\zeta^2)x_2+(\overline\zeta^3-\zeta^3)x_3   }\partial_{x_1}\mathcal CT (\phi^{(k-2)}-\phi^{(k-3)}  )\, d\overline\zeta\wedge d\zeta\  |_{L^\infty}.\nonumber
\end{align}
 
Moreover, by Fubini's theorem,   \eqref{E:d-admissible-u-1}, \eqref{E:d-admissible-u-new}, Theorem \ref{T:CIO}, and an induction,
\begin{align}
&|\iint  s_{c}(\zeta) e^{ (\overline\zeta-\zeta)x_1+(\overline\zeta^2-\zeta^2)x_2+(\overline\zeta^3-\zeta^3)x_3   }\partial_{x_1}\mathcal CT (\phi^{(k-2)}-\phi^{(k-3)})  \, d\overline\zeta\wedge d\zeta \ |_{L^\infty}\label{E:2-2-final}\\
\le &\frac 1{2\pi}\ | \iint   \, d\overline\zeta_1\wedge d\zeta_1 \  s_{c}(\zeta_1) e^{ (\overline\zeta_1-\zeta_1)x_1+(\overline\zeta_1^2-\zeta_1^2)x_2+(\overline\zeta_1^3-\zeta_1^3)x_3   }{(\overline\zeta_1 -\zeta_1 )} \nonumber\\
\times& (\phi^{(k-2)}(x,\overline\zeta_1)-\phi^{(k-3)}(x,\overline\zeta_1))\iint \frac{s_{c}(\zeta )}{ \zeta _1-\overline\zeta}e^{ (\overline\zeta -\zeta )x_1+(\overline\zeta ^2-\zeta^2)x_2+(\overline\zeta ^3-\zeta ^3)x_3   } d\overline\zeta \wedge d\zeta \   |_{L^\infty}\nonumber \\
+&\frac 1{2\pi}\ | \iint   \, d\overline\zeta_1\wedge d\zeta_1 \ { T[\partial_{x_1}\phi^{(k-2)}(x,\overline\zeta_1 )-\partial_{x_1}\phi^{(k-3)}(x,\overline\zeta_1)] }\nonumber\\
\times& \iint \frac{s_{c}(\zeta )}{ \zeta_1 -\overline\zeta}e^{ (\overline\zeta -\zeta )x_1+(\overline\zeta ^2-\zeta^2)x_2+(\overline\zeta ^3-\zeta ^3)x_3   } d\overline\zeta \wedge d\zeta \   |_{L^\infty} \nonumber\\
\le&(C\epsilon_0)^k+\frac 1{2\pi}\ |  \iint \, d\overline\zeta_1\wedge d\zeta_1 \ { T[\partial_{x_1}\phi^{(k-2)} -\partial_{x_1}\phi^{(k-3)} ] }\nonumber\\
\times& \iint \frac{s_{c}(\zeta )}{ \zeta_1 -\overline\zeta}e^{ (\overline\zeta -\zeta )x_1+(\overline\zeta ^2-\zeta^2)x_2+(\overline\zeta ^3-\zeta ^3)x_3   } d\overline\zeta \wedge d\zeta\    |_{L^\infty} \nonumber\\
\le&(C\epsilon_0)^k+C\epsilon_0(C\epsilon_0)^{k-1} \le(C\epsilon_0)^k .  \nonumber
\end{align}

Combining \eqref{E:2-2} and \eqref{E:2-2-final}, we obtain
\be\label{E:2-2-est}
| \iint  T [\partial_{x_1} \phi^{(k-1)} -\partial_{x_1} \phi^{(k-2)} ]  \, d\overline\zeta\wedge d\zeta \ |_{L^\infty}\le  (C\epsilon_0)^k.
\ee 

Hence \eqref{E:CT-u-est} follows from \eqref{E:rep-continuous}, \eqref{E:2}, \eqref{E:2-1-est}, and \eqref{E:2-2-est}. Estimate \eqref{E:CT-u-est-0} can be proved similarly.

To prove \eqref{E:CT-u-est-1}, thanks to  Theorem \ref{T:basic} and Theorem \ref{T:CIO}, it suffices to prove $\lambda\partial_{x_1}\phi^{(k)}\in W$. In view of \eqref{E:recursion-iteration-N} and Theorem \ref{T:CIO}, write 
\begin{align*}
 \lambda\partial_{x_1}\phi^{(k)} 
=&
    \sum_{n=1}^N\frac{\lambda\partial_{x_1}\phi^{(k)}_{z_n,\res} (x   )}{\lambda -z_n}  +\lambda\partial_{x_1}\mathcal CT \phi^{(k-1)}(x, \lambda)\\
    =&
    \sum_{n=1}^N\frac{\lambda\partial_{x_1}\phi^{(k)}_{z_n,\res} (x   )}{\lambda -z_n}  +\frac1{2\pi i}\partial_{x_1}\iint \frac \zeta{\zeta-\lambda}T \phi^{(k-1)}(x, \zeta)d\bar\zeta\wedge d \zeta\\
=&
    \sum_{n=1}^N\frac{\lambda\partial_{x_1}\phi^{(k)}_{z_n,\res} (x   )}{\lambda -z_n}  +\frac1{2\pi i} \iint \frac 1{\zeta-\lambda}T(\zeta-\bar\zeta) \bar \zeta\phi^{(k-1)}(x, \zeta)d\bar\zeta\wedge d \zeta    \\
+&\frac1{2\pi i}\iint \frac 1{\zeta-\lambda}T \bar \zeta\partial_{x_1}\phi^{(k-1)}(x, \zeta)d\bar\zeta\wedge d \zeta .    
\end{align*} By induction, $\phi^{(k)}\in W$ and, in an entirely similar way,
 \[
|\mathcal CT  (-\nabla_2+\nabla_1^2  ) \phi^{(k)}-\mathcal CT  (-\nabla_2+\nabla_1^2  ) \phi^{(k-1)}|_W\le C(C\epsilon_0)^{k+1}.
\]Hence \eqref{E:CT-u-est-0} is justified. 
\end{proof}

Lemma \ref{L:J}, \ref{L:commutator}, and Theorem \ref{T:CIO}   imply that \eqref{E:cauchy}-\eqref{E:rep-discrete} hold  rigorously and converge to
\begin{gather} 
 (-\nabla_2+\nabla_1^2  )m 
=  \left[-\nabla_2+\nabla_1^2,
  J  +
  \mathcal CT  \right] m 
+      
  (J    +\mathcal CT)(-\nabla_2+\nabla_1^2) m \quad\in W
  \label{E:inverse-lax}\\
 -u(x)\equiv     
 \left[-\nabla_2+\nabla_1^2 ,
  J +\mathcal CT\right]  m
=  -\frac i\pi\partial_{x_1}\iint  T  m \ d\overline\zeta\wedge d\zeta  +2  \partial_{x_1}\sum_{n=1}^N m_{z_n,\res}(x) \label{E:lax-t-cauchy}
\end{gather}with 
\be\label{E:u-est}
\sum_{0\le l_1+2l_2+3l_3\le d+4}\ |\partial_x^l\left[u(x)-u_s(x)\right]\,|_{L^\infty}\le C\epsilon_0.
\ee
Combining  with Theorem \ref{T:phi-k-N}, we obtain
\begin{align}
&(-\nabla_2+\nabla_1^2 )m=-(1-J-\mathcal CT)^{-1}u(x)1=-u(x)(1-J-\mathcal CT)^{-1}1=-u(x)m(x,\lambda).\label{E:formal-derivation}
\end{align}

To summarize,
\begin{theorem} \label{T:u-inverse} Given a $d$-admissible scattering data $\mathcal S=(\{z_n\},\{\kappa_j\}, \mathcal D,s _c)$, let $  m(x,\lambda) $ be the solution to \eqref{E:intro-CIE}, \eqref{E:intro-sym} and $u(x)$ be define by \eqref{E:lax-t-cauchy}. 
Hence we prove \eqref{E:intro-Lax-u}-\eqref{E:intro-Lax-u-asymp}. 
\end{theorem}

\begin{definition}\label{D:inverse-scattering-transform}
Given a $d$-admissible scattering data $(\{z_n\},\{\kappa_j\}, \mathcal D,s _c)$, we define the inverse scattering transform $\mathcal S^{-1}$  by \eqref{E:intro-u-rep} (i.e., \eqref{E:lax-t-cauchy}) and
\[
\mathcal S^{-1}( \{z_n,\kappa_j, \mathcal D,s _c(\lambda)\})=-\frac i\pi\partial_{x_1}\iint  T  m \ d\overline\zeta\wedge d\zeta  +2  \partial_{x_1}\sum_{n=1}^N m_{z_n,\res}(x).\] 
\end{definition}

Equipping the topology of scattering data $(\{z_n\},\{\kappa_j\}, \mathcal D,s _c)$  by the norms on the right hand side of  \eqref{E:epsilon-0-BPP-new},  and the topology of $u(x)$ by the $C^l$-norms, we obtain
\begin{theorem}
\label{T:inverse-transf-continuity} The inverse scattering transform 
 $\mathcal S^{-1}$ is continuous at each $d$-admissible scattering data with trivial continuous scattering data $(\{z_n\},\{\kappa_j\}, \mathcal D,s _c)$.
\end{theorem}

\section{The evolution}\label{S:cauchy-kp}

In this section, applying the IST established in previous sections, \cite{Wu20, Wu21}, we shall   
\begin{itemize}
\item justify the evolution equation for $u(x)=\mathcal S^{-1}(\{z_n\},\{\kappa_j\}, \mathcal D,s _c)$ is the KPII equation; 
\item solve the Cauchy problem for perturbed $ {\mathrm{Gr}(N, M)_{> 0}}$ KP solitons;
\item prove a uniqueness theorem; 
\item derive an $L^\infty$-stability theorem of $ {\mathrm{Gr}(N, M)_{> 0}}$ KP solitons.
\end{itemize}

\begin{theorem}\label{T:inverse-KP}
If $ (\{z_n\},\{\kappa_j\}, \mathcal D,s _c)$ is $d$-admissible  then  $u(x)=\mathcal S^{-1}(\{z_n\},\{\kappa_j\}, \mathcal D,s _c) $ maps $\RR\times\RR\times\RR^+$ to $ \RR$ and satisfies  the KPII equation \eqref{E:intro-KP-ist}.

\end{theorem} 
 
\begin{proof}
We are proving the reverse theorem of \cite[Theorem 5]{Wu21}. Letting $\Phi(x, \lambda)= e^{ \lambda  x_1+ \lambda ^2x_2}  m(x, \lambda)$, by the representation formula \eqref{E:intro-Lax-u-asymp}, we define the evolution operators 
\be\label{E:evo-operator}
\begin{split}
\mathcal M =&- \partial_{x_3}+  \partial_{x_1}^3+\frac 32u\partial_{x_1}+\frac 34u_{x_1}+\frac 34\partial_{x_1}^{-1}u_{x_2}+\tau(\lambda) ,\ \
\tau(\lambda)= -\lambda ^3 \\
\mathcal M \Phi(x,\lambda)=&e^{ \lambda  x_1+ \lambda ^2x_2}\left(\mathcal M+3\lambda\partial_{x_1}^2+3\lambda^2\partial_{x_1}+\lambda^3+\frac 32u\lambda\right)m(x,\lambda) \\
\equiv  & e^{ \lambda  x_1+ \lambda ^2x_2}\left(\mathfrak Mm\right)(x,\lambda). 
\end{split}
\ee 

Applying the $\mathcal D$-symmetry, the form of \eqref{E:intro-sym-N-D}, 
 and the linear flow properties of the scattering data \eqref{E:linearization-D-evol}, we obtain, 
\begin{align*}
 \,&\partial_{\overline\lambda}\left[\mathcal M_\lambda\Phi(x, \lambda)\right]=\mathcal M_\lambda\left[ \partial_{\overline\lambda}\Phi(x, \lambda)\right]=\mathcal M_\lambda\left[  s_c(\lambda,x_3)\Phi(x, \overline\lambda)\right]\\
=\,&\Phi(x, \overline\lambda)\left[ {-\partial_{x_3}}+\tau(\lambda)\right] s_c(\lambda,{x_3})+s_c(\lambda,{x_3})\left[ {\mathcal M_\lambda}-\tau(\lambda)\right] \Phi(x, \overline\lambda)\\
=\,&\Phi(x, \overline\lambda)\left[-\partial_{x_3}+\tau(\lambda)\right] s_c(\lambda,{x_3})+ s_c(\lambda,{x_3})\left[{  \mathcal M_{ \overline \lambda}-\tau( \overline \lambda)}\right]\Phi(x, \overline\lambda)\\
=\, &\Phi(x, \overline\lambda)\left[-\partial_{x_3}+\tau(\lambda)-\tau( \overline \lambda)\right] s_c(\lambda,{x_3})+s_c(\lambda,{x_3}) \mathcal M_{\overline\lambda}\Phi(x,\overline\lambda)\\
=&\, s_c(\lambda,{x_3}) \mathcal M_{\overline\lambda}\Phi(x,\overline\lambda),
\end{align*} which yields 
\be\label{E:inverse-map-T-tilde-p}
\partial_{\overline\lambda}\left( \mathfrak Mm \right)(x,\lambda)=s_c(\lambda)e^{(\overline\lambda-\lambda)x_1+(\overline\lambda^2-\lambda^2)x_2 +(\overline\lambda^3-\lambda^3)x_3 }\left( \mathfrak Mm \right)(x, \overline\lambda).
\ee

On the other hand, for $j=1,\cdots,N$, 
\begin{align}
&- \mathcal M_{\kappa_j}(\kappa_j^N \Phi(x,\kappa_j^+)) =\sum_{n=N+1}^M\mathcal M_{\kappa_j}(  {\mathcal D}_{nj}(x_3)\Phi(x,\kappa_n^+))\nonumber\\
 =& \sum_{n=N+1}^M( \Phi(x,\kappa_n^+)[{ -\partial_{x_3}}+\tau(\kappa_j)]  {\mathcal D}_{nj}(x_3) +  {\mathcal D}_{nj}(x_3)[{  \mathcal M_{\kappa_j}}-\tau(\kappa_j)]\Phi(x,\kappa_n^+))\nonumber\\
 =& \sum_{n=N+1}^M( \Phi(x,\kappa_n^+)[-\partial_{x_3}+\tau(\kappa_j)]  {\mathcal D}_{nj}(x_3) +  {\mathcal D}_{nj}(x_3)[{ \mathcal M_{\kappa_n}-\tau(\kappa_n)}] \Phi(x,\kappa_n^+))\nonumber\\
 =& \sum_{n=N+1}^M  \Phi(x,\kappa_n^+)[{ -\partial_{x_3}-\tau(\kappa_n)+\tau(\kappa_j)}]  {\mathcal D}_{nj}(x_3)+ {\mathcal D}_{nj}(x_3) \mathcal M_{\kappa_n} \Phi(x,\kappa_n^+)\nonumber\\
 = &\sum_{n=N+1}^M  {\mathcal D}_{nj}(x_3) \mathcal M_{\kappa_n} \Phi(x,\kappa_n^+). \nonumber
\end{align} That is,
\be\label{E:soliton-debar-p}
(e^{\kappa_1x_1+\kappa_1^2x_2+\kappa_1^3x_3}\mathfrak Mm(x,\kappa _1^+) ,\cdots,e^{\kappa_Mx_1+\kappa_M^2x_2+\kappa_M^3x_3}\mathfrak Mm(x,\kappa _2^+)){ \mathcal D  }=0. 
\ee  

As $ |\lambda|\to\infty$, letting
\begin{align}
m(  x,\lambda)\sim&\sum_{j=0}^\infty\frac{M_j(  x )}{\lambda^{j}} , \, \textit{} \label{E:soliton-debar-0-p-new} \\
\mathfrak Mm\sim& q_2(x)\lambda^2+q_1(x)\lambda+q_0(x)+\frac {q_{0,\res}}\lambda+\cdots ,\label{E:soliton-debar-0-p}
 \end{align} from the Lax equation \eqref{E:intro-Lax-u},  
\begin{gather}
\qquad    2\partial_{x_1}M_{j+1}
=      (\partial_{x_2}-\partial_{x_1}^2-u)M_j,\nonumber\\
\qquad M_0=  1,\ 
M_1=  - \frac {1}2 \partial_{x_1}^{-1}u,\  
M_2=   -\frac 14 \partial_{x_2} \partial_{x_1}^{-1}u+ \frac 14 u+\frac 14 \partial_{x_1}^{-1}\left(u \partial_{x_1}^{-1} u\right),\cdots\label{E:v-evolution-KPII-new-p} 
\end{gather}  

As a result, as $\lambda\to\infty$,  
\begin{align}
&\mathfrak Mm\label{E:v-evolution-asm-s-p}\\
\to&\frac34 u_{x_1}+\frac 34 \partial_{x_1}^{-1}u_{x_2} 
+3\lambda\partial_{x_1}^2(1+\frac{M_1}{\lambda}) +3\lambda^2\partial_{x_1}(1+\frac{M_1}{\lambda}+\frac{M_2}{\lambda^2})+\frac 32u\lambda \nonumber\\
=&\frac34 u_{x_1}+\frac 34 \partial_{x_1}^{-1}u_{x_2}+\left(
-\frac32u_{x_1}+3\partial_{x_1} (-\frac 14\partial_{x_2}\partial_{x_1}^{-2}u+\frac{u}{4}+\frac14\partial_{x_1}^{-1}[u\partial_{x_1}^{-1}u])
\right) \nonumber\\
+&\lambda\left( 3\partial_{x_1}M_1+\frac{3}{2}u\right)+\frac32u(-\frac 12\partial_{x_1}^{-1}u)\nonumber\\
=&0.\nonumber\end{align}
Therefore, $q_2=q_1=q_0\equiv 0$ and $\mathfrak M\in W$. Together with \eqref{E:inverse-map-T-tilde-p}, \eqref{E:soliton-debar-p}, and the {unique solvability of the system of the CIE and $\mathcal D$ symmetry},  yields $\mathfrak M m(x,\lambda)=0$ and  the Lax pair. Using the $d$-admissible condition and taking cross differentiation, we derive the KPII equation \eqref{E:intro-KP-ist}.
\end{proof}

We conclude the paper by solving the Cauchy problem with a uniqueness result for perturbed $  {\mathrm{Gr}(N, M)_{> 0}}$ KP solitons. It yields Theorem \ref{T:intro-stability}  
\begin{theorem} \label{T:u-WP} Given $u_0 $  and $ \{z_n \}$ satisfying \eqref{E:intro-ini-data}, then $u=\mathcal S^{-1}\circ\mathcal S(u_0,\{z_n\})$  satisfies, $u:\RR\times\RR\times\RR^+\to \RR$ and \eqref{E:theorem-1.1}, \eqref{E:theorem-1.1-1}. 
Moreover,  suppose  $\det(\frac 1{\kappa_j-z_n^+})_{1\le j, n\le N} \ne 0$, $z^+_1= 0,\{z^+_n,\kappa_j\}$  distinct real, and $\det(\frac 1{\kappa_j-z_n^-})_{1\le j, n\le N} \ne 0$,   $ z^-_1=0,\{z^-_n,\kappa_j\}$  distinct real. Then
\be\label{E:u-uniq}
\mathcal S^{-1}\circ\mathcal S( u_0,\{z^+_n\} )=\mathcal S^{-1}\circ\mathcal S( u_0,\{z^-_n \}).
\ee
\end{theorem}
\begin{proof} If $u_0 $  and $ \{z_n \}$ satisfy \eqref{E:intro-ini-data} then, from \eqref{E:ini-admissible} (or elaborate argument in \cite{Wu20,Wu21}), one concludes that  $\mathcal S(u_0,\{z_n\})=(\{z_n\},\{\kappa_j\}, \mathcal D,s _c)$ is $d$-admissible. 
Therefore the Cauchy problem \eqref{E:theorem-1.1}, and the $L^\infty$-stability \eqref{E:theorem-1.1-1} are derived by applying  Theorem \ref{T:phi-k-N}, \ref{T:u-inverse}, and \ref{T:inverse-KP}.

It remains to prove \eqref{E:u-uniq}, i.e., $u^+(x)\equiv u^-(x)$ denoted as 
\be\label{E:u-plus-minus}
u^\pm=\mathcal S^{-1}\circ\mathcal S( u_0,\{z^\pm _n\} )=\mathcal S^{-1}( \{z^\pm_n\},\{\kappa_j\}, \mathcal D^\pm,s^\pm _c \}),
\ee  
where
\be\label{E:lax-inverse-WP-pm}
\begin{split}
&\left(-\partial_{x_2}+\partial_{x_1}^2+2 \lambda\partial_{x_1}+{u^\pm (x)}\right){  m}^\pm(x ,\lambda)=0 ; 
\end{split}\ee
\begin{gather}
m^\pm=1+\sum_{n=1}^N\frac{   m^\pm_{z_n,  \res }(x  )}{\lambda -z_n }  +\mathcal C  T^\pm m^\pm,\label{E:m-pm}\\
(e^{\kappa_1x_1+\kappa_1^2x_2+\kappa_1^3x_3}m^\pm(x,\kappa^+_1),\cdots,e^{\kappa_Mx_1+\kappa_M^2x_2+\kappa_M^3x_3}m^\pm(x,\kappa^+_M))\mathcal D^\pm=0;\label{E:w-w-0-pm} \\
\nonumber 
\end{gather}
\begin{align}
 T^\pm \phi  (x,\lambda)
= &s^\pm_c(\lambda  )e^{(\overline\lambda-\lambda)x_1+(\overline\lambda^2-\lambda^2)x_2+(\overline\lambda^3-\lambda^3)x_3  }\phi(x, \overline\lambda),\label{E:cauchy-operator-pm} \\
  s_c^\pm(\lambda) =& \frac { \Pi_{2\le n\le N}(\overline\lambda-z_n^\pm)}{(  \overline\lambda-z_1^\pm)^{N-1}   }\frac {\sgn(\lambda_I)}{2\pi i} \iint e^{-[(\overline\lambda-\lambda)x_1+(\overline\lambda^2-\lambda^2)x_2 ] } \label{E:cauchy-operator-pm-s-c} \\
  \times&\widetilde\xi(x_1,x_2,0, \overline\lambda) v_0(x_1,x_2)m^\pm(x_1,x_2,0, \lambda)dx_1dx_2,  
\nonumber\\
   \mathcal D ^\pm= & \widetilde{\mathcal D}^\pm \times \left({\tiny\begin{array}{ccc}\widetilde{\mathcal D}^\pm_{11}&\cdots &\widetilde{\mathcal D}^\pm_{1N}\\\vdots&\cdots &\vdots\\
\widetilde{\mathcal D}^\pm_{N1}&\cdots &\widetilde{\mathcal D}^\pm_{NN}\end{array}}\right)^{-1}\textit{diag}\,(\kappa_1^N,\cdots,\kappa_N^N) ,\label{E:BPP-new-pm}\\
 \widetilde{ \mathcal D}^\pm=&\textit{diag}\,(\frac{\Pi_{2\le n\le N}(\kappa_1-z_n^\pm)}{(\kappa_1-z_1^\pm)^{N-1}}, \cdots, \frac{\Pi_{2\le n\le N}(\kappa_M-z_n^\pm)}{(\kappa_M-z_1^\pm)^{N-1}})\mathcal D^{\sharp} \label{E:BPP-tilde-pm}\\
 { \mathcal D}^{\sharp } \ =&\left({\mathcal  D}_{ji}^\sharp\right)= \left(\mathcal  D^\flat_{ji}+\sum_{l=j}^M\frac{c_{jl}\mathcal  D^\flat_{li}}{ 1-c_{jj}} \right), \label{E:intro-sym-N-D-flat-new-pm-sharp}\\
{ \mathcal D}^{\flat } \ =&\textit{diag}\,(   
\kappa^N_1 ,\cdots,\kappa^N_M )\, A^T,\label{E:intro-sym-N-D-flat-new-pm}
\end{align} and   $c_{jl}=-\int\Psi _j(x_1,x_2,0 ) v_0(x_1,x_2)\varphi_l(x_1,x_2,0  )dx_1dx_2 $, $\Psi_j(x)$, $\varphi_l(x)$ are residue of the adjoint eigenfunction at $\kappa_j$ and values of the Sato eigenfunction at $\kappa_l$ \cite[Theorem 4]{Wu21}.

In view of \eqref{E:lax-inverse-WP-pm}, to prove   $u^+(x)\equiv u^-(x)$, it suffices to establish
\be\label{E:m-pm-uniq}
\begin{split} 
m^+(x,\lambda)
=
\Lambda m^-(x,\lambda),\qquad \Lambda=\Lambda(\lambda)=\frac{\Pi_{2\le n\le N}(\lambda-z^-_n)}{\Pi_{2\le n\le N}(\lambda-z^+_n)}.
\end{split} 
\ee Applying \eqref{E:m-pm}, \eqref{E:w-w-0-pm}, and Theorem \ref{T:phi-k-N}, it yields to showing 
\begin{gather}
\partial_{\overline\lambda}\left(\Lambda m^-\right)=    T^+\left(\Lambda m^-\right),\label{E:m-pm-0}\\
(e^{\kappa_1x_1+\kappa_1^2x_2+\kappa_1^3x_3}(\Lambda m^-)(x,\kappa^+_1),\cdots,e^{\kappa_Mx_1+\kappa_M^2x_2+\kappa_M^3x_3}(\Lambda m^-)(x,\kappa^+_M)) \mathcal D^+=0.\label{E:w-w-0-pm-0} 
\end{gather}  

Applying \cite[Theorem 4]{Wu21},  one has
\be\label{E:s-pm-0}
s_c^+ =\Lambda \overline{\Lambda}^{-1} s_c^-.
\ee
Along with \eqref{E:m-pm}, \eqref{E:cauchy-operator-pm}, yields
\begin{align}
&[\partial_{\overline\lambda}\left(\Lambda m^-\right)](x,\lambda)\label{E:m-pm-0-pf}\\
=&\Lambda(\lambda)[ \partial_{\overline\lambda}m^-](x,\lambda)=\Lambda(\lambda)[ T^-m^-](x,\lambda)\nonumber\\
=&\Lambda (\lambda)s_c^-(\lambda)e^{(\overline\lambda-\lambda)x_1+(\overline\lambda^2-\lambda^2)x_2+ (\overline\lambda^3-\lambda^3)x_3}m^-(x,\overline\lambda)\nonumber\\
=&\Lambda (\lambda)s_c^-(\lambda)\overline \Lambda^{-1} e^{(\overline\lambda-\lambda)x_1+(\overline\lambda^2-\lambda^2)x_2+ (\overline\lambda^3-\lambda^3)x_3}\Lambda(\overline\lambda)m^-(x,\overline\lambda)\nonumber\\
=&\ s_c^+(\lambda)  e^{(\overline\lambda-\lambda)x_1+(\overline\lambda^2-\lambda^2)x_2+ (\overline\lambda^3-\lambda^3)x_3}\Lambda(\overline\lambda)m^-(x,\overline\lambda)\nonumber\\
=&T^+\left(\Lambda m^-\right).\nonumber
\end{align}So \eqref{E:m-pm-0} is justified.

On the other hand, from \eqref{E:BPP-new-pm} and \eqref{E:BPP-tilde-pm}, proving of \eqref{E:w-w-0-pm-0} is equivalent to showing
\be\label{E:w-w-0-pm-0-new}
(e^{\kappa_1x_1+\kappa_1^2x_2+\kappa_1^3x_3}(\Lambda m^-)(x,\kappa^+_1),\cdots,e^{\kappa_Mx_1+\kappa_M^2x_2+\kappa_M^3x_3}(\Lambda m^-)(x,\kappa^+_M))\widetilde {\mathcal D}^+=0.
\ee 

From \eqref{E:BPP-tilde-pm},
\be\label{E:tilde-pm-mathcal-D}
\widetilde{\mathcal D}^-=\textit{diag}\,(\Lambda(\kappa_1) ,\cdots,\Lambda(\kappa_M) ) \widetilde{\mathcal D}^+.
\ee 
Along with \eqref{E:BPP-new-pm} and \eqref{E:BPP-tilde-pm}, yields
\begin{align}
&(e^{\kappa_1x_1+\kappa_1^2x_2+\kappa_1^3x_3 }\Lambda(\kappa_1) m^-(x,\kappa^+_1),\cdots,e^{\kappa_Mx_1+\kappa_M^2x_2+\kappa_M^3x_3}\Lambda (\kappa_M)m^-(x,\kappa^+_M))\widetilde{\mathcal D}^+\label{E:ini-uniq}\\
=&(e^{\kappa_1x_1+\kappa_1^2x_2+\kappa_1^3x_3 } m^-(x,\kappa^+_1),\cdots,e^{\kappa_Mx_1+\kappa_M^2x_2+\kappa_M^3x_3} m^-(x,\kappa^+_M))\widetilde{\mathcal D}^-\nonumber\\
=&0.\nonumber
  \end{align}
Therefore \eqref{E:w-w-0-pm-0-new} is proved.
\end{proof}


\begin{thebibliography}{ABC999}

\bibitem{ABF83}
 Ablowitz, M., Bar Yaacov, D., Fokas, A.: On the inverse scattering transform for the Kadomtsev Petviashvili equation. \textit{Stud. Appl. Math. }69 (1983), no. 2, 135-143. 

\bibitem{AMV13} 
Alejo, M.A., Muñoz, C., Vega, L.: The Gardner equation and the L2-stability of the N-soliton solution
of the Korteweg–de Vries equation. \textit{Trans. Am. Math. Soc.} (2013) 365(1), 195–212. 

\bibitem{Ben72}
Benjamin, T. B.: The stability of solitary waves. \textit{Proc. Roy. Soc. London Ser. A} 328 (1972), 153–183. 

\bibitem{BC07}
Biondini, G., Chakravarty, S.: Elastic and inelastic line-soliton solutions of the Kadomtsev Petviashivili II equation. \textit{Mathematics and Computers in Simulation }74 (2007), no.2-3, 237-250. 

\bibitem{BK03}
Biondini, G.,  Kodama, Y.: On a family of solutions of the Kadomtsev-Petviashvili equation which also satisfy the Toda lattice hierarchy. \textit{J. Phys. A }36 (2003), no. 42, 10519-10536. 


\bibitem{BP301} 
 Boiti, M.,  Pempinelli, F.,  Pogrebkov, A. K.,  Prinari, B.: Towards an inverse scattering theory for non-decaying potentials of the heat equation. \textit{Inverse problems}  \textbf{17}  (2001),  937-957.  

\bibitem{BP302}
Boiti, M.,  Pempinelli, F.,  Pogrebkov, A. K.,  Prinari, B.: Inverse scattering theory of the heat equation for a perturbed one-soliton potential. \textit{J. Math. Phys.}{\bf 43} (2002), no. 2, 1044-1062. 

\bibitem{BP309}
Boiti, M.,  Pempinelli, F.,  Pogrebkov, A. K.,  Prinari, B.: Building an extended resolvent of the heat operator via twisting transformations. \textit{Theoretical and Mathematical Physics} {\bf 153} (2009), no. 3, 721-733. 




\bibitem{BP310}
Boiti, M.,  Pempinelli, F.,  Pogrebkov, A. K.,  Prinari, B.: The equivalence of different approaches for generating multisoliton solutions of the KPII equation. \textit{Theoretical and Mathematical Physics} {\bf 165} (2010), no. 1, 1237-1255.  

\bibitem{BP211}
Boiti, M.,  Pempinelli, F.,  Pogrebkov, A. K.:
Green's function of heat operator with pure soliton potential. \textit{ArXiv:1201.0152v1}, 1-10,[nlin.SI] 30 Dec 2011.

\bibitem{BP212} 
Boiti, M.,  Pempinelli, F.,  Pogrebkov, A. K.: 
Extended resolvent of the heat operator with a multisolution potential. \textit{Teoret. Mat. Fiz.} 172 (2012), no. 2, 181–197. \textit{Theoret. and Math. Phys.} 172 (2012), no. 2, 1037–1051.  

\bibitem{BP214}
Boiti, M.,  Pempinelli, F.,  Pogrebkov, A. K.:  IST of KPII equation for perturbed multisoliton solutions. \textit{Topology, geometry, integrable systems, and mathematical physics}, 49-73, Amer. Math. Soc. Transl. Ser. 2, 234, Amer. Math. Soc., Providence, RI, 2014. http://dx.doi.org/10.1090/trans2/234
 
 
\bibitem{Bo75}
Bona, J.: On the stability theory of solitary waves \textit{Proc. R. Soc. Lond. Ser. A} 344(1975), 363-374. 


\bibitem{Bur84} 
Burtsev, S. P.: Damping of soliton oscillations in media with a negative dispersion law. \textit{Soviet Phys. JETP} 61 (1985), no. 2, 270–274.
 
\bibitem{CK09}
 Chakravarty, S.,  Kodama, Y.: Soliton solutions of the KP equation and application to shallow water waves. \textit{Stud. Appl. Math.} 123 (2009), no. 1, 83-151.  

 

\bibitem{D91}
Dickey, L. A. : Soliton equations and Hamiltonian systems. \textit{Advanced Series in Mathematical Physics,} 12, (1991). World Scientific Publishing Co., Inc., River Edge, NJ.

\bibitem{Ga66} 
 Gakhov, F. D.: Boundary value problems, (1966) Oxford, New York, Pergamon Press; Reading, Mass., Addison-Wesley Pub. Co.

\bibitem{Gr97}  Grinevich, P. G.: Nonsingularity of the direct scattering transform for the KP II equation with a real exponentially decaying-at-infinity potential. \textit{Lett. Math. Phys. }40 (1997), no. 1, 59-73. 

\bibitem{GN88}Grinevich, P. G., Novikov, S. P.: A two-dimensional "inverse scattering problem'' for negative energies, and generalized-analytic functions. I. Energies lower than the ground state. (Russian) \textit{Funktsional. Anal. i Prilozhen. }22 (1988), no. 1, 23–33, 96; translation in \textit{Funct. Anal. Appl.} 22 (1988), no. 1, 19–27.  





\bibitem{IK05} Ionescu, A., Kenig, C.: Local and global well-posedness of periodic KP-I   equation. \textit{Preprint}   (2005).  
\bibitem{KP70}
 Kadomtsev, B. B., Petviashvili, V. I.: On the stability of solitary waves in weakly dispersive media. \textit{Sov. Phys. Dokl.} 15 (1970), 539-541.

\bibitem{KV22}
Killip, R., Visan, M.: Orbital stability of KdV multisolitons in $H^{-1}$. \textit{ Comm. Math. Phys.}389(2022), no.3, 1445-1473. 
 
\bibitem{KS21}
Klein, C., Saut, J.: Nonlinear dispersive equations—inverse scattering and PDE methods. \textit{Applied Mathematical Sciences} 209 (2021), Springer, Cham.  

\bibitem{K17}
 Kodama, Y.: KP solitons and the Grassmannians. Combinatorics and geometry of two-dimensional wave patterns. \textit{SpringerBriefs in Mathematical Physics,} 22, (2017). Springer, Singapore. 

\bibitem{K18}
Kodama, Y.: Solitons in two-dimensional shallow water. \textit{CBMS-NSF Regional Conference Series in Applied Mathematics}, 92. Society for Industrial and Applied Mathematics (SIAM), Philadelphia, PA, (2018). 


\bibitem{KW13}
Kodama, Y., Williams, L. K.: The Deodhar decomposition of the Grassmannian and the regularity of KP solitons. \textit{Adv. Math.} 244 (2013), 979-1032. 

\bibitem{Li86}
V. D. Lipovskiĭ:
The Hamiltonian structure of the Kadomtsev-Petviashvili-II equation in a class of decreasing Cauchy data. (Russian)
\textit{Funktsional. Anal. i Prilozhen. }20 (1986), no. 4, 35-45.

\bibitem{MS93}
{Maddocks, J. H., Sachs, R. L.: On the stability of KdV multi-solitons. \textit{Comm. Pure Appl. Math.} 46 (1993), no. 6, 867–901.} 

\bibitem{MM05}
Martel, Y., Merle, F.: Asymptotic stability of solitons of the subcritical gKdV equations revisited. \textit{Nonlinearity} 
18(1), 55–80 (2005) 


\bibitem{MMT02}
Martel, Y., Merle, F., Tsai, T.P.: Stability and asymptotic stability in the energy space of the sum of N
solitons for subcritical gKdV equations. \textit{Commun. Math. Phys.}  (2002) 231(2), 347–373. 

\bibitem{MV03}
Merle, F., Vega, L.: $L^2$ stability of solitons for KdV equation. \textit{Int. Math. Res. Not.} (2003), no. 13, 735–753. 

\bibitem{M15} 
 Mizumachi, T.: Stability of line solitons for the KP-II equation in $\mathbb R^2$, \textit{Mem. Amer. Math. Soc.},  238 (2015), no. 1125, vii+95 pp. 


\bibitem{M19} 
Mizumachi, T.: The phase shift of line solitons for the KP-II equation, \textit{Nonlinear dispersive partial differential equations and inverse scattering}, 83  (2019), Fields Institute Communications, Springer Science + Business Media, part of Springer Nature. 433-495.

\bibitem{MT12} 
Mizumachi, T., Tzvetkov, N.: Stability of the line soliton of the KP-II equation under periodic transverse perturbations. \textit{Math. Ann.} 352 (2012), no. 3, 659–690. 



\bibitem{MST07} 
 Molinet, L., Saut, J., Tzvetkov, N.: Global well-posedness for the KP-I equation on the background of a non-localized solution. \textit{Comm. Math. Phys.} 272 (2007), no. 3, 775–810. 

\bibitem{MST11} 
 Molinet, L., Saut, J., Tzvetkov, N.: Global well-posedness for the KP-II equation on the background of a non-localized solution. \textit{Ann. Inst. H. Poincare Anal. Non Lineaire } \textbf{28}  (2011),  no. 5, 653-676.  



\bibitem{P00} 
 Prinari, B.: On some nondecaying potentials and related Jost solutions for the heat conduction equation. \textit{Inverse Problems} 16 (2000), no. 3, 589-603.   
 
 
\bibitem{RT09} 
  Rousset, F.; Tzvetkov, N.: Transverse nonlinear instability for two-dimensional dispersive models. \textit{Ann. Inst. H. Poincaré C Anal. Non Linéaire} 26 (2009), no. 2, 477–496. 
  
  
\bibitem{RT12} 
  Rousset, F.; Tzvetkov, N.: Stability and instability of the KdV solitary wave under the KP-I flow. Comm. Math. Phys. 313 (2012), no. 1, 155–173. 
  
    

\bibitem{S81} 
 Sato, M.: Soliton equations as dynamical systems on infinite-dimensional Grassmann manifold. \textit{RIMS Kokyuroku}  439, (1981).

\bibitem{SS82}  Sato, M., Sato, Y.: 
Soliton equations as dynamical systems on infinite-dimensional Grassmann manifold. \textit{Nonlinear partial differential equations in applied science} (Tokyo, 1982), 259-271, North-Holland Math. Stud., 81, Lecture Notes Numer. Appl. Anal., 5, North-Holland, Amsterdam, 1983.

\bibitem{S89a} 
Sato, M.: The KP hierarchy and infinite-dimensional Grassmann manifolds. Theta functions--Bowdoin 1987, Part 1 (Brunswick, ME, 1987), 51-66, \textit{Proc. Sympos. Pure Math.}, 49, Part 1, Amer. Math. Soc., Providence, RI, 1989. 



\bibitem{S89b}
Sato, M.: D-modules and nonlinear integrable systems. \textit{Algebraic analysis, geometry, and number theory }(Baltimore, MD, 1988), 325-339, Johns Hopkins Univ. Press, Baltimore, MD, 1989.


\bibitem{V62}  Vekua, I. N.: \textit{Generalized analytic functions}, 1962, Pergamon Press, London-Paris-Frankfurt; Addison-Wesley Publishing Co., Inc., Reading, Mass..

\bibitem{VA04} 
 Villarroel, J., Ablowitz, M. J.: On the initial value problem for the KPII equation with data that do not decay along a line. \textit{Nonlinearity}  \textbf{17}  (2004),  no. 5, 1843-1866. 

\bibitem{Wein86} 
Weinstein, M. I.: Lyapunov Stability of Ground States of
Nonlinear Dispersive Evolution Equations. \textit{Communications on Pure and Applied Mathematics,}   Vol. XXXIX (1986) 51-68.

\bibitem{W85} 
 Wickerhauser, M. V.: Nonlinear evolutions of the heat operator. \textit{Dissertation, Yale University }  (1985).


\bibitem{W87} 
Wickerhauser, M. V.: Inverse scattering for the heat operator and evolutions in $2+1$ variables. \textit{Comm. Math. Phys.}  108  (1987),  no. 1, 67-89. 

\bibitem{Wu18} 
 Wu, D.: The direct problem for the perturbed  Kadomtsev-Petviashvili II one line solitons. 
https:// doi.org/10.48550/arXiv.1807.01420


\bibitem{Wu20} 
Wu, D.: The direct scattering problem for the perturbed $Gr(1,2)_{> 0}$  Kadomtsev-Petviashvili II solitons. \textit{Nonlinearity}  33 (2020), no. 12, 6729–6759.   

\bibitem{Wu21} 
Wu, D.: The direct scattering problem for perturbed Kadomtsev–Petviashvili multi line solitons. \textit{J. Math. Phys.}   62 (2021), no. 9, Paper No. 091513, 19 pp.  
\bibitem{Wu24} 
Wu, D.: The inverse scattering theory of Kadomtsev-Petviashvili II equations. http://doi.org/10.48550 /arXiv.2408.07868





\bibitem{Z75}
Zakharov, V. E.: Instability and nonlinear oscillations of solitons. \textit{Journal of Experimental and Theoretical Physics Letters} Vol. 22 (1975),  172-174.  
 

\end{thebibliography}
\end{document}